\documentclass[10pt]{amsart}
\usepackage{fullpage, amsfonts,amsmath,amscd,amssymb}
\usepackage{amssymb, amsbsy, amsthm, amsmath, amstext, amsopn, verbatim,
multicol}
\usepackage{calc}
\usepackage{ifthen}

\newtheorem*{thm}{Theorem}
\newtheorem*{prop}{Proposition}
\newtheorem*{lem}{Lemma}
\newtheorem*{cor}{Corollary}

\theoremstyle{definition}

\newtheorem*{question}{Question}

\theoremstyle{remark}

\newtheorem*{rem}{Remark}

\DeclareMathOperator{\res}{cont}
\DeclareMathOperator{\diag}{diag}

\newcommand{\im}{\operatorname{im}}

\newcommand{\ep}{\epsilon}

\newcommand{\h}{\mathfrak{h}}
\newcommand{\D}{\mathcal{D}}

\newcommand{\CC}{\mathbb{C}}

\newcommand{\M}{\mathcal{M}}

\newcommand{\NN}{\mathbb{N}}
\newcommand{\Z}{\mathbb{Z}}
\newcommand{\Q}{\mathbb{Q}}
\newcommand{\R}{\mathbb{R}}
\newcommand{\zl}{\mathbb{Z}_0^{\ell}}
\newcommand{\C}{\mathbb{C}}
\newcommand{\tx}{\tilde{{X}}} \newcommand{\ty}{\tilde{{Y}}}
\newcommand{\tbx}{\tilde{{\bf X}}} \newcommand{\tby}{\tilde{{\bf Y}}} \newcommand{\tv}{\tilde{v}} \newcommand{\tw}{\tilde{w}}
\DeclareMathOperator{\id}{id} \DeclareMathOperator{\gr}{gr}

 \DeclareMathOperator{\md}{--mod}
\DeclareMathOperator{\hilb}{Hilb}
\newcommand{\Hr}{{\bf H}^{\text{reg}}}
\newcommand{\hi}[1]{\hilb^{#1}\CC^2} \newcommand{\hil}[1]{\hilb^n(#1)}
 \DeclareMathOperator{\Tr}{Tr}

\DeclareMathOperator{\spc}{Spec} 
\DeclareMathOperator{\rp}{Rep} \DeclareMathOperator{\mt}{Mat}
 
 \DeclareMathOperator{\lie}{Lie}

\DeclareMathOperator{\re}{Re} \DeclareMathOperator{\ide}{Id}
\DeclareMathOperator{\ima}{Im} \DeclareMathOperator{\grad}{grad}
\DeclareMathOperator{\proj}{Proj} \DeclareMathOperator{\coh}{Coh}
\DeclareMathOperator{\supp}{supp} \DeclareMathOperator{\cont}{N}
    \newcommand{\bep}{\boldsymbol{\ep}}
    \newcommand{\btheta}{\boldsymbol{\theta}}  
    \newcommand{\bpsi}{\boldsymbol{\psi}}
    \newcommand{\blambda}{\boldsymbol{\lambda}}
    \newcommand{\btau}{\boldsymbol{\tau}}
    \newcommand{\bmu}{\boldsymbol{\mu}}
    \newcommand{\bnu}{\boldsymbol{\nu}}
    \newcommand{\bempty}{\boldsymbol{\emptyset}}
\DeclareMathOperator{\irr}{\textsf{Irr}} \DeclareMathOperator{\ch}{\bf Char}
    \DeclareMathOperator{\rch}{\bf rCh}

\begin{document}
\title{Quiver varieties, category $\mathcal{O}$ for rational Cherednik algebras, and Hecke algebras.}
\author{I.G.~Gordon} \address{School of Mathematics and Maxwell Institute for Mathematical Sciences,
University of Edinburgh,
James Clerk Maxwell Building,
Kings Buildings, Mayfield Road, Edinburgh EH9 3JZ,
U.K.}
\email{igordon@ed.ac.uk}
 \thanks{ I thank Cedric Bonnaf\'e, Ali Craw, Toshiro Kuwabara, Hiraku Nakajima, Richard Vale and Xavier Yvonne for useful comments and conversations. I am grateful to the Leverhulme Trust, the Glasgow Mathematical Journal Trust Fund, Glasgow University, Kyoto University, the University of Chicago, the Isaac Newton Institute and EPSRC for support while parts of this paper were written.}
\begin{abstract}
We relate the representations of the rational Cherednik algebras associated with the complex reflection group $\mu_{\ell} \wr \mathfrak{S}_n$ to sheaves on Nakajima quiver varieties associated with extended Dynkin graphs via a $\Z$-algebra construction. This is done so that as the parameters defining the Cherednik algebra vary, the stability conditions defining the quiver variety change. 

This construction motivates us to use the geometry of the quiver varieties to interpret the ordering function (the $c$-function) used to define a highest weight structure on category $\mathcal{O}$ of the Cherednik algebra. This interpretation provides a natural partial ordering on $\mathcal{O}$ which we expect will respect the highest weight structure. This partial ordering has appeared in a conjecture of Yvonne on the composition factors in $\mathcal{O}$ and so our results provide a small step towards a geometric picture for that. 

We also interpret geometrically another ordering function (the $a$-function) used in the study of Hecke algebras. (The connection between Cherednik algebras and Hecke algebras is provided by the KZ-functor.) This is related to a conjecture of Bonnaf\'e and Geck on equivalence classes of weight functions for Hecke algebras with unequal parameters since the classes should (and do for type $B$) correspond to the G.I.T. chambers defining the quiver varieties. As a result anything that can be defined via the quiver varieties, including the $a$-function, will be constant on these classes. 

\end{abstract}

\maketitle
\tableofcontents
\section{Introduction}
\subsection{} In this paper we point out a simple relationship between the combinatorics of certain rational Cherednik algebras and the geometry of certain Nakajima quiver varieties. We also show relations to the cell combinatorics of certain Iwahori-Hecke algebras with unequal parameters.

These connections all arise from an attempt to find a geometric model for the category $\mathcal{O}$ of rational Cherednik algebras which could be used to understand the composition factors of standard modules. There is nothing earthshaking in this, but it opens up a field of speculation about  Hecke algebras and quiver varieties that we should like to graze around in.

The results here show that rational Cherednik algebras can be degenerated to some specific Nakajima quiver varieties and that these quiver varieties still contain much combinatorial information which is relevant to category $\mathcal{O}$. Then, thanks to the KZ-functor, this transfers to combinatorial information on a corresponding cyclotomic Hecke algebra. In particular the $c$-function on $\mathcal{O}$ and the $a$-function on the Hecke algebra correspond to Morse functions on quiver varieties:  previously studied orderings built from these functions are then unified by stratifications of explicit subvarieties of the quiver varieties, and these subvarieties can be studied geometrically. However, the finer structure of $\mathcal{O}$ and the Hecke algebras are not immediately visible to the geometry;  to see that some some further rigid structure will probably be needed. 

\subsection{The cast} We concentrate on the complex reflection group $G = G_n(\ell) = \mu_{\ell} \wr \mathfrak{S}_n$ where $\ell$ and $n$ are natural numbers. This group acts naturally on its reflection representation $\h = \C^{n}$. Associated to $G$ there are two algebras depending on an $\ell$-dimensional parameter space.
\begin{enumerate}
\item Rational Cherednik algebras $H_{\bf h}$ (with parameter $t=1$) are deformations of the differential operator ring $\D (\h) \ast G$ where the deformation depends on parameters ${\bf h}= (h,H_1, \ldots, H_{\ell-1}) \in \Q^{\ell}$.
\item The Iwahori-Hecke algebras $\mathcal{H}_{\bf q}(G)$ are deformations of the group algebra $\C G$ where the deformation depends on parameters ${\bf q}  \in (\C^*)^{\ell}$.
\end{enumerate}
There is also a family of varieties depending on the same parameter space.
\begin{enumerate}
\item[(3)] Nakajima quiver varieties $\mathcal{M}_{\btheta}(n)$ generically resolve the singular space $(\h\oplus \h^*)/G$ where the family depends on stability parameters ${\btheta} \in \Q^{\ell}$.
\end{enumerate}
\subsection{Cherednik algebras} We study the full subcategory $\mathcal{O}_{\bf h}$ of $H_{\bf h}\md$ introduced in \cite{duop} and studied further in \cite{GGOR}, see \ref{O} for the definition. This category has a highest weight structure. Its simple objects are labelled by $\irr G$, and this set is in natural bijection with $\ell$-multipartitions of $n$, ${\blambda} = (\lambda^{(1)}, \ldots , \lambda^{(\ell)})$. The ordering on $\irr G$ is given according to the value of the {\it $c$-function} which assigns to $\blambda\in \irr G$ the scalar $c_{\bf h}(\blambda)$ by which the deformed Euler operator ${\bf z}\in H_{\bf h}$ acts on the highest weight of the simple object in $\mathcal{O}_{\bf h}$ corresponding to $\blambda$. This function first appeared in this context in \cite[Lemma 2.5]{duop} and has also played a role in Kazhdan-Lusztig theory and representations of finite groups of Lie type. Rouquier showed in \cite{rou} that much about the Morita equivalence classes of $\mathcal{O}_{\bf h}$ can be understood by studying the regions of the parameter space in which the values of the $c$-function induce the same ordering on multipartitions. We call these regions {\it $c$-chambers}. They are finite in number.

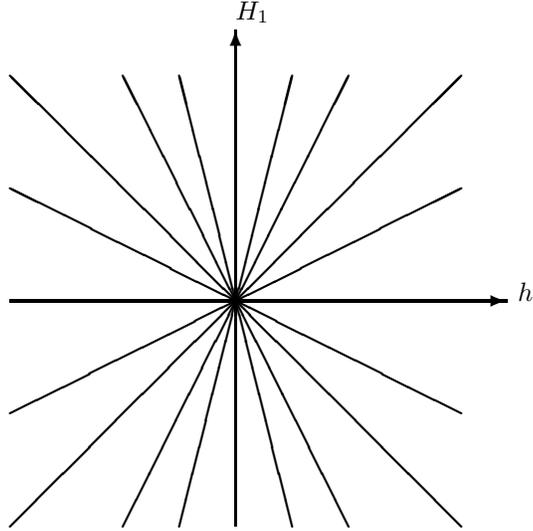
\begin{figure}[h]
\setlength{\unitlength}{0.75cm}
$$
\begin{picture}(9,9)
\thicklines
\put(0,0){\line(1,1){8}}
\put(2,0){\line(1,2){4}}
\put(3,0){\line(1,4){2}}
\put(4,0){\vector(0,1){8.8}}
\put(5,0){\line(-1,4){2}}
\put(6,0){\line(-1,2){4}}
\put(8,0){\line(-1,1){8}}
\put(0,2){\line(2,1){8}}
\put(0,4){\vector(1,0){8.8}}
\put(8,2){\line(-2,1){8}}
\put(9,4){$h$}
\put(4,9){$H_1$}
\end{picture}
$$
\caption{$c$-chambers for $\ell=2$, $n =3$. In each chamber there is a total ordering on $\ell$-multipartitions.}
\label{cch}
\end{figure}

It is our aim to study a little more of $\mathcal{O}_{\bf h}$ in terms of these $c$-chambers. In particular, in the interior of any $c$-chamber the $c$-function induces a total ordering on the set of multipartitions; we believe this ordering is too coarse for representation theory and so we introduce a partial ordering which should govern the combinatorics of $\mathcal{O}_{\bf h}$. To do this, we require geometry.   

\subsection{Quiver varieties} \label{intquiv}The quiver varieties $\mathcal{M}_{\btheta}(n)$ are G.I.T. quotients  equipped with canonical projective morphisms $\pi_{\btheta} : \mathcal{M}_{\btheta}(n) \longrightarrow (\h\times \h^*)/G$. For generic choices of stability parameter $\btheta \in \Q^{\ell}$ these provide symplectic resolutions of singularities. Results of Crawley-Boevey and of LeBruyn can be used to describe the G.I.T. chamber structure on $\Q^{\ell}$ in terms of the combinatorics of the affine root system of type $\tilde{A}_{\ell -1}$. 

To relate these varieties with Cherednik algebras we recall that when $\ell =1$ $H_{\bf h}$ provides a quantisation of the Hilbert scheme of $n$ points on the plane, the relevant quiver variety in this special case,  \cite{GS}. This quantisation is constructed by showing that the Opdam-Heckman shift functors for $H_{\bf h}$ are noncommutative analogues of powers of an ample line bundle that appears naturally in the quiver theoretic description of the Hilbert scheme. The quantisation procedure then works effectively whenever the shift functors induce equivalences of categories.  For general $\ell$ we follow an analogous procedure. We use the naive shift functors introduced in \cite{Gdiff} to construct a functor between $H_{\bf h}\md$ and $\coh \mathcal{M}_{\btheta}(n)$ for {\it any} ${\bf h}$ and ${\btheta}$. However for this functor to be useful we would like the naive shift functors to be equivalences. We make an ansatz based on the equivalences between category $\mathcal{O}_{\bf h}$'s constructed by Rouquier. This shows that there should be a simple relation between ${\bf h}$ and ${\btheta}$:
$$\btheta = ( - h - H_1 - \cdots - H_{\ell -1}, H_1, \ldots H_{\ell-1}).$$
Having made this identification, we show easily that the walls of the G.I.T. chambers are walls of the $c$-chambers.

\begin{figure}[h]
\setlength{\unitlength}{0.75cm}
$$
\begin{picture}(9,9)
\thicklines
\put(0,0){\line(1,1){8}}
\put(2,0){\line(1,2){4}}
\put(4,0){\vector(0,1){8.8}}
\put(6,0){\line(-1,2){4}}
\put(8,0){\line(-1,1){8}}
\put(0,4){\vector(1,0){8.8}}
\put(9,4){$h$}
\put(4,9){$H_1$}
\end{picture}
$$
\caption{G.I.T. chambers for $\ell=2$, $n =3$. In each chamber there is a partial ordering on $\ell$-multipartitions.}
\label{GITch}
\end{figure}
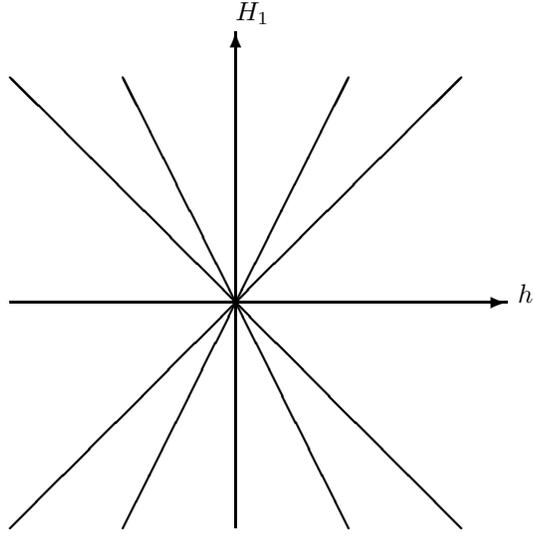

\subsection{} There is a $\C^*$-action on $\mathcal{M}_{\btheta}(n)$ lifting the usual hyperbolic action on $(\h \oplus \h^*)/G$. The attracting set of this action is $\mathcal{Z}_{\btheta} = \pi_{\btheta}^{-1}((\h\times 0)/G)$. The Bialynicki-Birula decomposition then provides a partial ordering on the irreducible components of $\mathcal{Z}_{\btheta}$. In case $\mathcal{M}_{\btheta}(n)$ is smooth we show that there is a {\it natural} labelling of these components by the $\ell$-multipartitions of $n$. This is achieved by using the hyper-K\"ahler structure on $\mathcal{M}_{\btheta}(n)$: rotation of the complex structure provides an equivariant identification with a {\it generalised Calogero-Moser space}, $\mathcal{X}_{\btheta}(n)$, which is a moduli space for representations of rational Cherednik algebras wih parameter $t=0$. The fixed points on these affine varieties can be represented by {\it baby Verma modules}. These modules occur in families for all $\bf h$ and are labelled by $\ell$-multipartitions of $n$. We write $x_{\btheta}(\blambda)$ for the  fixed point of $\mathcal{M}_{\btheta}(n)$ corresponding to $\blambda$. The irreducible component corresponding to $\blambda$ is then the closure of the set of points in $\mathcal{Z}_{\btheta}$ attracted to $x_{\btheta}(\blambda)$ under the $\C^*$-action. We have therefore a {\it geometric ordering} on $\ell$-multipartitions of $n$ which depends on the parameter $\btheta$.

\subsection{}This description gives us more. The Bialynicki-Birula decomposition is an algebraic analogue of a Morse theoretic decomposition. The hyper-K\"ahler structure on $\mathcal{M}_{\btheta}(n)$ involves three real symplectic forms, two of which are used to make the complex symplectic form mentioned above. Taking the moment map for the $U(1) < \C^*$-action with respect to the third form produces a Morse function $$f_{\btheta}: \mathcal{M}_{\btheta}(n) \longrightarrow \lie(U(1))^* = \R$$ and the ordering induced by the values of this function at the critical (i.e. fixed) points refines the geometric partial ordering. Rotating the complex structure then lets us relate this function with a moment map for the complex symplectic form on $\mathcal{X}_{\btheta}(n)$, and in turn with the representation theory of the rational Cherednik algebra, giving our first theorem.
\begin{thm}
\label{firstthm}
Let $\btheta = ( - h - H_1 - \cdots - H_{\ell -1}, H_1, \ldots H_{\ell-1})$ be in the interior of a G.I.T. chamber. Then $c_{\bf h}(\blambda) = f_{\btheta}(x_{\btheta}(\blambda))$. In particular the geometric partial ordering on $\ell$-multipartitions of $n$ is refined by the $c$-ordering.
\end{thm}

\subsection{} \label{onthewall} We can extend the result on the partial ordering to all choices of $\btheta$, although we do not have a topological interpretation of the $c$-function. In this case there are fewer fixed points and we find a non-trivial partition of the set of $\ell$-multipartitions by naturally associating a fixed point to each multipartition. 

\subsection{} In order to determine explicitly the geometric partial ordering we follow the lead of Haiman, \cite{hai}. We can reduce to the case where $h=-1$ by a combination of rescaling $\bf h$ by positive rationals and applying a simple duality which swaps $h$ and $-h$. The parameter space has now essentially $\ell-1$ dimensions. There is a natural action of the affine Weyl group $\tilde{\mathfrak{S}}_{\ell}$ on this space which makes the walls of the G.I.T. chambers into a subset of the reflecting hyperplanes and so the G.I.T. chambers are unions of alcoves. As $n$ tends to infinity the chamber decomposition converges to the alcove decomposition. In general there will be some regions which are unions of infinitely many alcoves. One of these is the {\it asymptotic region} -- it gives rise to the dominance ordering on $\ell$-multipartitions.

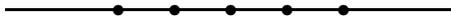
\begin{figure}[h]

\setlength{\unitlength}{0.75cm}
$$
\begin{picture}(8,1)
\put(0,0){\line(1,0){8}}
\put(2,0){\circle*{0.2}}
\put(3,0){\circle*{0.2}}
\put(4,0){\circle*{0.2}}
\put(5,0){\circle*{0.2}}
\put(6,0){\circle*{0.2}}
\end{picture}
$$
\caption{The line through $h=-1$: alcoves appearing for $\ell=2$, $n =3$}
\label{ach}
\end{figure}

The action of $\tilde{\mathfrak{S}}_{\ell}$ on the parameter space is mirrored geometrically by Nakajima's reflection functors which allow us to study subvarieties of Hilbert schemes instead of $\mathcal{M}_{\btheta}(n)$. As explained by Haiman, \cite{hai}, the combinatorics of these varieties is easier to understand thanks to Nakajima's geometric construction of representations of Kac-Moody Lie algebras. It leads us to our second theorem. 
\begin{thm}
\label{secondthm}
The geometric partial ordering on $\ell$-multipartitions is induced from the dominance ordering on partitions under a version of the classical bijection between $\ell$-multipartitions and partitions with a given core. The choice of core depends on $\btheta$.
\end{thm}
This theorem can also be extended to deal with values of $\btheta$ on the walls of G.I.T. chambers.

\subsection{}
Putting the two theorems above together shows that there is a simply described, geometrically defined partial ordering on $\ell$-multipartitions which is refined by the $c$-ordering. We hope that this partial ordering is the real ordering for the highest weight structure on $\mathcal{O}_{\bf h}$. Because of its geometric nature we are able to pose a straightforward question at the end of the paper on the nature of the characteristic cycles in $\mathcal{M}_{\btheta}(n)$ of standard modules in $\mathcal{O}_{\bf h}$ which would imply our hope. 

This ordering has also appeared in the work of Yvonne, \cite{yvo}, on the decomposition matrices of cyclotomic $q$-Schur algebras or more generally of category $\mathcal{O}_{\bf h}$'s. His conjectures  depend on the work of Uglov on higher level Fock spaces, \cite{uglov}. They predict an ordering on $\mathcal{O}_{\bf h}$ which depends on the choice of a multi-charge depending on $\bf h$; this multi-charge corresponds to our choice of core corresponding to $\btheta$. The geometric ordering is then the ordering Yvonne predicts with, in the notation of his paper, $\xi = 1$. (We expect the ordering he has for other choice of $\xi$ to be refined the ordering we give here. In other words, the geometric ordering is the mother of all orderings, just as the dominance ordering on partitions is the mother of all orderings for the symmetric group.) Yvonne studies the case of a  ``dominant" multi-charge in most detail: that corresponds to the asymptotic case here.

\subsection{Hecke algebras} 
The combinatorics appearing here seems to be related to recent work on cells in Hecke algebras with unequal parameters. This should not be too surprising since the categories $\mathcal{O}_{\bf h}$ are expected to play the role of generalised $q$-Schur algebras for cyclotomic Hecke algebras with the KZ-functor playing the role of the Schur functor. (In fact in the asymptotic parameter case Rouquier has shown in \cite{rou} that $\mathcal{O}_{\bf h}$ is equivalent to the module category of an appropriate cyclotomic $q$-Schur algebra associated to $G$). 

Jacon has studied an $a$-function associated to each irreducible representation of $G$ and used the ordering induced by the values it takes, called the $a$-ordering, to label the irreducible representations of $\mathcal{H}_{\bf q}(G)$, \cite{Jac}. We can describe this function in terms of a Morse function associated to half of the hyperbolic $\C^*$-action introduced above, giving us the following result.
\begin{thm}
\label{thirdthm}
Jacon's $a$-function is linear on the G.I.T. chambers. Moreover, the geometric partial ordering on $\ell$-multipartitions of $n$ is refined by the $a$-ordering.
\end{thm}
The $a$-ordering and the $c$-ordering are at first sight unrelated: in general they induce different total orderings of $\ell$-multipartitions. Unlike the $c$-ordering, the $a$-ordering is not linear on the whole parameter space. However we show that the description in terms of half the $\C^*$-action produces an extension of a theorem of Brou\'e and Michel, \cite{brmi}, which presents the $c$-function as a sum of two related $a$-functions.

\subsection{} If $\ell =2$ then we have that $G$ is the Weyl group of type $B_n$. There is then a partition of $\irr G$ into two-sided cells and an ordering on these cells which depends on the choice of parameter $\bf q$, \cite{lus1}. Geck has conjectured that the parameter space for $\mathcal{H}_{\bf q}(G)$ should decompose into a finite number of chambers and that the cell structure should then depend only on the chamber or wall in which $\bf q$ lies, \cite{geck}. These conjectural chambers are precisely the G.I.T. chambers for $\ell =2$. Moreover, in all known examples the partition into two-sided cells agrees with the partition by fixed points mentioned in \ref{onthewall}, and the partial ordering on two-sided cells agrees with the geometric partial ordering. We suspect this may be a consequence of the properties of KZ-functor, \cite{GGOR}, combined with conjectural geometric ordering on $\mathcal{O}_{\bf h}$. 

\subsection{Calogero-Moser spaces and cells} We note in passing that our constructions here allow us to give a combinatorial description of the blocks of restricted rational Cherednik algebras. Details can be found in \cite{GM} where a conjectural link to two-sided cells for finite Weyl groups is also made.


\subsection{The case $G = \mathfrak{S}_n$}
For $\ell = 1$ the results of this paper reduce to the well-studied case of $G = \mathfrak{S}_n$. Here the parameter space is one-dimensional. There is a wall (i.e. a point) at the origin, and it is known that the shift functors which send $h$ to $h+1$ encounter a problem when they cross that wall. Rouquier's theorem agrees with this.  The simple modules of $\mathcal{O}_{h}$ are labelled by partitions. For negative values of $h$, the natural ordering on partitions is the dominance ordering; for positive values things switch around and it is the anti-dominance ordering. Geometrically this corresponds on one side to the quiver variety being the Hilbert scheme with its usual tautological bundle, and on the other side the Hilbert scheme with the dual bundle. The $a$-function is $n(\lambda)$, the $c$-function is, up to the addition of a constant, $h(n({}^t\lambda) - n(\lambda))$: both of these functions respect the dominance ordering. 
 
\subsection{Organisation} In Sections 2 and 3 we recall some definitions and theorems concerning rational Cherednik algebras and quiver varieties respectively. In Section 4 we compare the $c$-chambers and the G.I.T. chambers, whilst in Section 5 we label the $\C^*$-fixed points and then show that the $c$-function can be interpreted topologically. We recall some combinatorial notions associated to partitions in Section 6. In Section 7 we calculate the geometric ordering explicitly and in Section 8 we extend our results to the walls of the chambers. We relate the geometric combinatorics to Hecke algebras in Section 9. In Section 10 we pose a number of questions concerning characteristic cycles, derived categories, the generalised $n!$-conjecture and $q$-Schur algebras. There are then a couple of appendices which deal with a pair of gruesome calculations that I didn't know how to do properly. 

\subsection{Notation} Throughout $n$ will denote a fixed positive integer and $\ell$ will be a positive integer. Given two families $\{x_{\lambda}\}_{\lambda \in \Lambda}$ and $\{y_{\lambda}\}_{\lambda\in \Lambda}$ of real numbers indexed by the same set $\Lambda$, we will write $$x_\lambda \doteq y_\lambda \qquad \text{for all }\lambda\in\Lambda$$ to indicate that $x_{\lambda} = y_{\lambda} + C$  for all $\lambda$ where the constant $C$ is independent of $\lambda\in \Lambda$.

\section{The $c$-function and rational Cherednik algebras}

\subsection{}
Let $\mu_{\ell}$ be the cyclic subgroup of $SL_2(\C)$ generated by $\sigma = \diag (\eta, \eta^{-1})$ where $\eta = \exp (2\pi \sqrt{-1}/ \ell)$. The vector space $V = (\C^2)^n$ admits an action of \label{Gnl-defn}$G = G_n(\ell) = \mathfrak{S}_n \ltimes (\mu_{\ell})^n$: $(\mu_{\ell})^n$ acts by extending the natural action of $\mu_{\ell}$ on $\C^2$, whilst $\mathfrak{S}_n$ acts by permuting the $n$ copies of $\C^2$. For an element $\gamma\in \mu_{\ell}$ and an integer $1\leq i\leq n$ we write $\gamma_i$ to indicate the element $(1,\ldots, \gamma, \ldots ,1) \in \mu_{\ell}^n$ which is non--trivial only in the $i$--th factor.

\subsection{Partitions and multipartitions}
\label{irreds}
A partition of degree $n$ is a non--increasing sequence (finite or infinite) $\lambda = ( \lambda_1 \geq \lambda_2 \geq \ldots)$ of non--negative integers with sum $n$. We write $|\lambda | = n$. We identify two partitions that differ only by zeroes. We denote by $\mathcal{P}(n)$\label{part-defn} the set of all partitions of $n$. We denote the dominance ordering on $\mathcal{P}(n)$ by \label{dom-defn}$\unlhd$, so that $\mu \unlhd \lambda$ precisely when $\mu_1 + \cdots + \mu_i \leq \lambda_1 + \cdots + \lambda_i$ for all $i$.

Given a partition $\lambda$ we let ${}^t\lambda$ \label{trpart-defn}be the transposed partition and note that $\mu\unlhd \lambda$ if and only if ${}^t\lambda \unlhd {}^t\mu$. Set $n(\lambda) = \sum_{i} \lambda_i (i-1)$, the partition statistic. 

An $\ell$-multipartition of degree $n$ is an $\ell$--tuple of partitions $\blambda  = (\lambda^{(1)}, \ldots , \lambda^{(\ell)})$ with $\sum |\lambda^{(i)}| = n.$ We let $\mathcal{P}(\ell,n)$ \label{multpart-defn}denote this set of multipartitions. There is also a dominance ordering on $\mathcal{P}(\ell,n)$ where $\bmu \unlhd \blambda$ if and only if $\sum_{k=1}^{j-1} |\mu^{(k)}| + (\mu_1^{(j)} + \cdots + \mu_i^{(j)}) \leq  \sum_{k=1}^{j-1} |\lambda^{(k)}| + (\lambda_1^{(j)} + \cdots + \lambda_i^{(j)})$ for all $i$ and for all $1\leq j \leq \ell$.

\subsection{Irreducible representations of $G$} The set of isomorphism classes of complex irreducible representations of $G_n({\ell})$ are labelled by $\mathcal{P}(\ell ,n)$. We follow the natural labelling presented in \cite[6.1.1]{rou}: here the trivial representation corresponds to $((n), \emptyset, \ldots, \emptyset)\in \mathcal{P}(\ell ,n)$.

\subsection{Parameter space} Throughout the paper we will be using the set of rational parameters, $\bf H$, \label{H-defn}which consists of $\ell$-tuples ${\bf h} = (h, H_1, \ldots , H_{\ell -1})\in \Q^{\ell}$\label{h-defn}. We define $H_0\in \Q^{\ell}$ by $H_0 + \cdots + H_{\ell -1} = 0$.

\subsection{The $c$--function}  \label{cdef}
The following function $c: \mathcal{P}(\ell , n) \times \bf{H} \rightarrow \mathbb{Q}$ will be the central throughout:
$${c}_{\bf h}({\blambda}) = \ell \sum_{r=2}^{\ell} |\lambda^{(r)}|(H_1+ \cdots + H_{r-1}) - \ell \left( \frac{n(n-1)}{2} + \sum_{r=1}^{\ell} n(\lambda^{(r)}) - n({}^t\lambda^{(r)}) \right) h.$$
Given ${\bf h}\in {\bf H}$, the $c$-function induces an ordering on $\mathcal{P}(\ell ,n)$ by the rule:
$${\blambda} >_{\bf h} {\bmu} \Leftrightarrow c_{\bf h}({\bmu}) > c_{\bf h}({\blambda}).$$\label{cord-defn}
We call this the $c$-order. The dependence of this order on the parameters decomposes $\mathbf{H}$ into a finite number of chambers defined by the linear equations $c_{\bf h}({\blambda}) = c_{\bf h}(\bmu)$ for all ${\blambda}, {\bmu} \in \mathcal{P}(\ell, n)$. We call these the $c$-chambers, see Figure \ref{cch}. In the interior of a $c$-chamber the $c$-order is a total order on $\mathcal{P}(\ell ,n)$.

\subsection{Rational Cherednik algebras, \cite{EG}} There is a symplectic form on $V$ which is induced from $n$ copies of the  standard symplectic form $\omega$ on $\C^2$. If we pick a basis $\{ x, y\}$ for $\C^2$ such that $\omega (x,y) = 1$ then we can extend this naturally to a basis $\{ x_i , y_i : 1\leq i\leq n\}$ of $V$ such that the $x$'s and the $y$'s form Lagrangian subspaces and $\omega (x_i, y_j) = \delta_{ij}$. We let $T(V^*)$ denote the tensor algebra on $V^*$: with our choice of basis this is just the free algebra on generators $X_i, Y_i$ for $1\leq i\leq n$ where $X_i$ and $Y_i$ are the dual basis to $x_i$ and $y_i$. The rational Cherednik algebra $H_{t, {\bf h}}$ \label{RCA-defn}associated to $G$ is the quotient of the smash product $T(V^*) \ast G$ by the following relations:
\begin{eqnarray*}
&X_iX_j = X_jX_i, \qquad Y_iY_j = Y_jY_i \qquad & \text{for all } 1\leq i,j\leq n \\
&[Y_i, X_i ] = t + h\sum_{j\neq i} \sum_{t=0}^{\ell -1} s_{ij}\sigma_i^t \sigma_j^{-t} + \sum_{t=0}^{\ell-1} \big( \sum_{j=0}^{\ell -1} \eta^{-tj} H_j \big) \sigma_i^t & \text{for } 1\leq i \leq n \\
& [Y_i,X_j] = -h \sum_{t=0}^{\ell -1} \eta^t s_{ij}\sigma_i^t \sigma_j^{-t}. & \text{for } i\neq j
\end{eqnarray*}

\subsection{} \label{goingbetweenparams} There is another presentation of $H_{t, {\bf h}}$ given in terms of conjugation invariant functions on the complex reflections in $G$, \cite{EG} and \cite{Gdiff}. To get to that presentation send $h$ to $k$ and set $c_{\sigma^t} = \sum_{j=0}^{\ell -1} \eta^{-tj}H_j$.

\subsection{Category $\mathcal{O}_{\bf h}$} \label{O} Let $\mathcal{O}_{\bf h}$ be the category of finitely generated $H_{1, {\bf h}}$-modules on which all the variables $Y_i$ act locally nilpotently. This is a highest weight category with simple objects $\{ L_{\bf h}(\blambda) : \blambda \in \mathcal{P}(\ell, n)\}$ and ordering given by $L_{\bf h}(\blambda) < L_{\bf h}(\bmu)$ if $\blambda <_{\bf h} \bmu$ and $c_{\bf h}({\blambda}) - c_{\bf h}({\bmu}) \in \Z$.

The following was proved in \cite[Theorem 5.5]{rou}.
\label{rouq}
\begin{thm} Let $\mathbf{h}, \mathbf{h'} \in \mathbf{H}$ belong to the same $c$-chamber and differ by an element of $\mathbb{Z}^{\ell}$. Then $\mathcal{O}_{\mathbf{h}}$ and $\mathcal{O}_{\mathbf{h'}}$ are equivalent.
\end{thm}

\subsection{The $t=0$ case} There is a significant difference between $H_{t,{\bf h}}$ with $t\neq 0$ and $H_{0,{\bf h}}$: when $t\neq 0$ the centre is trivial; when $t=0$ the algebra is module-finite over its centre. We will let $Z_{0,{\bf h}}$ denote the centre of $H_{0,{\bf h}}$\label{Z-defn}. By \cite[p.267]{EG} $Z_{0, {\bf h}}$ is a Poisson algebra, its bracket arising from its quantisation by a subalgebra of $H_{t,{\bf h}}$

\section{Quiver varieties}

\subsection{} Let $Q$ \label{Q-defn}be the cyclic quiver with $\ell$ vertices and cyclic orientation. Let $\overline{Q}$ be the double quiver of $Q$, obtained by inserting an arrow $a^{\ast}$ in the opposite direction to every arrow $a$ in the quiver. By the McKay correspondence the vertices $0, \ldots , \ell -1$ are in a sensible 1-1 correspondence with the irreducible representations of $\mu_{\ell}$. Vertex $i$ corresponds to the representation $\eta_i$ which sends $\sigma\in \mu_{\ell}$ to $\eta^i$. In this way we can label representations of $\mu_{\ell}$ by positive roots of (the affine Lie algebra associated to) $Q$, i.e. of type $\tilde{A}_{\ell -1}$: for $0\leq i \leq \ell -1$ the simple root $\alpha_i$ corresponds to vertex $i$ and hence to the irreducible representation $\eta_i$; the fundamental root $\delta = \sum_{i=0}^{\ell -1} \alpha_i$ corresponds to the regular representation of $\mu_{\ell}.$

\subsection{} Choose an extending vertex of $Q$: in this case it could be any vertex; we take it to be $0$. Then let $Q_{\infty}$ be the quiver obtained by adding one vertex named $\infty$ to $Q$ that is joined to $0$ by a single arrow and let $\overline{Q}_{\infty}$ denote the double quiver of $Q_{\infty}$.

We will consider representation spaces of the quiver \label{Qi-defn}$\overline{Q}_{\infty}$. Let ${\bf d} = (d_0, \ldots, d_{\ell-1}) \in \Z_{\geq 0}^{\ell}$ be a dimension vector for $\overline{Q}$ and set ${\bf d'} = e_{\infty} + {\bf d}$, a dimension vector for $\overline{Q}_{\infty}$.
Recall that \begin{eqnarray*} \rp(\overline{Q}_{\infty}, {\bf d'}) &=& \left(\bigoplus_{i=0}^{\ell -1}  \mt_{d_{i+1},d_i}(\C)\right) \oplus  \left(\bigoplus_{i=0}^{\ell -1} \mt_{d_i,d_{i+1}}(\C)\right)
\oplus \C^{d_0} \oplus (\C^{d_0})^{\ast} \\ & = &\{ (X_0,\ldots , X_{\ell-1}, Y_0, \ldots ,Y_{\ell -1}; v, w) \} \\ &=& \{ ({\bf X},{\bf Y}; v, w)\}.\end{eqnarray*}
We set \label{rep-defn}$R({\bf d'}) =  \rp(\overline{Q}_{\infty}, {\bf d'})$. Let ${\bf G}({\bf d}) = \prod_{i=0}^{\ell -1} GL_{d_i}(\C)$ be the base change group. If ${\bf g} = (g_0,\ldots , g_{\ell-1})\in {\bf G}({\bf d})$ then ${\bf g}$ acts on $R({\bf d'})$ by
$${\bf g}\cdot (X_0, \ldots, X_{\ell-1},Y_0, \ldots ,Y_{\ell -1},
; v, w) = (g_1X_0g_0^{-1}, \ldots , g_0X_{\ell-1}g_{\ell -1}^{-1}, g_0Y_0g_1^{-1}, \ldots , g_{\ell -1}Y_{\ell -1}g_0^{-1}; g_0v, wg_0^{-1} ).$$ Orbits of ${\bf G}({\bf d})$ on $R({\bf d'})$ are in $1$--$1$ correspondence with isomorphism classes of representations of $\overline{Q}_{\infty}$ of dimension ${\bf d'}$.

\subsection{Hyper-K\"ahler structure, \cite[Section 2]{Nak}}
There is a quaternionic structure on the complex space $R({\bf d'})$ given by letting the quaternion $J$ act by $$J ({\bf X}, {\bf Y}; v,w) = ({\bf Y}^{\dag}, - {\bf X}^{\dag}; w^{\dag}, -v^{\dag})$$ where the daggers denote the Hermitian adjoint. There is also an inner product on $R({\bf d'})$ given by \begin{equation} \label{metric}g(({\bf X},{\bf Y};v,w), (\tbx, \tby; \tv, \tw)) = \re \left[ \sum_{r=0}^{\ell -1} \left( \Tr (X_r \tx_r^{\dag}) + \Tr (Y_r \ty_r^{\dag})\right) + \Tr(v \tv^{\dag}) + \Tr (w \tw^{\dag}) \right].\end{equation} Thus $R({\bf d'})$ has a hyper-K\"ahler structure. We have three associated real symplectic forms on $R({\bf d'})$ $$\omega_I (\cdot , \cdot) = g(I\cdot, \cdot), \quad \omega_J (\cdot , \cdot) = g(J\cdot, \cdot), \quad \omega_K (\cdot , \cdot) = g(K\cdot, \cdot ),$$ which we split into the {real} symplectic form $\omega_{\R} = \omega_I$ and the {complex} symplectic form $\omega_{\C} = \omega_J + \sqrt{-1} \omega_K$. The subgroup ${\bf U}({\bf d}) \leq {\bf G}({\bf d})$ acts on $R({\bf d'})$ preserving the forms.

\subsection{} Associated to the action of ${\bf U}({\bf d})$ on each of the symplectic forms we have moment maps $\mu_{\R} = \mu_I : R({\bf d'}) \longrightarrow (\lie {\bf U}({\bf d}) )^*$ and $\mu_{\C} = \mu_J + \sqrt{-1}\mu_K : R({\bf d'}) \longrightarrow (\lie {\bf G}({\bf d}))^*$. Using the trace pairing, we can identify $\lie {\bf U}({\bf d})$ and $\lie {\bf G}({\bf d})$ with their duals and hence write the maps explicitly as
$$ \mu_{\R} ({\bf X},{\bf Y};v,w) = \frac{\sqrt{-1}}{2} \left( [{\bf X},{\bf X}^{\dag}] + [{\bf Y}, {\bf Y}^{\dag}] + vv^{\dag} - w^{\dag}w \right)$$
and
$$ \mu_{\C} ({\bf X},{\bf Y};v,w) = [{\bf X},{\bf Y}] + vw.$$
The first map is ${\bf U}({\bf d})$-equivariant, the second ${\bf G}({\bf d})$-equivariant.

\subsection{Quotient varieties} Given $\btheta = (\theta_0, \ldots , \theta_{\ell -1})\in \Q^{\ell}$ we can introduce two complex varieties.

The first is the algebro-geometric quotient \label{X-defn}$$\mathcal{X}_{\btheta}({\bf d}) = \mu_{\C}^{-1}({\btheta})//{\bf G}({\bf d})$$ where we abuse notation by letting ${\btheta}$ also denote $(\theta_0 \ide_{d_0}, \ldots, \theta_{\ell -1} \ide_{d_{\ell-1}})\in \lie {\bf G}({\bf d})$. It is, by \cite[Theorem 3.2.2]{GG}, an affine variety; its points parametrise the isomorphism classes of semisimple representations of dimension ${\bf d'}$ of deformed preprojective algebras $\Pi^{\btheta}(Q_{\infty})$, \cite{CB}.

The second, a {\it quiver variety}, is the geometric invariant theory quotient \label{M-defn}$$\mathcal{M}_{\btheta}({\bf d}) = \mu_{\C}^{-1}(0)//_{\btheta} {\bf G}({\bf d}).$$ It is, by \cite[Section 2]{king}, a variety projective over its base $\mathcal{X}_{0}({\bf d})$; its points parametrise the isomorphism classes of $\theta$-polystable representations of dimension ${\bf d'}$ of the preprojective algebra, \cite[Proposition 3.2]{king} and \cite{CB}.

A more algebraic description of $\mathcal{M}_{\btheta}({\bf d})$ is given by \begin{equation}\label{semi-inv} \mathcal{M}_{\btheta}({\bf d}) = \proj \bigoplus_{i\geq 0}\C[\mu_\C^{-1}(0)]^{\chi_{\btheta}^i}.\end{equation} Here ${\chi_{\btheta}}: {\bf G}({\bf d}) \longrightarrow \C^*$ is the (fractional) character which sends ${\bf g} = (g_0, \ldots , g_{\ell -1})$ to $\prod (\det g_i)^{\theta_i}$ and $\C[\mu_\C^{-1}(0)]^{\chi_{\btheta}^i}$ denotes the space of semi-invariant functions on $\mu_\C^{-1}(0)$, i.e. those which transform as $${\bf g}\cdot f = \chi_{\btheta}^i({\bf g}) f.$$ By definition the space $\C[\mu_\C^{-1}(0)]^{\chi_{\btheta}^i}$ is zero if $i\btheta \notin \Z^{\ell}$.
\subsection{}
\label{KempfNess}
These varieties can be described in terms of the hyper-K\"ahler structure. There is a homeomorphism between $\mathcal{M}_{\btheta}({\bf d})$ and  $\mu^{-1}_{\C}(0)\cap \mu^{-1}_{\R}(\frac{\sqrt{-1}}{2}{\btheta}) / {\bf U}({\bf d})$ and also between $\mathcal{X}_{\btheta}({\bf d})$ and $\mu^{-1}_{\C}({\btheta})\cap \mu^{-1}_{\R}(0) / {\bf U}({\bf d})$, see \cite[Section 6]{king} and \cite[Theorem 3.24]{nakhilb}. This shows that both $\mathcal{M}_{\btheta}({\bf d})$ and $\mathcal{X}_{\btheta}({\bf d})$ are hyper-K\"ahler reductions of $R({\bf d'})$ and hence themselves (possibly singular) hyper-K\"ahler manifolds with real and complex forms $\omega_{\R}$ and $\omega_{\C}$ induced from those on $R({\bf d'})$. When the manifolds are smooth these forms are symplectic; in general they give rise to Poisson structures.

\subsection{} \label{rotate} The hyper-K\"ahler description of \ref{KempfNess} allows us to compare $\mathcal{M}_{\btheta}({\bf d})$ and $\mathcal{X}_{\btheta}({\bf d})$. There is a mapping by ``rotating the complex structure", given by multiplication by $(-I-K)/\sqrt{2}$:
\begin{equation} \label{rotation}\Psi: \mathcal{M}_{\btheta}({\bf d}) = \mu^{-1}_{\C}(0)\cap \mu^{-1}_{\R}(\frac{\sqrt{-1}}{2}\btheta) / {\bf U}({\bf d}) \longrightarrow \mu^{-1}_{\C}(\frac{1}{2}{\btheta}) \cap \mu^{-1}_{\R}(0) / {\bf U}({\bf d}) = \mathcal{X}_{\frac{1}{2}\btheta}({\bf d}).\end{equation} It is a diffeomorphism in the smooth case.

The following lemma will be useful later on.
\begin{lem}
The real form $\omega_{\R}$ on $\mathcal{M}_{\btheta}({\bf d})$ is sent to the imaginary part of the complex form $\omega_{\C}$ on $\mathcal{X}_{\frac{1}{2}\btheta}({\bf d})$ under the diffeomorphism $\Psi$.
\end{lem}
\begin{proof}
This is a simple general fact. Let $z_1, z_2\in R({\bf d'})$ and set $u = (-I-K)/\sqrt{2}$. Then
\begin{eqnarray*}
(\Psi \omega_{\R}) (z_1,z_2) = \omega_{\R}(\Psi^{-1} z_1, \Psi^{-1} z_2) & = & g(Iu^{-1} z_1, u^{-1} z_2) \\ & = &g( uIu^{-1}z_1, z_2) \\ & = &  g(Kz_1, z_2) \\ & = & \omega_K (z_1,z_2) = (\ima \omega_{\C})(z_1,z_2), \end{eqnarray*} where we used the quaternionic invariance of the form in the third equality.
\end{proof}

\subsection{$\C^*$-action} \label{circleaction}There is a $\C^*$-action on both $\mathcal{M}_{\btheta}({\bf d})$ and $\mathcal{X}_{\btheta}({\bf d})$, induced by the following hyperbolic action on $R({\bf d'})$ $$\lambda \circ ( {\bf X},{\bf Y};v,w) =  (\lambda {\bf X}, \lambda^{-1}{\bf Y} ; v,w).$$ This restricts to a $U(1)$-action on $ \mu^{-1}_{\C}(0)\cap \mu^{-1}_{\R}(\frac{\sqrt{-1}}{2}{\btheta}) / {\bf U}({\bf d})$ and $\mu^{-1}_{\C}({\btheta})\cap \mu^{-1}_{\R}(0) / {\bf U}({\bf d})$. The mapping $\Psi$ of \eqref{rotation} is $U(1)$-equivariant. Indeed letting $\lambda \in U(1)$ we see that \begin{eqnarray*} \lambda \circ \Psi(({\bf X},{\bf Y}; v,w)) \!\!\!&=& \!\!\!\frac{1}{\sqrt{2}} (\lambda \sqrt{-1} ({\bf X}-{\bf Y}^*) , \lambda^{-1} \sqrt{-1}({\bf Y}+{\bf X}^*); \sqrt{-1}(v-w^*), \sqrt{-1}(w-v^*)) \\ &=&  \frac{1}{\sqrt{2}}(\sqrt{-1} (\lambda {\bf X} - (\lambda^{-1}{\bf Y})^*), \sqrt{-1} (\lambda^{-1}{\bf Y} + (\lambda {\bf X})^*); \sqrt{-1}(v-w^*), \sqrt{-1}(w-v^*))\\ & =& \Psi( \lambda \circ ({\bf X},{\bf Y};v,w)).\end{eqnarray*} In particular we deduce that the $\C^*$--fixed points on $\M_{\btheta}({\bf d})$ correspond naturally to the $\C^*$--fixed points on $\mathcal{X}_{\frac{1}{2}\btheta}({\bf d})$.

Observe that the metric on $R({\bf d'})$ is $U(1)$-stable and hence so too are the symplectic forms $\omega_I, \omega_J$ and $\omega_K$. Moreover the complex symplectic form $\omega_{\C}$ is $\C^*$-equivariant.

\subsection{Resolutions} \label{resol}We specialise to the case ${\bf d} = n\delta$ where $\delta = (1,\ldots ,1)$ is the affine dimension vector of $Q$. In this case we simplify our notation, writing \label{X-2defn}$\mathcal{X}_{\btheta}(n)$ and \label{M-2defn}$\mathcal{M}_{\btheta}(n)$ for $\mathcal{X}_{\btheta}({\bf d})$ and $\mathcal{M}_{\btheta}({\bf d})$ respectively. Thanks to \cite[Theorem 1.1]{CBdecomp} combined with \cite[Lemma 9.2]{CB} we have $\mathcal{X}_{0}(n) \cong V/G$ as a Poisson variety and so $\mathcal{M}_{\btheta}(n)$ is projective over $V/G$. When $\mathcal{M}_{\btheta}(n)$ is smooth this gives a symplectic resolution \label{pi-defn}$$\pi_{\btheta}: \mathcal{M}_{\btheta}(n) \longrightarrow V/G.$$ (To see that $\mathcal{M}_{\btheta}(n)$ is  connected apply the mapping $\Psi$ from \eqref{rotation} and then observe that $\mathcal{X}_{
\frac{1}{2}\btheta}(n)$ is connected by Theorem \ref{mo's}.)

\subsection{First relation to rational Cherednik algebras} \label{mo's} The varieties $\mathcal{X}_{\btheta}(n)$ appear in the study of rational Cherednik algebras with $t=0$.
\begin{thm}[\cite{EG}, \cite{mo}] Let ${\bf h}\in {\bf H}$ and set $\btheta = (-h+H_0, H_1, \ldots , H_{\ell -1})$. Let $Z_{0,{\bf h}}$ denote the centre of the rational Cherednik algebra $H_{0,{\bf h}}$. Then there is a $\C^*$-equivariant isomorphism between the complex Poisson varieties $$\spc Z_{0, {\bf h}} \longrightarrow \mathcal{X}_{\btheta}(n).$$
\end{thm}
\begin{proof} The result is stated in \cite[Proposition 6.6 and Theorem 7.4]{mo} for a different labelling of parameters. It is an elementary calcuation to go between the parameters in \cite[Section 7]{mo} (where his $c_1$ is $2h$) and the parameters here using \ref{goingbetweenparams} and \cite[6.2]{mo}: this shows that there is an isomorphism between $\spc Z_{0, {\bf h}}$ and $\mathcal{X}_{\ell\btheta}(n)$ which is a Poisson mapping up to a scalar multiple. In the proof of \cite[Theorem 11.16]{EG} this scalar is shown to be $1/\ell$. Rescaling from $\mathcal{X}_{\ell\btheta}(n)$ down to $\mathcal{X}_{\btheta}(n)$ rescales the Poisson bracket by $1/\ell$ and thus provides the Poisson isomorphism. The $\C^*$-equivariance is evident from the construction of the isomorphism.
\end{proof}

\section{Chamber decompositions} \label{chpict}

\subsection{Second relation to rational Cherednik algebras} \label{Zalg} We would like to relate the representation theory of $H_{1,{\bf h}}$ with the geometry of the spaces $\mathcal{M}_{\btheta}(n)$. To do this we combine the $\Z$-algebra formalism of \cite{GS} with the differential operator approach of \cite{Gdiff}. We will use the definitions and notation of $\Z$-algebras from \cite{GS} without further comment: the interested reader should consult that paper for details, in particular Section 5 there.
\begin{thm} Let ${\bf h}\in {\bf H}$ and set $\btheta = ( -h+H_0,  H_1, \ldots , H_{\ell -1})$.
Then there is a noncommutative filtered $\Z$-algebra $B_{\bf h}$ such that the following properties hold.
\begin{enumerate}
\item[(i)] There is functor from $H_{1,{\bf h}}\md$ to $\coh B_{\bf h}$ which preserves filtrations.
\item[(ii)] (Vale, \cite{vale}) For generic ${\bf h}$ this is an equivalence.
\item[(iii)] The $\Z$-algebra $\gr B_{\bf h}$ is isomorphic to the $\Z$-algebra associated to the homogeneous coordinate ring of $\mathcal{M}_{\btheta}(n)$ introduced in \eqref{semi-inv}.
\end{enumerate}
\end{thm}
\begin{proof}
We need a little notation before beginning. Let
$$ e= \frac{1}{\ell^n n!} \sum_{g\in G} g\in \C G$$ be the symmetrising idempotent. The subalgebra $eH_{1, {\bf h}}e$ is denoted by $U_{\bf h}$ and called the {\it spherical algebra}.

The $\Z$-algebra $B_{\bf h}$ is constructed from by gathering ``shift functors" between module categories for spherical subalgebras associated to various parameters ${\bf h'}$. Here we use the functors defined in \cite[Lemma 4.4]{Gdiff}: to any fractional character $\Lambda$ of ${\bf G}(n\delta)$ we associate a filtered $(U_{\bf h'}, U_{\bf h''})$-bimodule $B_{{\bf h'}, {\bf h''}}^{\Lambda}$. The statement of \cite[Lemma 4.4]{Gdiff} is given in terms of the parameter space $(k,c)$ mentioned in \ref{goingbetweenparams}, but a simple calculation shows that ${\bf h'}$ and ${\bf h''}$ are related by the rule $$h' = h''-\Lambda_0  -\cdots - \Lambda_{\ell -1}, \quad\text{and}\quad H_i' = H_i'' + \Lambda_i \quad\text{for all }1\leq i \leq \ell -1.$$ The shift functor $U_{{\bf h''}}\md \longrightarrow U_{{\bf h'}}\md$ is then given by tensoring  by $B_{{\bf h'}, {\bf h''}}^{\Lambda}\otimes_{U_{\bf h''}} - $. By \cite[Lemma 4.1]{Gdiff} the associated graded module is \begin{equation}\label{grsemi-inv}\gr B_{{\bf h'}, {\bf h''}}^{\Lambda} = \C[\mu_{\C}^{-1}(0)]^{\Lambda}.\end{equation}

Now given {\it any} $\btheta' \in \Q^{\ell}$ we could construct the $\Z$-algebra $B_{\bf h}(\btheta')$ as the following direct sum $$B_{\bf h}(\btheta') = \bigoplus_{i\geq j\geq 0} B^{(i-j)\btheta'}_{{\bf h}_i, {\bf h}_j}$$ where $${\bf h}_i = {\bf h} + i(-\theta_0' - \cdots -\theta_{\ell -1}', \theta_1', \ldots , \theta'_{\ell -1}) \quad \text{and} \quad {\bf h}_j = {\bf h} + j(-\theta_0' - \cdots -\theta_{\ell -1}',  \theta'_1, \ldots , \theta'_{\ell -1}).$$ Then formally following the arguments of \cite[Sections 5.4 and 5.5]{GS} and using \eqref{grsemi-inv} would yield an isomorphism between $\gr B_{\bf h}(\btheta')$ and the $\Z$-algebra associated to the homogeneous coordinate ring of $\mathcal{M}_{\btheta}(n)$, as well as a functor from $U_{\bf h}\md$ to $\coh B_{\bf h}(\btheta')$. Composing this functor with the idempotent functor $M\mapsto eM$ from $H_{1, {\bf h}}\md$ to $U_{\bf h}\md$ completes a proof of (i) and (iii).

Using this general $\Z$-algebra would be unsatisfactory, however, as we could not expect that the functor from $H_{1, {\bf h}}\md$ to $\coh B_{\bf h}(\btheta')$ would be an equivalence. For this we make an ansatz which explains our choice of stability condition. By analogy with \cite[3.16]{GS} we expect that the shift functor $B^{\Lambda}_{\bf h', h''}\otimes_{U_{\bf h''}}-$ sends $\mathcal{O}_{\bf h''}$ to $\mathcal{O}_{\bf h'}$. For the functor above to be an equivalence we would in particular expect the restriction of the shift functor to category $\mathcal{O}$ to be an equivalence. However, Rouquier's Theorem \ref{rouq} shows that we cannot pick ${\bf h'}$ and ${\bf h''}$ independently; they should lie in the same chamber. The most obvious way to ensure this is to take all ${\bf h'}$ to live on the positive part of the line from the origin through ${\bf h}$. For this to happen we pick $\btheta = (-h+H_0, H_1, \ldots , H_{\ell -1})$ which gives ${\bf h}_i = (i+1) {\bf h}$. We take $B_{\bf h}$ to be $B_{\bf h}(\btheta)$.

The generic equivalence of (ii) is proved in \cite{vale}.
\end{proof}
\subsection{} It is an important problem to calculate for which ${\bf h'}$ the shift functor $B_{\bf h', \bf h''}^{\Lambda}\otimes_{U_{\bf h''}} -$ is an equivalence. For $\ell =1$ this is answered in \cite[Section 3]{GS} for another construction of the shift functors; however, the main result of \cite{GGS} shows that the functors of \cite{GS} agree with the definition given here.

\subsection{G.I.T. chambers} \label{GITcham} For the rest of the paper we will take ${\bf h}\in\bf H$ and enforce the relation
 $\btheta = ( -h+H_0 , H_1, \ldots , H_{\ell -1})$. The following lemma is closely related to \cite[Theorem 2.8]{Nak}.
\begin{lem} The variety $\mathcal{M}_{\btheta}(n)$ is smooth if $\btheta$ does not lie on one of the following hyperplanes
\begin{equation} \label{smooth} (H_i+\cdots +H_j) + mh = 0 \qquad \text{or } \qquad h=0,\end{equation} where $1\leq i \leq j \leq \ell -1$ and $1-n\leq m \leq n-1$.
\end{lem}
\begin{proof}
By \cite[Theorem 1.2]{leb} $\mathcal{X}_{\frac{1}{2}\btheta}(n)$ is smooth if and only if all representations in the parameter space are simple. Applying the rotation $\Psi^{-1}$ of \eqref{rotation} matches the simple representations to the stable representations, \cite[Section 7]{BlB}. Then by either \cite[Proposition 8.8]{Lb@2} or the argument of 
\cite[Lemma 3.10(2)]{Nak2} if every representation in the space $\mathcal{M}_{\btheta}(n)$ is stable implies then 
$\mathcal{M}_{\btheta}(n)$ is a smooth variety. Thus the smoothness of $\mathcal{X}_{\frac{1}{2}\btheta}(n)$ implies the smoothness of $\mathcal{M}_{\btheta}(n)$.

We now describe the precise condition that ensures all the representations of $\Pi^{\frac{1}{2}\btheta}(Q_{\infty})$ are simple, and hence by \cite[Theorem 1.2]{leb} the precise condition that ensures that $\mathcal{X}_{\frac{1}{2}\btheta}(n)$ is smooth. Let $R$ denote the root system of type $A_{\ell -1}$ corresponding to the subquiver of $Q$ with vertices $1, \ldots , \ell -1$. Then all representations of $\mathcal{X}_{\frac{1}{2}\btheta}(n)$ are simple if and only if \begin{equation} \label{simples} \btheta \cdot (\beta + m\delta) \neq 0\end{equation} for some $\beta \in R\cup \{ 0\}$, $0\leq m \leq n-1$ such that $\beta + m\delta$ is positive.

To prove this we must show that there are no simple representations of the deformed preprojective algebra $\Pi^{{\btheta}}(Q_{\infty})$ of dimension $\gamma < {\bf d'} = e_{\infty} + n\delta$ if and only if \eqref{simples} is satisfied. If such a representation exists then obviously $\gamma$ either involves $e_{\infty}$ or it doesn't, and so by factoring out the representation corresponding to $\gamma$ if necessary, we can assume without loss of generality that $\gamma$ does not involve $e_{\infty}$ and so is supported entirely on $Q$. Now we apply \cite[Theorem 1.2]{CB} to see that $\btheta \cdot \gamma = 0$ and that $\gamma$ must be a positive root of the root system associated to $Q$. In particular since $\gamma < n\delta$ we have $\gamma = \beta + m\delta$ as required. Conversely, suppose that $\btheta \cdot (\beta + m\delta) = 0$ for some $\beta$ and some $m\delta$. Since $\btheta \cdot {\bf d'} = 0$ this means that a decomposition of ${\bf d'} = \gamma_1 + \cdots + \gamma_r$ into roots of $Q_{\infty}$ which are minimal with respect to the condition that ${\btheta}\cdot \gamma_i= 0$ is non-trivial, i.e. $r>1$. By \cite[Theorem 1.2]{CB} $\Pi^{\frac{1}{2}{\btheta}}(Q_{\infty})$ has a semisimple representation whose components have dimension $\gamma_1, \cdots , \gamma_r$, as required.

Since $R$ is a root system of type $A_{\ell -1}$, $\beta$ is either zero or has the form $\pm (\alpha_i + \cdots +\alpha_j)$ for some $1\leq i \leq j \leq \ell-1$. Thus \eqref{simples} becomes either $m h = 0$ in the case $\beta = 0$ or $$\mp (H_i + \cdots +H_j) + m h =0.$$ This defines the set of hyperplanes in \eqref{smooth}.
\end{proof}
We expect the condition in the Lemma is actually necessary too. For $\btheta$ such that $\btheta \cdot \delta \neq 0$ arguments similar to \cite[Proposition 8.10]{Lb@2} should show that $\mathcal{M}_{\btheta}(n)$ is singular, but if $\btheta \cdot \delta = 0$ then there are no stable representations of dimension vector $e_{\infty} + n\delta$ and we don't know whether the results of {\it loc. cit.} apply. Of course, if the identification of \ref{KempfNess} were a diffeomorphism then the last two paragraphs of the proof would prove the necessity.

\subsection{} We let \label{Hr-defn}$\Hr$ denote the open subset of ${\bf H}$ obtained by removing the hyperplanes occuring in \eqref{smooth}. The above lemma shows us that the hyperplanes occuring in \eqref{smooth} contain the hyperplanes which define the G.I.T. chambers. That is, inside any one of the chambers defined by these hyperplanes the corresponding varieties are isomorphic and have the same associated tautological bundle. This is clear since the proof of the lemma shows that the set of $\btheta$-stable representations and $\btheta'$-stable representations of $\mu_{\C}^{-1}(0)$ are exactly the same if $\btheta$ and $\btheta'$ belong to the same chamber. Of course, it is possible that there are too many hyperplanes specified in \eqref{smooth}, but the condition in the lemma does turn out to be necessary too then they will be exactly the walls of the G.I.T. chambers. By abusing language we will now call the hyperplanes in \eqref{smooth} the G.I.T. walls and the chambers they define the G.I.T. chambers. See Figure \ref{GITch} in \ref{intquiv} for a picture.

\subsection{G.I.T. chambers versus $c-$chambers} \label{RouqvsGIT} Now we can relate the above criterion to the $c$-chamber decomposition. The following lemma states that the set of G.I.T. walls is a subset of the set of $c$-walls and so, in particular, the $c$-chambers are contained inside the G.I.T. chambers.

\begin{thm} The $c$-chamber decomposition of ${\bf H}$ refines the G.I.T. chamber decomposition. 
\end{thm}
\begin{proof}
This is very straightforward. For each hyperplane described in \eqref{smooth} we must find a corresponding pair $\blambda , \bmu \in \mathcal{P}(\ell , n)$ such that $c_{\bf h}({\blambda}) = c_{\bf h}({\bmu})$ defines the hyperplane.

We begin by realising all G.I.T. walls with non--negative $h$--coefficient. Consider  $$\blambda(a,b,j):=( (a), \emptyset, \emptyset, \ldots, (b), \emptyset , \ldots , \emptyset)\in \mathcal{P}(\ell ,n)$$ where $a,b$ are non-negative integers with $a+b = n$ and the partition $(b)$ appears in the $j$-th entry for some $2\leq j \leq \ell $. Then \begin{eqnarray*}c_{\bf h}(\blambda(a,b,j)) &=& \ell \left[ b(H_1 + \cdots + H_{j-1}) - \left(\frac{n(n-1)}{2} + 0 - \frac{a(a-1)}{2} + 0 - \frac{b(b-1)}{2}\right)h\right] \\ & = & \ell b {[} (H_1 + \cdots + H_{j-1}) - ah{]}. \end{eqnarray*} Therefore, letting $j$ vary between $2$ and $\ell$ and $a$ vary between $0$ and $n-1$ yields all the hyperplanes in \eqref{smooth} which feature $H_1$ and have $h$ appearing with a non--negative coefficient.

Now we observe that starting with $\blambda = (\lambda^{(1)} , \ldots, \lambda^{(\ell)})\in \mathcal{P}(\ell ,n)$ if we consider its conjugate ${}^{t}\blambda = ({}^t\lambda^{(0)} , \ldots, {}^t\lambda^{(\ell-1)})$ then the only difference between $c_{\bf h}({\blambda})$ and $c_{\bf h}({}^t\blambda)$ is in the coefficient of $h$. Call this coefficient $h_{\blambda}$ and $h_{{}^t\blambda}$ respectively. From the definition we have $$h_{{}^t\blambda} + \ell\frac{n(n-1)}{2} = - \left(h_{\blambda} + \ell\frac{n(n-1)}{2}\right)$$ so that $h_{{}^t\blambda} = -h_{\blambda} - \ell n(n-1)$. It follows that if $\bmu$ is another element of $\mathcal{P}(\ell ,n)$ then $$h_{{}^t\blambda} - h_{{}^t\bmu} = - (h_{\blambda} - h_{\bmu}).$$ Thus given any $c$-wall we can find another $c$-wall whose $h$--coefficient has been multiplied by $-1$. This together with the previous paragraph yields all hyperplaces in \eqref{smooth} featuring $H_1$.

Now we need to find the G.I.T. walls that don't include $H_1$. To do so we induct on $\ell$, using the notation $c_{\bf h}^{(\ell)}(\blambda)$ to describe the $c$-functions dependence on $\ell$. The induction begins because the previous paragraphs have dealt with $H_1$. Now suppose $\ell >1$ and consider $\lambda := (\emptyset, \lambda^{(2)}, \ldots , \lambda^{(\ell)})\in \mathcal{P}(\ell , n)$ and $\overline{\lambda} := ( \lambda^{(2)}, \ldots , \lambda^{(\ell)})\in \mathcal{P}(\ell -1 , n).$ Then \begin{eqnarray*}c_{\bf h}^{(\ell)}(\blambda) &=& \ell (\sum_{i=2}^{\ell} |\lambda^{(i)}|(H_1 + \cdots + H_{i-1}) - h_{\lambda}h) \\ & = & \ell (\sum_{i=3}^{\ell} |\lambda^{(i)}| (H_2 + \cdots +H_{i-1}) - h_{\lambda}h) +  \ell nH_1 \\ & = & c_{\bf h'}^{(\ell -1)}(\overline{\blambda}) + \ell nH_1,\end{eqnarray*} where ${\bf h'} = (h, H_2, \ldots H_{\ell -1})$. Thus for any $\blambda , \bmu$ of the above form we find $$c_{\bf h}^{(\ell)}({\blambda}) - c_{\bf h}^{(\ell)}(\bmu) = c_{\bf h'}^{(\ell-1)}({\overline{\blambda}})- c_{\bf h'}^{(\ell-1)}({\overline{\bmu}}).$$ Thus, by the induction hypothesis, we find all combinations that don't involve $H_1$. This completes the proof.
\end{proof}

\section{The $c$-function, topologically}\label{morch}
In this section we will show that the $c$-function can be interpreted as the value of a Morse function on $\mathcal{M}_{\btheta}(n)$ at the $U(1)$-fixed points. We use this to give the $c$-ordering a geometric significance. Recall that throughout we will let $\btheta$ and $\bf h$ be related by formula of Theorem \ref{Zalg}: $$\btheta = (-h+H_0,H_1, \ldots , H_{\ell -1}).$$
\subsection{Fixed points} \label{fix} Recall the $\C^*$-action on $\mathcal{M}_{\btheta}(n)$ from \ref{circleaction}.
\begin{lem} Let ${\bf h}\in \Hr$. Then the $\C^*$-fixed points on $\mathcal{M}_{\btheta}(n)$ are naturally labelled by the $\ell$-multipartitions of $n$.
\end{lem}
\begin{proof} Since $\btheta$ belongs to the interior of a G.I.T. chamber, the varieties $\mathcal{M}_{\btheta}(n)$ and $\mathcal{X}_{\frac{1}{2}\btheta}(n)$ are smooth by Lemma \ref{GITcham}. By \ref{circleaction} the $\C^*$-fixed points of $\mathcal{M}_{\btheta}(n)$ correspond under the mapping $\Psi$ of \eqref{rotation} to the $\C^*$-fixed points of $\mathcal{X}_{\frac{1}{2}\btheta}(n)$. Moreover, by Theorem \ref{mo's}, there is a $\C^*$-equivariant isomorphism between $\mathcal{X}_{\frac{1}{2}\btheta}(n)$ and $\spc Z_{0,{\bf h}/2}$. Thus it is enough to describe the $\C^*$-fixed points of $\spc Z_{0, {\bf h}/2}$.

By \cite[Section 3.6]{babyV} there is a $\C^*$--equivariant morphism \begin{equation} \label{upsilon} \Upsilon: \spc Z_{0,{\bf h}/2} \longrightarrow \h/G \times \h^*/G.\end{equation} Here the $\C^*$-action on the codomain is induced from the action $\lambda \circ (z_1, z_2) = (\lambda z_1, \lambda^{-1} z_2)$ on $\h\times \h^*$. Thus any fixed point of $\spc Z_{0,{\bf h}/2}$ must be mapped by $\Upsilon$ to a fixed point of $\h/G \times \h^*/G$. However, the only fixed point of $\h/G \times \h^*/G$ is the origin. Thus the fixed points of $\spc Z_{0, {\bf h}/2}$ belong to the fibre $\Upsilon^{-1}(0)$.

The fibre $\Upsilon^{-1}(0)$ is described in \cite[Section 5]{babyV}. In particular, its (closed) points are labelled by the isomorphism classes of simple $G$--modules, i.e. by elements of $\mathcal{P}(\ell ,n)$. Morever, since the fibre is finite and $\C^*$ is connected, each point is fixed by the $\C^*$-action.
\end{proof}

\subsection{} \label{nat} Let us make explicit what the word ``natural" means in the statement of Lemma \ref{fix}.
\begin{cor}
Let $\blambda \in \mathcal{P}(\ell ,n)$ and let \label{xt-defn}$x_{\btheta}(\blambda) \in \mathcal{M}_{\btheta}(n)$ be the $\C^*$-fixed point corresponding to $\blambda$ which is constructed in Lemma \ref{fix}. Then the assignment $\btheta \mapsto x_{\btheta}(\blambda)$ extends to a continuous section of the family of quiver varieties $\{ \mathcal{M}_{\btheta}(n): \btheta\in \Q^{\ell}\}$ over the parameter space $\Q^{\ell}$. In other words, the fixed points of the varieties $\mathcal{M}_{\btheta}(n)$ vary continuously in families and extend to the walls of the G.I.T. chambers.
\end{cor}
\begin{proof}
When $\btheta$ belongs to the interior of a G.I.T. chamber then the fixed point of $\spc Z_{0, {\bf h}/2}$ labelled by $\blambda$ corresponds to the baby Verma module $M_{{\bf h}/2}(\blambda)$ constructed in \cite[Section 4.2]{babyV}. However, the construction in \cite[Section 4.2]{babyV} exists for all choices of ${\bf h'} \in \C^{\ell}$ and shows that the modules $\{ M_{\bf h'}(\blambda) : {\bf h'}\in \C^{\ell} \}$ vary in a polynomial family. As shown in \cite[Section 5.4]{babyV} this family gives rise to a polynomial section of $\C^*$-fixed points in the family $\{ \mathcal{X}_{\frac{1}{2}\btheta}(n) : \btheta \in \C^{\ell} \}$. Restricting this to $\btheta \in \Q^{\ell}$ and then applying the inverse to the mapping $\Psi$ of \eqref{rotation} proves the lemma.
\end{proof}
From now on we will always use this labelling of $\C^*$-fixed points by $\mathcal{P}(\ell, n)$.
\subsection{Morse function} \label{morse}
The $U(1)$-action on $\mathcal{M}_{\btheta}(n)$ is hamiltonian with respect to the real symplectic form $\omega_{\R}$ and has moment map $\mu:\mathcal{M}_{\btheta}(n)\longrightarrow (\lie U(1))^*$. If we evaluate this moment map at $-2\sqrt{-1} \in \lie U(1)$ then the corresponding mapping is given by $$f_{\btheta}({\bf X},{\bf Y};v,w) = \sum_{r=0}^{\ell -1} \Tr (X_rX_r^{\dag} - Y_r Y_r^{\dag})$$ for $({\bf X},{\bf Y};v,w) \in \mathcal{M}_{\btheta}(n)$.

\begin{lem} For any ${\bf h}\in \Hr$ we have $$f_{\btheta}( x_{\btheta}(\blambda)) \doteq c_{\bf h} ( {\blambda})\qquad \text{for all }\blambda\in \mathcal{P}(\ell ,n).$$
\end{lem}
\begin{proof}
By definition, the moment map $\mu$ above is induced from the moment map for the $U(1)$-action on $R(e_{\infty} + n\delta)$. For $z\in \lie U(1)$ this can be written as $$\mu ({\bf X},{\bf Y};v,w) (z) = \frac{1}{2} \omega_{\R} ( z\cdot ({\bf X},{\bf Y};v,w) , ({\bf X},{\bf Y};v,w) ).$$ Hence by Lemma \ref{rotate} the function can be calculated on $\mathcal{X}_{\frac{1}{2}\btheta}(n)$ instead via the formula \begin{equation} \label{half} \mathcal{X}_{\frac{1}{2}\btheta}(n) \ni ({\bf X},{\bf Y};v,w) \mapsto \left(z \mapsto \frac{1}{2} (\ima \omega_{\C}) (z\cdot ({\bf X},{\bf Y};v,w), ({\bf X},{\bf Y};v,w) )\right).\end{equation} Now the $\C^*$-action on $\mathcal{X}_{\frac{1}{2}\btheta}(n)$ is hamiltonian with respect to the complex form $\omega_{\C}$ and has moment map $\tilde{\mu} : \mathcal{X}_{\frac{1}{2}\btheta}(n) \longrightarrow (\lie \C^*)^*$  induced from the mapping on $R(e_{\infty} + n\delta)$ given by $$({\bf X},{\bf Y};v,w) \mapsto \left(z \mapsto \frac{1}{2} \omega_{\C} (z\cdot ({\bf X},{\bf Y};v,w), ({\bf X},{\bf Y};v,w) )\right)$$ for $z\in \lie \C^*$. Evaluating this at $-2\sqrt{-1}  \in \lie U(1)$ gives a function $\tilde{\mu}(-2\sqrt{-1}): \mathcal{X}_{\frac{1}{2}\btheta} (n) \longrightarrow \C$ whose imaginary part equals $f_{\btheta}$. We now calculate this function.

By Theorem \ref{mo's} there is a Poisson isomorphism between $\mathcal{X}_{\frac{1}{2}\btheta}(n)$ and $\spc Z_{0 , {\bf h}/2}$. Hence the function $\tilde{\mu}(-2\sqrt{-1})$ corresponds to an element $F_{\btheta} \in Z_{0, {\bf h}/2}$ such that the derivation $ \{ F_{\btheta} , \cdot \}$ equals the derivation on $\C [\mathcal{X}_{\frac{1}{2}\btheta}(n)] = Z_{0,  {\bf h}/2}$ induced by the action of $-2\sqrt{-1} \in \lie U(1)$.  The element $-2\sqrt{-1}$ acts as follows: $(-2\sqrt{-1} )\cdot ({\bf X},{\bf Y};v,w) = (-2\sqrt{-1} {\bf X}, 2\sqrt{-1} {\bf Y}; v,w)$. Let ${\bf z}_t = - \sum_{r=1}^n X_rY_r - z \in H_{t, {\bf h}/2}$ where $$z = \frac{1}{2} \left[\sum_{g\in S_0} h(1-g) + \sum_{r=1}^{\ell -1}\sum_{g\in S_r} (\sum_{s=1}^{\ell -1} \eta^{-rs} (H_1 + \ldots + H_s))g \right],$$ where $S_0$ (respectively $S_r$) is the conjugacy class of elements $(i \, j) \gamma_i \gamma_j^{-1}\in G$ with $1\leq i < j \leq n$ and $\gamma \in \mu_{\ell}$ (respectively $\sigma_i^r \in G$ for $1\leq i \leq n$). By \cite[Section 3.1, (4)]{GGOR} this satisfies the relations
$$[{\bf z}_t , X_i] = tX_i \quad \text{and} \quad  [{\bf z}_t, Y_i] = -tY_i$$ for all $1\leq i \leq n$. Hence, by the definition of the Poisson bracket on $Z_{0, {\bf h}/2}$, we see that $F_{\btheta} = -2\sqrt{-1} {\bf z}_0$, up to the addition of a scalar.

The above shows that $f_{\btheta} (x_{\btheta}(\blambda))$ equals the imaginary part of $-2\sqrt{-1} {\bf z}_0$, in other words the real part of $-2 {\bf z}_0$, evaluated at the fixed point of $\spc Z_{0, {\bf h}/2}$ corresponding to $\blambda$. This evaluation is simply the calculation of the scalar by which $-2 {\bf z}_0$ acts on the baby Verma module $M_{{\bf h}/2}(\blambda)$. Since the elements $Y_i$ kill the generator of $M_{{\bf h}/2}(\blambda)$ it follows that the evaluation of $-2 {\bf z}_0$ is just the same as the evaluation of the real part of the scalar by which the central element $2 z\in \C G$ acts on the irreducible representation corresponding to $\blambda$. However, this is exactly the definition of the $c$-function, \cite[Sections 2,5 and 6]{rou}; taking the real part here is unnecessary as the scalar already belongs to $\R$.
\end{proof}
\subsection{Geometric ordering} \label{geomorder} We now stratify part of $\mathcal{M}_{\btheta}(n)$ by studying the attracting sets of the $\C^*$-action. To this end let ${\bf h}\in \Hr$ and recall from \ref{resol} the resolution $\pi_{\btheta} : \mathcal{M}_{\btheta}(n) \longrightarrow V/G.$ Set \label{Zt-defn}$$\mathcal{Z}_{\btheta} := \pi_{\btheta}^{-1} ( \h\times \{ 0 \} / G),$$ a closed subvariety of $\mathcal{M}_{\btheta}(n)$.
\begin{lem} Keep the above notation and let ${\bf h}\in \Hr$.
\begin{enumerate}
\item $\mathcal{Z}_{\btheta}$ is lagrangian (with respect to the
complex symplectic form) in $\mathcal{M}_{\btheta}(n)$. \item
$\mathcal{Z}_{\btheta}$ is the disjoint union of locally closed
$n$--dimensional affine spaces labeled by $\mathcal{P}(\ell , n)$
$$\mathcal{Z}_{\btheta} =  \coprod_{\blambda\in\mathcal{P}(\ell
,n)} \mathcal{Z}_{\blambda}^{o}.$$ In particular the irreducible
components of $\mathcal{Z}_{\btheta}$ are the Zariski closures of
$\mathcal{Z}_{\blambda}^{o}$, which we denote by
$\mathcal{Z}_{\blambda}$. \item Let $\prec_{{\bf h}}$ be the
partial order on $\mathcal{P}(\ell ,n)$ generated by the rule
$$\bmu \prec_{{\bf h}} \lambda \text{ if } \blambda \neq \bmu
\text{ and } \mathcal{Z}_{\blambda}\cap \mathcal{Z}_{\bmu}^{o}
\neq \emptyset.$$ Then $\bmu \prec_{{\bf h}} \blambda$ implies
that $\bmu <_{\bf h} \blambda.$
\end{enumerate}
\end{lem}
\begin{proof}
This argument follows \cite[Sections 5.1 and 7.1]{nakhilb} very closely.

Set $M = \mathcal{M}_{\btheta}(n)$. Let $x\in M$ be a $\C^*$--fixed point and $T_xM$ the tangent space of $M$ at $x$. There is an induced action of $\C^*$ on this space, so we can decompose it as $$T_xM = \bigoplus_{i \in \mathbb{Z}} T(i)$$ where $T(i) = \{ t\in T_xM: z \cdot t = z^i t \text{ for all }z\in \C^* \}.$ By \ref{circleaction} the complex symplectic form on $M$ is $\C^*$-equivariant. Therefore we see in the weight decomposition of $T_xM$ that $\omega$ must pair together $T_x(i)$ and $T_x(-i)$. Furthermore, since $x$ is an isolated fixed point by Lemma \ref{fix} we have $T_x(0)= 0$. Thus $T_x(M)$ decomposes under the $\C^*$--action into two halves: $\oplus_{i < 0} T(i)$ and $\oplus_{i>0} T(i)$.

Let $\blambda \in \mathcal{P}(\ell ,n)$ and let $x_{\btheta}(\blambda) \in M$ be the corresponding fixed point. The {\it attracting set} of $x_{\btheta}(\blambda)$ is defined as $$\mathcal{Z}_{\blambda}^{o} := \{ x\in M : \lim_{z\to 0} z \cdot x = x_{\btheta}(\blambda)\}.$$ By \cite[Theorem 4.1]{bial} $\mathcal{Z}_{\blambda}^{o}$ is an affine space of dimension $n$ since $\dim \oplus_{i>0} T(i) = n$. The argument of \cite[Proposition 7.1]{nakhilb} shows that $\mathcal{Z}_{\blambda}^{o}$ is lagrangian. Thus (1) is proved.

We now claim that $ \mathcal{Z} = \coprod_{\blambda} \mathcal{Z}_{\blambda}^{o}.$ It is clear that the right hand side is indeed a disjoint union, since a convergence point is unique if it exists. Moreover, if $x\in M$ then $\pi_{\btheta} (z\cdot x) = z\cdot \pi_{\btheta}( x)$ converges if and only if $\pi_{\btheta}(x) \in \C^n\times \{0\}/G$. Thus each $\mathcal{Z}_{\blambda}^{o}$ is contained in the left hand side. Finally, if $x\in \mathcal{Z}$ then $\lim_{z\to 0} z\cdot \pi_{\btheta}(x) =0$ and so $\mathcal{Z}$ contracts under the $\C^*$--action to the {\it projective} variety $\pi_{\btheta}^{-1}(0)$. Therefore every point in $\mathcal{Z}$ converges to a limit. Thus the left hand side is contained in the right hand side too. We have now proved (2).

Part (3) follows from Morse theory. The attracting set can be identified as $\{ x\in M: \lim_{t\to -\infty} e^{\sqrt{-1}t\xi}\cdot x = x_{\btheta}(\blambda) \}$ where $\xi = -2 \sqrt{-1} \in \lie U(1)$ is the infinitesimal generator of $U(1)$ we chose in \ref{morse}. By a standard argument, see for instance \cite[5.1]{nakhilb}, this set can then be identified with the stable manifold of $x_{\btheta}(\blambda)$ which is defined as $$\{ x \in M: \lim_{t\to -\infty} \phi_t(x) = x_{\btheta}(\blambda) \}$$ where $\phi_t$ is a gradient flow of $f_{\btheta}$ with respect to the metric $g$ on $M$. Since we can approximate the function $f_{\btheta}$ around $x_{\btheta}(\blambda)$ by $$T_{x_{\btheta}(\blambda)}M \ni v = \sum_{i} v_i \mapsto \frac{\ell}{2} \sum_i i||v_i||^2$$ we see that in a neighbourhood $U\subseteq \mathcal{Z}_{\blambda}^o$ of $x_{\btheta}(\blambda)$ $f_{\btheta}$ is minimised at $x_{\btheta}(\blambda)$. Now if we consider the function function $\Phi : M\times \R \longrightarrow \R$ sending $(x,t)$ to $f_{\btheta} (\phi_t(x))$ we see that $$\frac{d\Phi}{dt} = \grad f_{\btheta} \left(\frac{d\phi_t}{dt}\right) = g(\grad f_{\btheta}, \grad f_{\btheta}) = ||\grad f_{\btheta}||_{\phi_t(x)}^2.$$ By the above description $d\Phi/dt$ is non-negative in a neighbourhood of $(x_{\blambda}(\btheta),0) \in \mathcal{Z}_{\blambda}^o \times \R$. But since $x_{\blambda}(\btheta)$ is the unique critical point of $f_{\btheta}$ in $\mathcal{Z}_{\blambda}^o$ it follows that $d\Phi/dt$ is non-negative on all of $\mathcal{Z}_{\blambda}^0\times \R$ and hence that $f_{\btheta}(x_{\btheta}(\blambda)) \leq f_{\btheta}(y)$ for all $y\in \mathcal{Z}_{\blambda}^o$, with equality if and only if $y=x_{\btheta}(\blambda)$.

Now suppose that $\bmu \prec_{{\bf h}} \blambda$. Then necessarily
$x_{\btheta}(\mu) \in \mathcal{Z}_{\blambda}$. If
$f_{\btheta}(x_{\btheta}(\bmu)) <
f_{\btheta}(x_{\btheta}(\blambda))$ there would be a neighbourhood
$U'$ of $x_{\btheta}(\bmu)$ in $\mathcal{Z}_{\blambda}$ with
$f_{\btheta}(y) < f_{\btheta}(x_{\btheta}(\blambda))$ for all
$y\in U'$. This neighbourhood would necessarily intersect
$\mathcal{Z}_{\blambda}^o$, contradicting the above paragraph.
Thus, $f_{\btheta}(x_{\btheta}(\blambda)) \leq
f_{\btheta}(x_{\btheta}(\bmu)).$ By Lemma \ref{morse} this implies
that $c_{\bf h}({\blambda}) \leq c_{\bf h}({\blambda})$ and hence,
by definition $\bmu \leq \blambda$ in the $c$-ordering. Since
$\bmu \neq \blambda$ we must have that $c_{\bf h'}(\blambda)
\neq  c_{\bf h'}(\bmu)$ for all ${\bf h}'$ not lying on a
$c$-wall. Thus if $\bf h$ does not lie on such a wall we have
$c_{\bf h}({\blambda}) <  c_{\bf h}({\bmu})$ and (3) is proved.

Now suppose for a contradiction that $\bmu\prec_{{\bf h}}
\blambda$ and $c_{\bf h}({\blambda}) = c_{\bf h}({\bmu})$. Then
$\bf h$ lies on a $c$-wall which is in the interior of a G.I.T.
chamber. Thus there are points on both sides of the $c$-wall which
belong to the G.I.T. chamber containing $\bf h$. Such points on
one side of this wall must have $c_{\bf h'}(\blambda) <
c_{\bf h'}(\bmu)$; those on the other side $c_{\bf h'}(\blambda) > c_{\bf h'}(\bmu)$. However, by \ref{constancy} below, we
have $\bmu \prec_{{\bf h}'} \blambda$ and so by the above
paragraph we have $c_{\blambda}({\bf h}') < c_{\bmu}({\bf h}')$.
This contradiction concludes the proof of (3).
\end{proof}

\subsection{Constancy on G.I.T. chambers} \label{constancy} We finish this section by remarking that the geometric ordering is constant on the interior of G.I.T. chambers. This follows from the fact that thanks to the representation theoretic description of $\mathcal{M}_{\btheta}(n)$, the family of varieties $\mathcal{M}_{\btheta}(n)$ as $\btheta$ varies inside a chamber are all isomorphic as $\C^*$-varieties. Note that the $c$-ordering is in general not constant on G.I.T. chambers since such a chamber may be a union of a number of different $c$-chambers, see Figures \ref{cch} and \ref{GITch} in the introduction.

\section{Combinatorics} In this section we describe a classical combinatorial algorithm which sets up a bijection between $\ell$-multipartitions $\mathcal{P}(\ell, n)$ and partitions $\mathcal{P}(n)$, see \cite[2.7]{JAMESANDKERBER} and \cite{lecmiy}. The bijection depends on a {\it multi-charge} ${\bf s}\in \mathbb{Z}^{\ell}$.

\subsection{$\beta$-numbers of a partition}
Let $\lambda \in \mathcal{P}(n)$ and let $s\in \Z$. We associate a set of strictly decreasing positive integers which are called the $\beta$-numbers of $\lambda$ and depend on $s$:\label{beta-defn}
$$\beta_s(\lambda) = \{ \lambda_1 +s, \lambda_2  + s -1 , \ldots , \lambda_{j} + s +1-j, \ldots  \}.$$
It's clear that we can recover $\lambda$ from this set of distinct integers since the sequence eventually stabilises to $s+1-j$ for large $j$.
\subsection{Multipartitions to partitions} \label{multpart}Given ${\bf s} = (s_1, \ldots , s_{\ell})\in \Z^{\ell}$ set $s = \sum_{i=1}^{\ell} s_i$. We define a bijection between $\ell$-multipartitions and partitions (of various degrees $n$).

Let $\blambda =(\lambda^{(1)}, \ldots , \lambda^{(\ell)})$ be an $\ell$-multipartition. Define a a set of distinct integers as follows $$ \bigcup_{i=1}^{\ell} \{ \ell (x-1) + i : x\in \beta_{s_{i}}(\lambda^{(i)}) \}.$$ The elements of this set eventually stabilise to $s+1-j$ and so it equals the $\beta_s(\lambda({\bf s}))$ for some partition $\lambda({\bf s})$. This process yields a bijection \label{tau-defn}$$\mathbb{Z}^{\ell}_0 \times \coprod_n \mathcal{P}(\ell , n) \longrightarrow \coprod_n \mathcal{P}(n), \qquad ({\bf s}, \blambda) \mapsto \tau_{\bf s}(\blambda), $$ where  $\zl =\{ (s_1,\ldots , s_{\ell}) \in \mathbb{Z}^{\ell} : s_1+\cdots + s_{\ell}=0 \}.$

\subsection{$\ell$-cores} \label{lcores}(See \cite{lecmiy}.) Let ${\bf s}\in \zl$. Recall that an $\ell$-core is a partition from which no outer rim-hooks of length $\ell$ can be removed. If we take the trivial $\ell$-multipartition ${\bempty} :=(\emptyset, \ldots , \emptyset)$, then $\tau_{{\bf s}}(\bempty)$ is an $\ell$-core, \cite[2.1]{lecmiy}. This sets up a bijection between $\ell$-cores and $\zl$.

\subsection{Contents of a partition}\label{young}
We identify a partition $\lambda$ with its Young diagram $\lambda = \{ (p,q)\in \Z_{\geq 0} \times \Z_{\geq 0}: p\leq \lambda_q\}.$ For example $\lambda = (5,5,3,1,1)$ gives
$$
\begin{array}[c]{cccccc}
\bullet \\
\bullet \\
\bullet & {\bullet } & \bullet  \\
\bullet & \bullet  & \bullet &
    \bullet&    \bullet & \\
\llap{${}_{(0,0)}$} \bullet &   \bullet &   \bullet &
\bullet &   \bullet
\end{array}
\qquad \qquad\lambda =(5,5,3,1,1).$$Boxes of the diagram are called nodes. We label each node of $\lambda$ with an integer called its content. By definition, the content of the node $(p,q)$ is \label{cont-defn}$\res(p,q) = p-q$. For $0\leq i \leq \ell -1$ we define \label{N-defn}$\cont_i(\lambda)$ to be the number of nodes of $\lambda$ whose content equals $i$ modulo $\ell$. For the example above taken with $\ell =3$ we have $\cont_0 (\lambda) = 6, \cont_1(\lambda ) = 4$ and $\cont_2(\lambda) = 5$. By \cite[Theorem 2.7.41]{JAMESANDKERBER} if $\lambda$ and $\mu$ are partitions of the same degree then they have the same $\ell$-core if and only if $\cont_i (\lambda) = \cont_i(\mu)$ for all $0\leq i \leq \ell -1$.

\subsection{Monomial ideals} \label{monomial}We can associate to any partition $\lambda$ of $n$ an ideal \label{I-defn}$I_{\lambda}$ of colength $n$ in the polynomial ring $\C[A,B]$. By definition $I_{\lambda}$ is the ideal spanned by the monomials $\{ A^pB^q : (p,q)\notin \lambda\}$. If we let $\mu_{\ell}\leq SL_2(\C)$ act algebraically on $\C[A,B]$ by $\sigma \cdot A = \eta A, \sigma \cdot B = \eta^{-1} B$ we see that $\C[A,B]/I_{\lambda}$ is an $n$-dimensional representation of $\mu_{\ell}$ whose character is $\sum_{i=0}^{\ell -1} \cont_{i}(\lambda) \alpha_i$.

\section{Combinatorial description of geometric ordering}\label{combgeom}
Throughout this section we will take ${\bf h} \in {\bf H}$ and enforce the equality $\btheta = (-h+H_0, H_1, \ldots , H_{\ell -1}).$ We will give a combinatorial description of the geometric ordering in any G.I.T. chamber. We mostly follow the approach of Haiman, \cite[7.2]{hai}, where many of the results here first appeared.
\subsection{Reduction} \label{postoneg}
We begin with a simple lemma which allows to reduce to the case $h>0$. To describe it we need two pieces of notation. We let $\overline{{\bf h}} = ( -h, -H_{\ell -1}, \ldots , -H_{1} )$ and, given $\blambda \in \mathcal{P}(\ell ,n)$, we set $\overline{\blambda} = ({}^t\lambda^{(\ell)}, \ldots, {}^t\lambda^{(1)}).$
\begin{lem}
Suppose that ${\bf h}\in \Hr$. Then $\overline{{\bf h}}\in \Hr$ and for $\blambda ,\bmu \in \mathcal{P}(\ell ,n)$ we have $$\blambda \prec_{{\bf h}} \mu \text{ if and only if } \overline{\blambda} \prec_{\overline{{\bf h}}} \overline{\bmu}.$$
\end{lem}
\begin{proof}
It is trivial to check that $\overline{{\bf h}}\in \Hr$.

We have $\overline{\btheta} = (-\theta_0, -\theta_{\ell -1}, \ldots , -\theta_1)$. Consider the automorphism $\phi: R(e_{\infty}+n\delta) \longrightarrow R(e_{\infty}+n\delta)$ that is defined by $$  (X_0, \ldots, X_{\ell-1},Y_0,  \ldots ,Y_{\ell -1}
; v, w)\mapsto (X_{\ell-1}^T, \ldots , X_{0}^T,  -Y_{\ell-1}^T, \ldots , -Y_0^T; w^T, v^T).$$ It is straightforward to check that $\phi$ preserves the relation $[{\bf X},{\bf Y}]+vw = 0$ and so it restricts to a automorphism from $\mu_\C^{-1}(0)$ to $\mu^{-1}_\C(0)$. Moreover if $\mu_{\R} (({\bf X}, {\bf Y}; v,w)) = \frac{\sqrt{-1}}{2}\btheta$ then $\mu_{\R} (\phi({\bf X}, {\bf Y}; v,w)) = \frac{\sqrt{-1}}{2}\overline{\btheta}$. Hence $\phi$ induces an isomorphism $$\phi: \mathcal{M}_{\btheta}(n)\longrightarrow \mathcal{M}_{\overline{\btheta}}(n).$$ By construction, $\phi$ is equivariant for the $\C^*$-actions on its domain and codomain. In particular fixed points get sent to fixed points, so there is an associated bijection $\phi^* : \mathcal{P}(\ell ,n) \longrightarrow \mathcal{P}(\ell,n)$. Furthermore, as attracting sets get sent to attracting sets it follows that $\blambda \prec_{\bf h} \bmu$ if and only if $\phi^*(\blambda) \prec_{\overline{{\bf h}}} \phi^*(\bmu).$

It now remains to show that $\phi^*(\blambda) = \overline{\blambda}$ for all $\blambda \in \mathcal{P}(\ell ,n)$. To do this we will first compare the $c$-functions on $\mathcal{M}_{\btheta}(n)$ and $\mathcal{M}_{\overline{\btheta}}(n)$. \begin{eqnarray*} c_{\overline{\bf h}}(\overline{\blambda}) & = & \sum_{r=2}^{\ell} |{}^t\lambda^{(\ell+1-r)}|(-H_{\ell -1} - \cdots - H_{\ell +1-r}) + \left( \frac{n(n-1)}{2} + \sum_{r=1}^{\ell} n({}^t\lambda^{(r)})-n(\lambda^{(r)})\right)h \\ & = & \sum_{r=2}^{\ell} |\lambda^{(\ell+1-r)}|(H_0 + \cdots + H_{\ell - r}) - \left( \frac{n(n-1)}{2} + \sum_{r=1}^{\ell } n(\lambda^{(r)})-n({}^t\lambda^{(r)})\right)h + n(n-1)h \\ & = & \sum_{r=1}^{\ell -1} |\lambda^{(r)}|(H_0 + \cdots + H_{r-1}) - \left( \frac{n(n-1)}{2} + \sum_{r=1}^{\ell} n(\lambda^{(r)})-n({}^t\lambda^{(r)})\right)h + n(n-1)h \\ & = & \sum_{r=2}^{\ell -1} |\lambda^{(r)}|(H_1 + \cdots + H_{r-1}) + \sum_{r=1}^{\ell -1} |\lambda^{(r)}|H_0 - \left( \frac{n(n-1)}{2} + \sum_{r=1}^{\ell} n(\lambda^{(r)})-n({}^t\lambda^{(r)})\right)h + n(n-1)h \\ & = & \sum_{r=2}^{\ell} |\lambda^{(r)}|(H_1 + \cdots + H_{r-1}) + nH_0 - \left( \frac{n(n-1)}{2} + \sum_{r=1}^{\ell} n(\lambda^{(r)})-n({}^t\lambda^{(r)})\right)h + n(n-1)h \\ & = & c_{\bf h}({\blambda}) + n(n-1)h + nH_0\end{eqnarray*} Thus by Lemma \ref{morse} we have that the Morse function of \ref{morse} satisfies $$f_{\btheta}(x_{\btheta}(\blambda)) \doteq  f_{\overline{\btheta}}(x_{\overline{\btheta}}(\overline{\blambda})) \qquad \text{for all }\blambda\in\mathcal{P}(\ell ,n).$$ On the other hand it follows from the definition of $f_{\btheta}$ and of the isomorphism $\phi$ that $$f_{\btheta}(x_{\btheta}(\blambda)) = f_{\overline{\btheta}}(x_{\overline{\btheta}}(\phi^*(\blambda))) \qquad \text{for all }\blambda\in\mathcal{P}(\ell, n).$$ Since the set of values of the $c$-function evaluated at elements of $\mathcal{P}(\ell ,n)$ are all distinct at a generic point of a G.I.T. chamber by Theorem \ref{RouqvsGIT}, it follows that $\phi^*(\blambda) = \overline{\blambda}$. \end{proof}

\subsection{Affine symmetric group and G.I.T. chambers} \label{stabtoweights}The spaces $\mathcal{M}_{\btheta}(n)$ are invariant under scaling $\btheta$ by $\Q_{+}$. This, together with Lemma \ref{postoneg}, means that without loss of generality we can assume that $- h = \theta_0 +\ldots + \theta_{\ell - 1} = 1.$ Thus we consider the space of stability parameters \label{slice-defn}$$\Theta_1 = \{ \btheta =  (\theta_0, \ldots , \theta_{\ell -1}) \in \Q^{\ell} : \theta_0 + \cdots + \theta_{\ell -1} = 1 \}$$ and the corresponding subset of ${\bf H}$. The set of walls \eqref{simples} then becomes \begin{equation} \label{affwalls} \{\btheta \cdot \beta = m : \beta \in R^+ \text{ and } 1-n\leq m \leq n-1\}, \end{equation} where $R^+$ is the set of positive roots in the root system of type $A_{\ell-1}$, that is elements of the form $e_i + \cdots + e_j$ for $1\leq i\leq j \leq \ell -1$.

Let $X\otimes_{\Z} \Q$ be the weight lattice of type $A_{\ell -1}$ spanned by fundamental weights $\varpi_i$ for $1\leq i \leq \ell-1$. Recall that we also have the coroot lattice $\Z R^{\vee} \subset X$. We fix a basepoint $(1,0,\ldots ,0) \in \Theta_1$ and thus identify $\Theta_1$ with the weight lattice $X\otimes_{\Z}\Q$ of type $A_{\ell -1}$ via $\btheta \mapsto \sum_{i=1}^{\ell-1}\theta_i\varpi_i$.

Let \label{affine-defn}$\tilde{\mathfrak{S}}_{\ell} = \mathfrak{S}_{\ell} \ltimes \Z R^{\vee}$ denote the affine symmetric group. It has a Coxeter presentation $$\tilde{\mathfrak{S}}_{\ell} = \langle \sigma_0, \ldots , \sigma_{\ell -1} : \sigma_i^2 = 1 \text{ and }\sigma_i\sigma_{i+1}\sigma_i = \sigma_{i+1}\sigma_i\sigma_{i+1}\text{ for $0\leq i\leq \ell-1$}\rangle$$ where subscripts are always counted modulo $\ell$. It acts naturally on $X\otimes_{\Z}\Q$ and hence it acts on $\Theta_1$ by the rule \begin{equation} \label{symaction}\sigma_i \cdot  (\theta_0 , \ldots , \theta_{\ell -1}) = (\theta_0 , \ldots, \theta_{i-1}+\theta_i, -\theta_i , \theta_{i}+\theta_{i+1}, \ldots , \theta_{\ell -1}) \qquad \text{for all }0\leq i \leq \ell -1.\end{equation} The reflecting hyperplanes for this action are determined by the roots of $\tilde{A}_{\ell -1}$ and so give the set $$\{ \btheta \cdot \beta = m  : \beta \in R^+ \text{ and } m\in \Z \}.$$ The connected components of the space (tensored over $\R$) obtained by removing these hyperplanes are called alcoves. The closure of any alcove is a fundamental domain for the action of $\tilde{\mathfrak{S}_{\ell}}$ on $\Theta_1$. We deduce the following result from Lemma \ref{GITcham}.

\begin{lem} \label{alcoves} The G.I.T. walls in $\Theta_1$ are a subset of the reflecting hyperplanes for the action of $\tilde{\mathfrak{S}}_{\ell}$. Moreover, in the limit as $n\to \infty$, the G.I.T. chambers and the alcoves agree.
\end{lem}

\subsection{Tactics} Lemma \ref{alcoves} combined with \ref{constancy} shows that in order to determine the geometric ordering in general it is sufficient to pick a point $\btheta \in \Theta_1$ in an alcove and then determine the ordering at all points in its orbit, $\tilde{\mathfrak{S}}_{\ell}\cdot \btheta$. The point we pick is $\btheta = {\bf 1}$, where ${\bf 1} = \frac{1}{\ell}(1, \ldots , 1)$. We will first describe the geometric ordering at the points $\Z R^{\vee} \cdot \bf 1$ using reflection functors, \cite{Nak3}, and then a simple argument will describe how the ordering varies under the remaining action of $\mathfrak{S}_{\ell}$.

\subsection{Translation by the coroot lattice} Fix an element of the coroot lattice $$\beta^{\vee} = \sum_{i=1}^{\ell -1} a_i \alpha_i^{\vee} \in \Z R^{\vee}$$ where $a_i\in \Z$ for all $i$, and let $\tau_{\beta^{\vee}} \in \tilde{\mathfrak{S}}_{\ell}$ denote the corresponding translation of $X\otimes_\Z \Q$, or equally $\Theta_1$.

Let $\Z\Phi$ be the root lattice of $Q$, i.e. the affine root lattice of type $\tilde{A}_{\ell -1}$.  There is an action on $\tilde{\mathfrak{S}}_{\ell}$ on $\Z\Phi$. We shift the origin of this action, defining  \label{ast-defn}$$\sigma_i \ast \gamma = \sigma_i\gamma + \delta_{i0}\alpha_0$$ for $\gamma \in \Z\Phi$ and for $0\leq i\leq \ell -1$.

\begin{lem} \label{reffun} Let $\btheta = \tau_{-\beta^{\vee}} \cdot {\bf 1}\in \Theta_1$. There is a $U(1)$-equivariant hyper-K\"ahler isometry between $\mathcal{M}_{\btheta}(n)$ and $\mathcal{M}_{{\bf 1}} (\tau_{\beta^{\vee}}\ast n\delta)$.
\end{lem}
\begin{proof}
Let $\tau_{\beta^{\vee}} = \sigma_{i_1}\sigma_{i_2} \cdots
\sigma_{i_k}$ be a reduced expression in
$\tilde{\mathfrak{S}}_{\ell}$. Repeated application of the simple
reflection functors of \cite[Theorem 6.1]{Nak3} yields a
hyper-K\"ahler isometry $$\mathcal{M}_{\btheta}(n) = \mathcal{M}_{\sigma_{i_k} \cdots \sigma_{i_2}\sigma_{i_1} \bf
1}(n\delta) \longrightarrow \mathcal{M}_{\bf
1}(\sigma_{i_1}\sigma_{i_2} \cdots \sigma_{i_k}\ast n\delta) = \mathcal{M}_{\bf 1}(\tau_{\beta^{\vee}}\ast n\delta).$$ Moreover
it follows immediately from their definition, \cite[4(iii)]{Nak3},
that the simple reflection functors are $U(1)$-equivariant.
\end{proof}

\subsection{} \label{basic} We need a little combinatorics now.
\begin{lem}
Let $\btheta = \tau_{-\beta^{\vee}}\cdot  {\bf 1}\in \Theta_1$ and $\gamma = \tau_{\beta^{\vee}}\ast n\delta$.
\begin{enumerate}
\item[(i)] There exists $\gamma_0 \in \NN\Phi_+$ such that $\gamma = \gamma_0 + n\delta$.
\item[(ii)] There is a unique partition $\nu$ such that $\gamma_0 = \sum_{i=0}^{\ell -1} \cont_i(\nu) \alpha_i$. It is an $\ell$-core.
\item[(iii)] Under the bijection of \ref{lcores} the element ${\bf s}\in \zl$ corresponding to $\nu$ is related to $\theta$ by $$\btheta = {\bf 1} + (s_1-s_{\ell} , s_2-s_1, \ldots , s_{\ell}-s_{\ell-1}).$$
\end{enumerate}
\end{lem}
\begin{proof}
(i) By \cite[Theorem 10.2]{Nak2} the space $\oplus_{{\bf d}} H^{\text{top}}(\mathcal{M}_{\bf 1}({\bf d}), \C)$ has a geometrically defined action of $U(\hat{\mathfrak{sl}}_{\ell})$ which makes it isomorphic to $V(\varpi_0)$ where $\varpi_0$ is the fundamental weight corresponding to the extending vertex of $Q$. This is called the basic representation of $U(\hat{\mathfrak{sl}}_{\ell})$. In this description $H^{\text{top}}(\mathcal{M}_{\bf 1}(\nu),\C)$ is the weight space $V(\varpi_0)_{\varpi_0 - \nu}$. By definition we have for any $w\in \tilde{\mathfrak{S}}_{\ell}$ $$w\ast 0 = \varpi_0 - w\varpi_0.$$ Moreover, since $\delta$ is an isotropic vector we see that $w\ast n\delta = (w\ast 0 ) + n\delta$. Thus $\gamma_0 = \varpi_0 - \tau_{\beta^{\vee}}\varpi_0$ must be an element of $\NN\Phi_+$ since all weights of $V(\varpi_0)$ differ from $\varpi_0$ by combinations of positive roots of $U(\hat{\mathfrak{sl}}_{\ell})$.

(ii) By \cite{mismiw} there is a combinatorial basis of $V(\varpi_0)$ labelled by $\ell$-regular partitions. In this description a partition $\lambda$ has weight $\sum_{i=0}^{\ell -1} \cont_i(\lambda) \alpha_i$, \cite[Section 2]{mismiw}. The weight space $V(\varpi_0)_{\varpi_0 - \gamma_0}$ is one-dimensional since it is a $\tilde{\mathfrak{S}}_{\ell}$-conjugate of $\varpi_0$ and so there is a unique $\ell$-regular partition $\nu$ associated to $\gamma_0$. By \cite[5.3]{LLT} it is an $\ell$-core.

(iii) Recall that $\beta^{\vee} = \sum_{i=1}^{\ell -1} a_i\alpha_i^{\vee}$. Explicit calculation following \ref{stabtoweights} shows that $$\btheta = \tau_{-\beta^{\vee}}\cdot {\bf 1} = {\bf 1} + (a_1+a_{\ell -1}, - 2a_1+a_2, a_1-2a_2 +a_3, \ldots , a_{\ell -3}-2a_{\ell -2} + a_{\ell -1}, a_{\ell - 2}-2a_{\ell -1}).$$ Analogously, by \cite[2.4]{lecmiy}, there is a transitive action of $\tilde{\mathfrak{S}}_{\ell}$ on $\zl$ which makes it isomorphic to $\tilde{\mathfrak{S}}_{\ell}\varpi_0$ as a $\tilde{\mathfrak{S}}_{\ell}$-set. By \cite[3.1]{lecmiy}, the partition corresponding an element $w\varpi_0$ via the correspondence in \cite{mismiw}  is the same as the partition corresponding to $w{\bf 0}$ under the bijection described in \ref{lcores}. Thus the translate $\tau_{\beta^{\vee}} {\bf 0}$ in $\zl$ is the set of integers which define the $\ell$-core $\nu$. To calculate $\tau_{\beta^{\vee}}{\bf 0}$ we observe from \cite[2.3]{lecmiy} that for $1\leq i \leq \ell-1$  $\tau_{\alpha_i^{\vee}}$ corresponds to translation by $e_i - e_{i+1}$. Thus we find $$\tau_{\beta^{\vee}} {\bf 0} = ( a_1, -a_1+a_2, -a_2+a_3, \ldots , -a_{\ell-2} + a_{\ell-1}, -a_{\ell -1})\in \zl.$$ Part (iii) now follows.
\end{proof}
\subsection{} We next recall Nakajima's construction of Hilbert schemes of points on the plane, \cite[Chapter 2 and Theorem 2.1]{nakhilb}. Let $K$ be any positive integer. Let $\tilde{H}(K)$ by the set of quadruples  $ (A, B; v,w) \in (\mt_K(\C)^{\oplus 2})\oplus \C^n \oplus (\C^n)^*$ which satisfy the condition $[A,B] + vw = 0$ and the following stability condition: there exists no proper subspace $S\subsetneq \C^n$ such that $A(S)\subseteq S$, $B(S) \subseteq S$ and $\im i \subseteq S$. The group $GL_K(\C)$ acts freely on $\tilde{H}(K)$ by the rule $g\cdot (A,B; v,w) = (gAg^{-1}, gBg^{-1}; gv, wg^{-1})$ and the quotient $\tilde{H}(K)/GL_K(\C)$ is the Hilbert scheme of $K$ points on the plane, denoted $\hi{K}$. This is a smooth variety of dimension $2K$.

The description above shows that $\hi{K}$ is a quiver variety for the quiver with one vertex and one loop  with dimension vector $K$ and stability parameter $-1$. Therefore $\hi{K}$ has a hyper-K\"ahler structure and a $\C^*$-action.

\subsection{} Usually the Hilbert scheme $\hi{K}$ is presented as the scheme whose underlying set consists of ideals of $\C [A,B]$ of colength $K$. Here every element of $\tilde{H}(K)$ becomes a cyclic $\C[A,B]$-module with generator $i(1)$. Taking the annihilator of this module then induces an isomorphism between $\hi{K}$ above and this more usual description.

\subsection{} \label{tohilb}The group $\mu_{\ell}$ acts on $\C[A,B]$ as explained in \ref{monomial} and hence induces an action on $\hi{K}$. Thus we can consider $(\hi{K})^{\mu_{\ell}}$: this is smooth and so a union of connected components. Each component inherits a $U(1)$-equivariant hyper-K\"ahler structure from $\hi{K}$.

We have the following lemma which is stated without proof in \cite[Proposition 7.2.8]{hai}.

\begin{lem} Let $\gamma = \tau_{\beta^{\vee}}\ast n\delta$ so that $\gamma = \gamma_0 + n\delta$. Let $\nu$ denote the $\ell$-core corresponding to $\gamma_0$ by Lemma \ref{basic} and set $K = |\nu| + \ell n$.
\begin{enumerate}
\item[(i)] Let $I_{\nu}$ be the monomial ideal of $\C[A,B]$ associated to $\nu$ as defined in \ref{monomial}. There is a component \label{hiv-defn}$\hil{\nu}$ of $(\hi{K})^{\mu_{\ell}}$ whose generic points have the form $V(I_{\nu}) \cup T$ where $T$ is a union of $n$ distinct free $\mu_{\ell}$-orbits in $\C^2$.
\item[(ii)] There is a $U(1)$-equivariant hyper-K\"ahler isometry between $\mathcal{M}_{\bf 1}(\gamma)$ and $\hil{\nu}$.
\end{enumerate}
\end{lem}
\begin{proof}
By \cite[Introduction]{CB} the quiver variety $\mathcal{M}_{\bf 1}(\gamma)$ is connected. We identify $R(\gamma)$ with a $\mu_{\ell}$ invariant subspace of $(\mt_{K}(\C)^{\oplus 2})\oplus \C^n \oplus (\C^n)^*$ by sending $({\bf X}, {\bf Y}; v,w)$ to $(A,B; v,w)$ where $A$ and $B$ are block matrices $$A = \left(\begin{array}{ccccc}0 & 0 & \cdots & 0 & X_{\ell-1} \\X_0 & 0 & \cdots & 0 & 0 \\0 & X_1 & \cdots & 0 & 0 \\\vdots & \vdots &  & \vdots & \vdots \\0 & 0 & \cdots & X_{\ell-2} & 0\end{array}\right), \qquad B = \left(\begin{array}{ccccc}0 & Y_0 & 0 & \cdots & 0 \\0 & 0 & Y_1 & \cdots & 0 \\\vdots & \vdots & \vdots &  & \vdots \\0 & 0 & 0 & \cdots & Y_{\ell-2} \\Y_{\ell-1} & 0 & 0 & \cdots & 0\end{array}\right).$$ Following the arguments of \cite[Theorem 4.4]{nakhilb} and \cite[Theorem 2]{wang} we see that there is an isomorphism from $\mathcal{M}_{\bf 1}(\gamma)$ to a connected component of $(\hi{K})^{\mu_{\ell}}$. The identification is clearly a $U(1)$-equivariant hyper-K\"ahler isometry.

Consider the monomial ideal $I_{\nu}$ -- it is supported on $0\in \C^2$ -- and take $n$ distinct generic $\mu_{\ell}$-orbits in $\C^2$ which we label $\mathcal{O}_1, \ldots , \mathcal{O}_n$. Then, in $\C^2$ the elements $0$ and $ \supp \mathcal{O}_i$ for $1\leq i \leq n$ are all distinct and thus there is a unique $0$-dimensional subscheme of $\C^2$ of colength $n$ associated to this data and hence a point of $\hi{K}$.  Moreover, since the ideal $I_{\nu}$ carries the representation of $\mu_{\ell}$ of dimension vector $\gamma_0 = \sum_{i=0}^{\ell} \cont_i(\nu)\alpha_i$, it follows that this is a point in the image of $\mathcal{M}_{\bf 1}(\gamma)$. The choice of distinct $\mu_{\ell}$-orbits forms a $2n$-dimensional space and hence this set of points is a dense set in the image of $\mathcal{M}_{\bf 1}(\gamma)$ and its closure determines the connected component of $(\hi{K})^{\mu_{\ell}}$. This completes the proof.
\end{proof}

\subsection{} \label{hilbdom}The $\C^*$-fixed points in $\hi{K}$ correspond to the monomial ideals $I_{\lambda} \lhd \C[A,B]$ with $|\lambda| = K$.  Moreover the geometric ordering arising from the attracting sets is the anti-dominance order on partitions, see for instance \cite[End of Section 4]{nakjack}.

\subsection{} \label{bij} We are finally in a position to make our first calculation of the ordering $\prec_{\bf h}$.
\begin{prop}
Let $\btheta = \tau_{-\beta^{\vee}}\cdot \bf{1}$ where $\beta^{\vee} = \sum_{i=1}^{\ell -1} a_i \alpha_i^{\vee}$. Let $\nu$ be the corresponding $\ell$-core and ${\bf s} = (a_1, -a_1+a_2, -a_2+a_3, \ldots , -a_{\ell -2} + a_{\ell-1}, -a_{\ell -1}) \in \zl$. Set $K = |\nu| + \ell n$.
\begin{enumerate}
\item[(i)] The mapping of $\tau_{\bf s}$ of \ref{multpart} restricts to a bijection between $\mathcal{P}(\ell ,n)$ and \label{partnu-defn}$\mathcal{P}_{\nu}(K)$, the set of partitions of degree $K$ having $\ell$-core $\nu$.
\item[(ii)] For all $\blambda, \bmu$ we have $$\blambda \prec_{{\bf h}} \bmu \Longleftrightarrow {}^t\tau_{\bf s}({}^t\blambda) \lhd {}^t\tau_{\bf s}({}^t \bmu).$$
\end{enumerate}
\end{prop}
\begin{proof}
(i) This is well-known, see for example \cite[Lemma 2.7.13 and Theorem 2.7.30]{JAMESANDKERBER}.

(ii) It follows from Lemma \ref{tohilb} that the fixed points that belong to $\hil{\nu}$ are precisely the ideals $I_{\lambda}$ whose $\ell$-core is $\nu$. The geometric ordering on $\hil{\nu}$ is the restriction of the anti-dominance ordering to $\mathcal{P}_{\nu}(K)$.

By Lemmas \ref{reffun} and \ref{tohilb}(ii) there is $U(1)$-equivariant hyper-K\"ahler isometry between $\mathcal{M}_{\btheta}(n)$ and $\hil{\nu} \subseteq (\hi{K})^{\mu_{\ell}}$. By taking fixed points this induces a natural bijection between $\phi: \mathcal{P}(\ell ,n)\longrightarrow \mathcal{P}_{\nu}(K)$ which, thanks to \ref{hilbdom}, intertwines the geometric ordering and the anti-dominance ordering. Thus we must show that $\phi (\blambda) = \tau_{\bf s}({}^t \blambda)$ for all $\blambda \in \mathcal{P}(\ell ,n)$.

Let ${\bep} = (\ep_0, \ldots , \ep_{\ell -1})\in \Q_0^{\ell}$ be sufficiently small so that $\btheta$ and $\btheta+\bep$ belong to the same G.I.T. chamber but sufficiently generic so that $\btheta+\bep$ does not lie on a $c$-wall. Then, as in Lemma \ref{reffun}, there is a $U(1)$-equivariant hyper-K\"ahler isometry between $\mathcal{M}_{\btheta + \bep}(n) = \mathcal{M}_{\tau_{-\beta^{\vee}} \cdot ({\bf 1} + {\bep})}(n)$ and $\mathcal{M}_{\bf 1 + \bep}(\gamma)$ where $\gamma = \tau_{\beta^{\vee}} \ast n\delta$. In particular there is a bijection between the fixed points of these two varieties.

Since ${\bep}$ is chosen to be small, we can assume that all entries of ${\bpsi} :=\bf 1 +  \bep$ are positive. It follows that immediately from the representation theoretic description of quiver varieties that $\mathcal{M}_{\bf 1}(\gamma)$ and $\mathcal{M}_{\bpsi}(\gamma)$ are isomorphic. Since the fixed points in these varieties are isolated, and hence rigid, we see therefore that the description of the fixed points in the two varieties must be the same  combinatorially. However, thanks to Lemma \ref{tohilb} we have an explicit description of the fixed points of $\mathcal{M}_{\bf 1}(\gamma)$ as $(\Gamma, \C^*)$-equivariant representations of the quiver $Q_{\infty}$, compare \cite[5.2]{nakhilb}.  The fixed points are labelled by $\rho \in \mathcal{P}_{\nu}(K)$. At the $i$th vertex of $Q$ we place the $\mu_{\ell}$-isotypic component of $\C[A,B]/I_{\rho}$ corresponding to the irreducible representation $\eta_i$. The grading is inherited from the natural grading on the quotient space. Finally, we set $v(1) = 1 \in \C[A,B]/I_{\rho}$ and $w = 0$.

Now we can calculate the Morse function $F_{\bpsi}$ of $\mathcal{M}_{\bpsi}(\gamma)$ at the fixed point $I_{\rho}$ corresponding to $\rho \in \mathcal{P}_{\nu}(K)$ in the same as way as \cite[Proposition 5.14]{nakhilb}. It is $$F_{\bpsi}(I_\rho) = \sum_{(p,q)\in \lambda} \psi_{\res (p,q)} \res(p,q).$$
 Thus for all $\blambda \in \mathcal{P}(\ell ,n)$ we have by Lemma \ref{morse}\begin{equation*} \sum_{(p,q)\in \phi(\blambda)} \!\!\!\!\!\psi_{\res (p,q)} \res(p,q) \doteq c_{\bf h}(\blambda)  \qquad \text{for all }\blambda\in\mathcal{P}(\ell ,n)\end{equation*} where ${\bf h}$ corresponds to ${\btheta} + {\bep}$ by the usual rule. Since ${\btheta}+{\bep}$ was chosen not to lie on a $c$-wall all values of $c_{\bf h}({\blambda})$ are distinct and thus to show $\phi (\blambda) = \tau_{\bf s}({}^t\blambda)$ it is sufficient to show that \begin{equation} \label{hilbequ}  \sum_{(p,q)\in \tau_{\bf s}({}^t\blambda)} \!\!\!\!\!\psi_{\res (p,q)} \res(p,q) \doteq c_{\bf h}({\blambda} ) \qquad \text{for all }\blambda\in\mathcal{P}(\ell ,n).\end{equation} Being a coward, I have moved this computation to the appendix. You have to be either brave or stupid to go through it.
\end{proof}

\subsection{Symmetric group action}
There is an action of $\mathfrak{S}_\ell$ on $\mathcal{P}(\ell ,n)$: $$w\cdot (\lambda^{(1)}, \ldots , \lambda^{(\ell)}) = (\lambda^{(w\cdot 1)}, \ldots , \lambda^{(w\cdot \ell)}) \qquad \text{for } w\in \mathfrak{S}_\ell$$ So given a partial ordering on $\mathcal{P}(\ell , n)$ and an element $w\in\mathfrak{S}_n$ we can construct a new partial ordering, its $w$-translate.

\begin{lem} Let ${\bf h} \in {\bf H}^{reg}$ and $w\in \mathfrak{S}_\ell$. Let ${\bf h}'\in \Hr$ denote the parameters corresponding to $w\cdot \btheta$. Then the geometric ordering $\prec_{\bf h'}$ is the $w$-translate of the geometric ordering $\prec_{\bf h}$.
\end{lem}

\begin{proof}
It is sufficient to prove this for the simple reflections $\sigma_j$ for $1\leq j \leq \ell -1$.
Moreover, since the geometric ordering is constant on G.I.T. chambers we can assume without loss of generality that $\btheta$ does not lie on a $c$-wall.

We apply the reflection functors again, \cite[Theorem 6.1]{Nak3}, to obtain a $U(1)$-equivariant hyper-K\"ahler isometry between $\mathcal{M}_{\btheta}(n)$ and $\mathcal{M}_{\sigma_j \cdot \btheta}(n)$. Here $\sigma_j \cdot \btheta$ is given by \eqref{symaction} and so $${\bf h'} = (h, H_1, \ldots , H_{j-1}+H_j, -H_j, H_{j}+H_{j+1}, \ldots , H_{\ell-1}).$$ The isometry induces a bijection between $U(1)$-fixed points on $\mathcal{M}_{\btheta}(n)$ and on $\mathcal{M}_{\sigma_j\cdot \btheta}(n)$, so a bijection $\phi_j$ from $\mathcal{P}(\ell,n)$ to itself. Moreover it provides an equality of Morse functions $$f_{\btheta}(x_{\btheta}(\blambda)) \doteq f_{\sigma_j\cdot \btheta}(x_{\sigma_j\cdot \btheta}(\phi_j(\lambda))) \quad \text{for all }\blambda\in\mathcal{P}(\ell ,n).$$ Thus, by Lemma \ref{morse}, to show that $\phi_j$ is given by $\sigma_j$-translation it is enough to show that $$c_{\bf h}({\blambda}) = c_{\sigma_j\cdot \blambda}({\bf h'})  \quad \text{for all }\blambda\in\mathcal{P}(\ell ,n).$$
This is obvious: \begin{eqnarray*} c_{\sigma_j\cdot \blambda}({\bf h'}) & = & \ell \sum_{r=2\atop r\neq j,j+1}^{\ell} |\lambda^{(r)}|(H_1+ \cdots + H_{r-1}) + \ell|\lambda^{(j+1)}|(H_1 + \cdots + H_{j-1}+H_j) + \\ && +\ell |\lambda^{(j)}|(H_1+ \cdots + (H_{j-1} + H_j)-H_j) - \ell \left( \frac{n(n-1)}{2} + \sum_{r=1}^{\ell} n(\lambda^{(r)}) - n({}^t\lambda^{(r)}) \right) h \\ & =& c_{\bf h}({\blambda}). \end{eqnarray*}
\end{proof}

\subsection{}\label{portmanteau} Since the action of $\tilde{\mathfrak{S}}_{\ell}$ is transitive on alcoves and hence reaches all G.I.T. chambers this completes the calculation of the geometric ordering. Let us present in a relatively succinct way. 

Let $\mathcal{C}_{n,\ell}$ be the set of G.I.T. chambers. We define a surjective map $$\alpha: \Z_{0}^{\ell} \times \mathfrak{S}_{\ell} \times \{\pm \} \longrightarrow \mathcal{C}_{n,\ell}$$ by sending $({\bf s}, w, \pm)$ to the chamber that contains $w^{-1}({\bf 1} + (s_1 - s_{\ell}, s_2-s_1, \ldots , s_{\ell} - s_{\ell -1}))$ in the $+$ case, and to the chamber that contains $w^{-1}(-{\bf 1} + (s_{\ell} - s_1, s_{\ell-1} - s_{\ell}, \ldots , s_1 - s_2))$ in the $-$ case. Recall from \ref{multpart} the mapping $\tau_{\bf s}$ and from \ref{postoneg} the automorphism of $\mathcal{P}(\ell,n)$ which sends $\blambda$ to $\overline{\blambda}$.

\begin{thm} Let ${\bf h}\in \Hr$ and let $\btheta$ be the corresponding stability condition. Let $\blambda, \bmu \in \mathcal{P}(\ell ,n)$. \begin{itemize}
\item[(i)] If $\btheta \in {\alpha^{-1}( {\bf s}, w, +)}$ then the geometric ordering is given by $$\bmu \prec_{\bf h} \blambda \text{ if and only if } \tau_{\bf s}({}^t (w\cdot \blambda)) \lhd \tau_{\bf s}({}^t (w\cdot \bmu)).$$
\item[(ii)] If $\btheta \in {\alpha^{-1}( {\bf s}, w, -)}$ then the geometric ordering is given by $$\bmu \prec_{\bf h} \blambda \text{ if and only if } \tau_{\bf s}({}^t (w\cdot \overline{\blambda})) \lhd \tau_{\bf s}({}^t (w\cdot\overline{\bmu})).$$   
\end{itemize}
\end{thm}

\subsection{Example: the asymptotic case} \label{asy}
 We tie the description of the geometric ordering to the calculation in \cite[Proposition 6.4]{rou}.
\begin{cor} Let $\ell >1$. Consider the G.I.T. chamber defined by $h<0 $ and $H_i \geq (1-n)h$ for $1\leq i \leq \ell -1$. The geometric ordering on $\mathcal{P}(\ell ,n)$ is the dominance ordering.
\end{cor}
\begin{proof} Throughout $\blambda, \bmu\in \mathcal{P}(\ell ,n)$. By Lemma \ref{GITcham} the stated inequalities do indeed give a G.I.T. chamber. To describe the geometric ordering here we take $h=-1$ and we recall from Lemma \ref{basic} that we should find ${\bf s}\in \Z_0^{\ell}$ such that $s_{i+1}-s_i + \frac{1}{\ell} = H_i$. In particular $s_{i+1}-s_i +\frac{1}{\ell} \geq n-1$ and since ${\bf s}$ is integral we even have \begin{equation} \label{dom} s_{i+1}-s_i \geq n-1.\end{equation}

By \ref{multpart} the $\beta$-numbers appearing in $\tau_{\bf s}(\blambda)$ are of the form $\ell  \lambda_j^{(i)} + \ell s_i - \ell j + i$ for $j\geq 1$ and $1\leq i \leq \ell$. In comparing $\beta$-numbers we are only interested in values of $j$ between $1$ and $n$ since $\blambda$ has degree $n$. Then for entries coming from $\lambda^{(i)}$ the maximum value is $\ell n + \ell s_i - \ell + i$. By \eqref{dom} this value is exceeded by any $\beta$-number coming from $\lambda^{(j)}$ with $j>i$  if $\lambda^{(j)}\neq \emptyset$.

Let $\tau_{\bf s}(\blambda)$ have degree $N$. We have ${}^t \tau_{\bf s}({}^t\blambda) \lhd {}^t \tau_{\bf s}({}^t\bmu)$ if and only if $$\sum_{j=t}^N \tau_{\bf s}({}^t\blambda)_j \leq \sum_{j=t}^N\tau_{\bf s}({}^t\bmu)_j$$ for all $1\leq t\leq N$. This is equivalent to the same inequalities holding for the $\beta$-numbers of $ \tau_{\bf s}({}^t\blambda)$ and $ \tau_{\bf s}({}^t\bmu)$. Thus, by the above paragraph, this is equivalent to $$\sum_{j=t}^{n} (\ell ( {}^t\lambda_j^{(i)}) + \ell s_i - \ell j + i )+ \sum_{k=1}^{i-1} \sum_{j=1}^n (\ell ( {}^t\lambda_j^{(k)}) + \ell s_k - \ell j + k )\leq  \sum_{j=t}^{n} (\ell ( {}^t\mu_j^{(i)}) + \ell s_i - \ell j + i )+ \sum_{k=1}^{i-1} \sum_{j=1}^n (\ell  ({}^t\mu_j^{(k)} )+ \ell s_k - \ell j + k)$$ for all $1\leq i \leq \ell$ and $1\leq t\leq n$. Now this is equivalent to $$\sum_{j=t}^{n} {}^t\lambda_j^{(i)} + \sum_{k=1}^{i-1} |\lambda^{(k)}| \leq \sum_{j=t}^{n} {}^t\mu_j^{(i)} + \sum_{k=1}^{i-1} |\mu^{(k)}|.$$ But this is equivalent to the rule for the dominance ordering on $\mathcal{P}(\ell ,n)$. The corollary follows by Proposition \ref{bij}.
\end{proof}

\section{Extension to the facets} \label{facex} In this section we show how to extend the partial order on $\mathcal{P}(\ell ,n)$ from $\Hr$ to ${\bf H}$. This will be useful in the context of Ariki-Koike algebras and the representation theory of $H_{0,{\bf h}}$, see \cite{GM}.

Throughout we will assume that ${\bf h}$ and $\btheta$ are related by the usual rule $\btheta = (-h+H_0, H_1, \ldots , H_{\ell-1}).$

\subsection{The type of a facet}
By Lemma \ref{postoneg} we can assume that $h=-1$ without loss of generality. As we have seen in \ref{stabtoweights} the G.I.T. chambers in the space $(-1, H_1, \ldots ,H_{\ell-1})$ are essentially alcoves for the action of $\tilde{\mathfrak{S}}_{\ell}$ introduced in \eqref{symaction}. Following \ref{portmanteau} we label these alcoves by $({\bf s}, w, +)\in \Z^{\ell}_0\times \mathfrak{S}_{\ell}\times \{\pm\}$.

We now consider the closures of the G.I.T. chambers. By the above these are described as facets of alcoves. Each of these facets has a \label{type-defn}{\it type} $J\subseteq \{0, \ldots , \ell -1\}$ which we now describe.

First consider the closure of the fundamental alcove $\overline{A}_0$, i.e. the closure of the alcove containing $\btheta = {\bf 1}$. This is a fundamental domain for the action of $\tilde{\mathfrak{S}}_{\ell}$. The stabiliser of a point ${\bf h} \in \overline{A}_0$ is a standard parabolic subgroup of $\tilde{\mathfrak{S}}_{\ell}$ generated by simple reflections $\{ \sigma_j: j\in J\}$ for some subset $J\subseteq \{0, \ldots , \ell - 1\}$. This subset is the type of ${\bf h}$, or of the facet contained ${\bf h}$. For the general case, let ${\bf h}$ be $\tilde{\mathfrak{S}}_{\ell}$-conjugate to ${\bf h}'\in \overline{A}_0$. We define the type of ${\bf h}$ just to be the type of ${\bf h}'$.

Note that the alcoves consist of the points whose type is
$\emptyset$.

\subsection{$J$-hearts and $J$-classes}
Let $J\subseteq \{ 0, \ldots , \ell -1\}$ and $\nu \in \mathcal{P}(K)$ for some $K>0$. Recall that a box $(p,q)$ of $\nu$ is said to be $i$-removable for some $0\leq i \leq \ell-1$ if $\res(p,q)$ is congruent to $i$ modulo $\ell$ and if $\nu \setminus \{(p,q)\}$ is the Young diagram of another partition, a predecessor of $\nu$. We then define the \label{heart-defn}{\it $J$-heart} of $\nu$ to be the sub-partition of $\nu$ which is obtained by removing as often as possible $j$-removable boxes with $j\in J$ from $\nu$ and its predecessors. We denote this by $\nu_J$. For instance if $J = \emptyset$ then the $J$-heart of $\nu$ is just $\nu$, whilst if $J = \{ 0, \ldots , \ell-1\}$ then $J$-heart of $\nu$ is $\emptyset$. A subset of $\mathcal{P}(K)$ whose elements are the partitions with a given $J$-heart is called a \label{Jclass-defn}{\it $J$-class}.

\subsection{The partial order}
Let $J\subseteq \{ 0, \ldots , \ell -1\}$. We define a partial order on $\mathcal{P}(K)$ which depends on $J$. Let $\nu ,\mu \in \mathcal{P}(K)$. We take \label{Jord-defn}$\mu \lhd_J \nu$ to be the transitive closure of the relation generated by the rule $\mu$ and $\nu$ have the same $\ell$-core, $\mu \lhd \nu$ and $\mu_J \neq \nu_J$. (But see \ref{hopepart}.)

Now suppose that ${\bf h}\in {\bf H}$ belongs to the closure of
the alcove associated to $({\bf s}, w, +)$ in the notation of \ref{portmanteau}. Then for $\blambda, \bmu
\in \mathcal{P}(\ell ,n)$ we define $\blambda \prec_{\bf h} \bmu$
if $\tau_{\bf s}({}^t(w\cdot \bmu)) \lhd_J \tau_{\bf
s}({}^t( w\cdot \blambda))$.

\begin{rem} This definition is ambiguous as stated. A point ${\bf h}$ may lie on
in the closure of more than alcove and in that case there will be
several definitions of the ordering $\blambda \prec_{\bf h} \bmu$
depending on the alcoves we choose. We expect, but have failed to
prove, that these are all the same so that the definition depends
only on $\bf h$ and not the choice of alcove. We could, however,
remove the ambiguity by insisting that $\bf h$ belongs to the
upper closure of an alcove, as defined in \cite[II.6.2]{janrep}.
In the meantime note that the proof of Part (ii) of the
Proposition below holds independently of the choice of alcove.
\end{rem}

\begin{prop}\label{gotofacs} Let $J\subseteq \{0, \ldots , \ell -1\}$ and suppose that ${\bf h}\in {\bf H}$ belongs to the closure of the G.I.T. chamber corresponding to an alcove $({\bf s}, w, +)$ of type $J$.
\begin{enumerate}
\item[(i)]  The $\C^*$-fixed points of $\mathcal{M}_{\btheta}(n)$ are naturally labelled by the $J$-classes in $\tau_{\bf s}(\mathcal{P}(\ell ,n))$.
\item[(ii)] For any $\blambda ,\bmu \in \mathcal{P}(\ell ,n)$, if $\blambda \preceq_{\bf h} \bmu$ then $c_{\bf h}(\blambda) \geq c_{\bf h}(\bmu)$.
\end{enumerate}
\end{prop}
\begin{proof}
To begin with assume that $w=\id$. Let $\nu$ be the $\ell$-core corresponding to ${\bf s}$ and let $\gamma$ be dimension vector described in Lemma \ref{basic}.

Since $w=\id$ we are studying the closure of an alcove which is a translation of the fundamental alcove $A_0$ and so the parameter $\btheta$ corresponding to ${\bf h}$ has the form $$\btheta = {\bf 1} + (s_1-s_{\ell}, s_2 -s_1, \ldots , s_{\ell} - s_{\ell -1}) + \bep$$ where $\bep \in \Q_{0}^{\ell}$ has each entry $\ep_i$ satisfying $-\frac{1}{\ell} \leq \ep_i \leq \frac{\ell -1}{\ell}$. Set $\bpsi = \bf 1+ \bep$ and note that $\psi_i \geq 0$ for all $i$. By construction we have $\{ j: \psi_j = 0\} = J$.

By applying (a generalisation of) Lemma \ref{reffun} we see that the $\C^*$-fixed points of $\mathcal{M}_{\btheta}(n)$ are in bijection with the $\C^*$-fixed points of $\mathcal{M}_{\bf 1 + \bep}(\gamma)$. The proof of Proposition \ref{bij} shows that the fixed points of $\mathcal{M}_{\bpsi}(\gamma)$ are labelled by $\rho \in \mathcal{P}_{\nu}(K)$ and given by the representations corresponding to $\C[A,B]/I_{\rho}$ as described in that proof. However, these representations need not be stable any longer and thus may give rise to the same polystable representation. By Lemma \ref{tohilb} we can describe these in terms of $-\bpsi$-polystable representations of $\tilde{H}(K)$. The socle of $\C[A,B]/I_{\mu}$ is the representation given by all the removable boxes of the Young diagram of $\rho$. However, the only removable boxes vertices which give $-\bpsi$-stable representations are the $j$-removable boxes for $j\in J$. Factoring out such vertices repeatedly gives rise to the polystable representation associated to $\rho$. But this process is just the passage from $\rho$ to $\rho_J$. Moreover, since partitions $\rho, \mu \in \mathcal{P}_{\nu}(K)$ have the same set of residues modulo $\ell$ we see that $\rho$ and $\mu$ give the same polystable representation if and only if they have the same $J$-heart and same number of boxes of content congruent to $j$ modulo $\ell$ removed for each $j\in J$ if and only if they have the same $J$-heart. This proves Part (i).

For the case $w\neq \id$ we apply the reflection functors of
\cite[3(i)]{Nak3}, which are bijective by \cite[Theorem
3.4]{Nak3}, to set up a correspondence between the fixed points
for the case $({\bf s},\id)$ described above and the case $({\bf
s}, w)$. Since the $J$ does not depend on $w$, this proves Part
(i) in general.

For Part (ii) observe first that by definition if $\blambda
\preceq_{\bf h} \bmu$ then $\tau_{\bf s}({}^t (w\cdot \bmu))
\unlhd \tau_{\bf s}({}^t(w\cdot \blambda))$ and so $\blambda
\preceq_{{\bf h}'} \bmu$ for any ${\bf h}'$ in the chamber
corresponding to ${\bf s}$. Thus $c_{{\bf h}'}(\blambda) \geq
c_{{\bf h}'}  (\bmu)$ by Lemma \ref{geomorder}(iii) and
Proposition \ref{bij}(ii). Since the fixed points vary
continuously in the parameters ${\bf h}$ by Corollary \ref{nat} it
follows that $c_{{\bf h}'}(\blambda) \geq c_{{\bf h}'}  (\bmu)$
for any value of ${\bf h'}$ in the closure of the chamber too, and
in particular for ${\bf h}$. \end{proof}

\subsection{Remarks}\label{hopepart} (i) We hope that it is unnecessary to take the transitive closure to define $\prec_{\bf h}$, i.e. that the displayed relation is already transitive. For example, if we knew that $\mathcal{M}_{\btheta}(n)$ was a normal variety then $\prec_{\bf h}$ would automatically be transitive. Indeed by \cite{sum} there is a $\C^*$-equivariant locally closed embedding $\iota : \mathcal{M}_{\btheta}(n) \longrightarrow \mathbb{P}^N$ for some $N$, where $\mathbb{P}^N$ has an action of the form $$t\cdot (x_0 : x_1 : \cdots :x_N) = (t^{w_0}x_0: t^{w_1}x_1 : \cdots : t^{w_N}x_N).$$ The argument of \cite[Lemma 1]{wl} now shows that if $\tau\lhd \mu\lhd \nu$ have the same $\ell$-core and $\tau_j \neq \mu_J$ and $\mu_J \neq \nu_J$ then $\iota(x_{\btheta}(\tau)),
 \iota(x_{\btheta}(\mu))$ and $\iota(x_{\btheta}(\nu))$ are three distinct points of $\mathbb{P}^N$. In particular, $x_{\btheta}(\tau) \neq x_{\btheta}(\nu)$ and so $\tau_J \neq \nu_J$ by Part (i) of the lemma above. It follows that $\prec_{\bf h}$ is transitive, as claimed.

(ii) We also hope that Part (ii) of the Proposition can be
strengthened to the statement that if $\blambda \prec_{\bf h}
\bmu$ then $\lambda <_{\bf h} \bmu$.

\section{The $a$-function and connections to Hecke algebras}
Our results appear to be related to current work on the two-sided cells of the Hecke algebra of type $B_n$ for unequal parameters, and more generally to Ariki-Koike algebras.
\subsection{Jacon's $a$-function}
Let $\bf h \in {\bf H}$ with $h>0$. Given $\blambda \in
\mathcal{P}(\ell,n)$ we set $$B_u^{(i)}({\bf h}, \blambda) = B_u^{(i)} = h(n+
\lambda^{(i)}_u - u) + (H_1 + \cdots + H_{i-1})$$ for $1\leq i
\leq \ell$ and $u>0$. We define $$a'_{\bf h}(\blambda) = \sum_{1\leq i \leq j \leq \ell \atop { 1\leq u,v\leq n \atop u<v
\text{ if } i=j}} \!\!\!\!\!\min \{ B_u^{(i)}, B_v^{(j)} \} -
\sum_{1\leq i, j \leq \ell \atop {1\leq u\leq n \atop 1\leq k \leq
\lambda^{(i)}_u}} \!\!\!\!\!\!\min\{h(n+k-u)+ (H_1 + \cdots
H_{i-1}), (H_1 + \cdots H_{j-1})\}.$$ Now we let \label{a-defn}$$a_{\bf
h}(\blambda) = a'_{\bf h}(\blambda) - a'_{\bf h}(((n), \emptyset,
\ldots , \emptyset)).$$ This function was studied in \cite[Section
3]{Jac} where it is used to parametrise simple modules for the
Ariki-Koike algebras. For that we need that
$$h=\frac{1}{e}, \qquad H_i = \frac{t_{i+1}-t_i}{e} -
\frac{1}{\ell}$$ for ${\bf t}\in \Z^{\ell}$, see \cite[6.5]{rou}.

\subsection{} There is another action of $\C^*$ on $\mathcal{M}_{\btheta}({\bf d})$ which does not pass over to an action on $\mathcal{X}_{\frac{1}{2}\btheta}(n)$. It is induced from the action on $R({\bf d}')$ given by
\begin{equation}
\label{2ndact}
\lambda \circ ({\bf X}, {\bf Y}; v,w) = ( {\bf X}, \lambda{\bf Y};  v, \lambda w).
\end{equation}
In other words, this is essentially associated to ``half" of the
action that we studied earlier. For any $\btheta \in \Theta_1$ it gives rise to an analogue of
$f_{\btheta}$, namely $A_{\btheta} : \mathcal{M}_{\btheta}(n)
\stackrel{r_{\btheta}}{\longrightarrow} \mathcal{M}_{\bpsi}(\tau_{\beta^{\vee}}\ast
n\delta) \longrightarrow \R$ which is defined by
$$A_{\btheta}(r_{\btheta}^{-1}({\bf X}, {\bf Y}; v,w)) = \sum_{r=0}^{\ell -1} \Tr
(Y_rY_r^{\dagger}),$$ where $r_{\btheta}$ is the composition of reflection functors producing an isomorphism between $\mathcal{M}_{\btheta}(n)$ and $\mathcal{M}_{\bpsi}(\tau_{\beta^{\vee}}\ast n\delta)$ discussed in Sections \ref{combgeom} and \ref{facex}. We extend this function to $\Theta_{>0} = \{ \btheta = (\theta_0, \ldots , \theta_{\ell -1}) \in \Q^{\ell} : \theta_0 + \cdots + \theta_{\ell-1} >0 \}$ by just follow the above geometric description, but allowing $\btheta \in \Theta_{>0}$ instead of just $\Theta_1$.  To define it on $\Theta_{<0}$ we follow \ref{postoneg} and define (in the notation of the proof of Lemma \ref{postoneg}) $$A_{\btheta}( {\bf X}', {\bf Y}'; v' , w') = A_{\overline{\btheta}}(\phi( {\bf X}', {\bf Y}'; v' , w'))$$ for $\btheta \in \Theta_{< 0}$ and $( {\bf X}', {\bf Y}'; v' , w') \in \mathcal{M}_{\btheta}(n)$. 
\begin{prop} \label{afnrel} Let ${\bf h}\in \Hr$ with $h\neq 0$ and let $\btheta$ be the corresponding stability condition. Let $\blambda \in \mathcal{P}(\ell ,n)$.
\begin{enumerate}
\item[(i)] The function $A_{\btheta}(\blambda)\equiv  A_{\btheta}(x_{\btheta}(\blambda)): \Theta_{\neq 0} \longrightarrow \Q$
is piecewise-linear on the parameter space $\Theta_{\neq 0} = \{\btheta \in \Q^{\ell} : \theta_0 + \cdots + \theta_{\ell -1} \neq 0\}$. In fact, it is linear on the closure of the positive cones of the alcoves. 
\item[(ii)] Let $h>0$. Then $$a_{\bf h} (\blambda)
\doteq A_{\btheta}(\blambda) \qquad \text{for all }\blambda \in
\mathcal{P}(\ell ,n).$$
\end{enumerate}
\end{prop}
\begin{proof}
Part (i) is straightforward. For $\btheta \in \Theta_{>0}$ Lemma \ref{tohilb} and the proof of Proposition \ref{bij} show that \begin{equation} \label{whatAis} A_{\btheta} (\blambda) \doteq \sum_{(p,q) \in \tau_{\bf s}({}^t \blambda)} \psi_{\res (p,q)} q.\end{equation} This is linear on the closure of the positive cone of the alcove containing $\btheta$. For $\btheta \in \Theta_{< 0}$ we then have $A_{c\btheta}(\blambda) = A_{\overline{c\btheta}}(\phi( x_{c\btheta}(\blambda))) = A_{c\overline{\btheta}}(x_{c\overline{\btheta}}(\overline{\blambda})) = c A_{\overline{\btheta}}(x_{\overline{\btheta}}(\overline{\blambda})) = cA_{\btheta} (\blambda),$ proving linearity there.

For Part (ii) we deal first with the alcoves that are translations of the fundamental alcove. For this note that $\btheta$ has the form $\tau_{-\beta^{\vee}}\cdot {\bf 1} + \bep$ where $\bep \in \Q_0^{\ell}$ has each entry bounded $-\frac{1}{\ell} \leq \ep_j \leq \frac{\ell -1}{\ell}$. Let $\gamma = \tau_{\beta^{\vee}}\ast n\delta$ and set $\bpsi = \bf 1 + \bep$. Finally let ${\bf s}\in \Z^{\ell}_0$ be defined as in Lemma \ref{lcores}.
 
We express the $A$-function in terms of the $\beta$-numbers of $\blambda$. Define $F: \NN \longrightarrow \R$ by the rule $F(i) = \psi_0 + \psi_{1} + \cdots + \psi_{i-1}.$ Observe that $F$ is increasing since $F(i+1) - F(i) = \psi_{i} = \frac{1}{\ell} + \ep_{i} \geq 0$. Then, given $\nu\in \mathcal{P}(K)$ and $N=\ell k$ large enough so that $\nu_{N+1}=0$, we can write \begin{eqnarray*} \sum_{(p,q)\in \nu} \psi_{\res(p,q)}q &=& \sum_{i=1}^N (i-1)(\psi_{-i} + \psi_{-i+1} + \cdots + \psi_{-i-1+ \nu_i}) \\ & = & \sum_{i=1}^N (i-1) (F(N+\nu_i - i +1) - F(N+1-i)).\end{eqnarray*}  Let $\beta_1 > \beta_2 > \cdots$ be the $\beta$-numbers of $\tau_{\bf s}({}^t\blambda)$ (with $s= \sum s_j = 0$). It follows from \eqref{whatAis} and the above that \begin{eqnarray*} A_{\btheta}(\blambda) &=&  \sum_{(p,q) \in \tau_{\bf s}({}^t\blambda)} \psi_{\res (p,q)} q \\ & = & \sum_{i=1}^N (i-1) (F(N+ \beta_i) - F(N+1-i))  \doteq  \sum_{i=1}^N (i-1) F(N+\beta_i) \\ &=& \sum_{1\leq i < j \leq N} \min\{ F(N+\beta_i), F(N+\beta_j)\}.\end{eqnarray*}The last equality holds because $F$ is an increasing function.

We let $N = \ell (n -S)$ where $S$ is the smallest entry of ${\bf s}$. This is the least possible value allowed for $N$, i.e. that insures that $N+ \beta_{N+i} = -i$ for all $i>0$. The definition of $\tau_{{\bf s}}(\blambda)$ from \ref{multpart} now gives
\begin{eqnarray} \label{whatAreallyis} A_{\btheta}(\blambda) &=& \!\!\!\!\!\!\!\!\!\!\sum_{1\leq i\leq j \leq \ell \atop {{1\leq u \leq n + s_{ i} -S} \atop {1\leq v \leq n+s_{j} -S \atop u<v \text{ if } i=j}}} \!\!\!\!\!\!\!\!\!\!\! \min\{F(N + \ell(({}^t\lambda^{(i)})_u + s_{i}-u) + i), F(N+ \ell(({}^t\lambda^{(j)})_v + s_{j} -v ) + j)\} \notag \\ & = & \!\!\!\!\!\!\!\!\!\!\sum_{1\leq i\leq j \leq \ell \atop {{1\leq u \leq n + s_{i} -S} \atop {1\leq v \leq n+s_{j} -S \atop u<v \text{ if } i=j}}} \!\!\!\!\!\!\!\!\!\!\! \min\{n+ ({}^t\lambda^{(i)})_u + s_{ i} -S  -u + F(i), n+ ({}^t\lambda^{(j)})_v + s_{j} -S  -v + F(j)\}. \end{eqnarray}
By the definition of $\bep$ we have $$F(i) +s_i = \frac{i}{\ell} + \ep_0 + \cdots + \ep_{i-1} +s_i = H_0 + \cdots + H_{i-1} +1 + s_{\ell}.$$ Thus we find that \begin{eqnarray*} A_{\btheta}(\blambda) &\doteq& \!\!\!\!\!\!\!\!\!\sum_{1\leq i < j \leq \ell \atop {{1\leq u \leq n + s_{ i} -S} \atop {1\leq v \leq n+s_{j} -S}}} \!\!\!\!\!\!\!\!\!\min \{(n-S+({}^t\lambda^{(i)})_u - u + H_0 + \cdots +  H_{i-1}, (n-S+({}^t\lambda^{(j)})_v-v + H_0 + \cdots + H_{j-1}\} \\ & +& \!\!\!\!\!\!\!\!\!\sum_{1\leq i \leq \ell \atop  1\leq u< v\leq n+s_i-S} \!\!\!\!\!\!\!\!\!\min \{ n-S+({}^t\lambda^{(i)})_u - u + H_0 + \cdots +  H_{i-1}, n-S+({}^t\lambda^{(i)})_v - v + H_0 + \cdots +  H_{i-1}\}. \end{eqnarray*}

Now we begin the comparison with $a_{\bf h}(\blambda)$. Here we
have ${\bf h} = (1, H_1, \ldots , H_{\ell-1})$ and so we have $A_{\btheta}(\blambda) = A_{\overline{\btheta}}(\overline{\blambda}).$ It follows from the definition of
$\overline{\blambda}$ that the values we must compare in
$A_{\btheta}(\blambda)$ are of the form $n - S + \lambda^{(\ell +
1 -i)}_u - H_0 - \cdots - H_{i-1} = n-S + \lambda^{(\ell + 1 -
i)}_u -u + H_1 + \cdots + H_{\ell - i}$. This equals $B^{(\ell + 1
- i)}_u- S$ and so we are left to prove that \begin{eqnarray*}
&&\sum_{1\leq i \leq j \leq \ell \atop { 1\leq u,v\leq n \atop u<v
\text{ if } i=j}} \!\!\!\!\!\min \{ B_u^{(i)}, B_v^{(j)} \} -
\sum_{1\leq i, j \leq \ell \atop {1\leq u\leq n \atop 1\leq k \leq
\lambda^{(i)}_u}} \!\!\!\!\!\!\min\{n+k-u+ (H_1 + \cdots H_{i-1}),
(H_1 + \cdots H_{j-1})\} \\&& \qquad \qquad \doteq
\!\!\!\!\!\!\!\!\!\sum_{1\leq i < j \leq \ell \atop {{1\leq u \leq
n + s_{ i} -S} \atop {1\leq v \leq n+s_{j} -S}}}
\!\!\!\!\!\!\!\!\!\min \{ B^{(i)}_u, B^{(j)}_v\} + \sum_{1\leq i
\leq \ell \atop  1\leq u< v\leq n+s_i-S} \!\!\!\!\!\!\!\!\!\min \{
B^{(i)}_u, B^{(i)}_v\}. \end{eqnarray*} We can cancel the common
terms in the first sum. We can also remove all the terms from the
bottom half of the equality that involve values $u$ and $v$
between $n+1$ and $n+s_i-S$ since these depend only on ${\bf s}$
and not on $\blambda$. Thus we are left to show that
\begin{eqnarray} \label{wanta}&& - \sum_{1\leq i, j \leq \ell
\atop {1\leq u\leq n \atop 1\leq k \leq \lambda^{(i)}_n}}
\!\!\!\!\!\!\min\{n+k-u+ (H_1 + \cdots H_{i-1}), (H_1 + \cdots
H_{j-1})\} \\ &&\qquad \qquad \doteq \sum_{1\leq i \neq j \leq
\ell \atop {1\leq u \leq n \atop n+1\leq v \leq n+s_j - S }}\min
\{ B_u^{(i)}, B_v^{(j)} \} + \sum_{1\leq i \leq \ell  \atop {1\leq
u \leq n \atop n+1\leq v \leq n+s_i - S}} \min \{ B_u^{(i)},
B_v^{(i)} \}. \notag\end{eqnarray} This is a simple calculation,
but it is a little involved; you can find it in the second
appendix.

To finish the proof of (ii) we have to deal with $w$-translates of alcoves where $w\in \mathfrak{S}_{\ell}$. By definition we have that $A_{w\cdot \btheta}(\blambda) = A_{\btheta}(w^{-1}\cdot \blambda)$ and so it is enough to prove the analogous equality for the $a$-function. We do this for the generators $\sigma_i\in \mathfrak{S}_{\ell}$ with $1\leq i \leq \ell -1$. 

Given $1\leq j \leq \ell$ let $H_{< j} = H_1 + \cdots + H_{j-1}$. Then the action of $\sigma_i$ defined by \eqref{symaction} transfers to the action on ${\bf H}$ that fixes $h$ and sends $H_{<j}$ to $H_{< \sigma_i(j)}$ if $i>1$ and sends $H_{< j}$ to $H_{<\sigma_i(j)} - H_1$ if $i=1$. Thus if $i>1$ we have that $B_u^{(i)}(\sigma_i \cdot {\bf h} , \blambda) = B_u^{(i)}({\bf h}, \sigma_i \cdot \blambda)$ and so it follows that $a_{\sigma_i \cdot \bf h}(\blambda) = a_{\bf h}(\sigma_i \cdot \blambda)$. Similarly if $i=1$ then we find that $a'_{\sigma_1 \cdot \bf h}(\blambda) \doteq a'_{\bf h}(\sigma_1 \cdot \blambda)$ for all $\blambda \in \mathcal{P}(\ell ,n)$ where the constant difference between these is a multiple of $H_1$. Since the $a$-function is defined as the difference of two $a'$-functions we we find in this case too that $a_{\sigma_1 \cdot \bf h}(\blambda) = a_{\bf h}(\sigma_1 \cdot \blambda)$. Thus we have $$a_{w\cdot \bf h}(\blambda) = a_{\bf h}(w^{-1}\cdot \blambda)$$ for all $w\in \mathfrak{S}_{\ell}$, as required.
\end{proof}
\subsection{} There is an ordering on $\mathcal{P}(\ell ,n)$ given by comparison of $a$-function on $\ell$-multipartitions which is introduced in \cite[Definition 4.3]{Jac}. We now show that this refines the partial order $\prec_{\bf h}$.
\begin{thm} \label{aorderthm} Let ${\bf h}\in {\Hr}$ and $\blambda ,\bmu \in \mathcal{P}(\ell ,n)$. Then $\bmu\prec_{\bf h} \blambda$ implies that $a_{{\bf h}}(\bmu) < a_{\bf h}(\blambda)$.
\end{thm}
\begin{proof}
There is a $T^2:= (\C^*)^2$-action on
$\mathcal{M}_{\bpsi}(\gamma)$ and analogously on $\hil{\nu}$ given by
$$(s,t)\cdot ({\bf X}, {\bf Y}; v,w) (s{\bf X}, t{\bf Y}; v, w)
\qquad \text{and} \qquad (s,t)\cdot ({ X}, {Y}; v,w) = (s{X},
t{Y}; v, w).$$ The original action we studied in
\ref{circleaction} is a specialisation of this to $s=t^{-1}$. Thus
the fixed points of $\mathcal{M}_{\bpsi}(\gamma)$ under this
$T^2$-action are fixed by the original $\C^*$-action; conversely
the monomial description in Proposition \ref{bij} of the
$\C^*$-fixed points shows that these fixed points are fixed by
$T^2$. Thus the $T^2$-fixed points are labelled by
$\mathcal{P}(\ell,n)$.

Let $p,q$ be positive integers and consider the subtorus $T^{p,q}:= \{ (t^p, t^{-q}): t\in \C^*\}$ of $T^2$. Note that $T^{1,1}$ is the one-dimensional torus we used earlier: we will call it the standard $\C^*$. We claim that the attracting sets on $\mathcal{M}_{\bpsi}(\gamma)$ of the $T^{p,q}$-action and the standard $\C^*$-action are the same. To prove this we begin by noting that since the $\C^*$-fixed points of $\mathcal{M}_{\bpsi}(\gamma)$ are fixed by $T^2$, the $\C^*$-attracting sets are $T^2$-stable. In particular they are $T^{p,q}$-stable. Thus, by the uniqueness claim of \cite[Theorem, p.492]{bial}, it is enough to prove that for any $\C^*$-fixed point $z$ we get the same decomposition of $T_z\mathcal{M}_{\bpsi}(\gamma)$ into positive and negative eigenspaces with respect to the $T^{p,q}$-action and the standard $\C^*$-action.

By \cite[(2.15) and Proof of Theorem 3.2]{haiqt} the eigenvalues of $T^2$ on the tangent space of $\hi{K}$ at the fixed point $I_{\lambda}$ ($\lambda \in \mathcal{P}(K)$) are given by the $2K$ monomials $\{ s^{1+l(x)}t^{-a(x)},  s^{-l(x)}t^{1+a(x)} : x\in \lambda\},$ where $a(x)$ and $l(x)$ are the arm and leg of the cell $x$ in the Young diagram of $\lambda$. The definition of $a(x)$ and $l(x)$ is given in \cite{haiqt}; what is vital here is that they are combinatorial quantities which are always non-negative. The eigenvalues of $T^{p,q}$ are given by $\{ t^{p(1+l(x)) +qa(x)}, t^{-pl(x)-q(1+a(x))} \}$. Since $l(x)$ and $a(x)$ are non-negative for all $x\in \lambda$ we see that the lines in the $(p,q)$-plane which produce zero eigenspaces in $T_z\mathcal{M}_{\bpsi}(\gamma)$ do not have positive gradient. Since we are assuming that both $p$ and $q$ are positive, it follows that the decomposition of $T_z\mathcal{M}_{\bpsi}(\gamma)$ is independent of the choice of $p$ and $q$, as required. This completes the proof of the claim.

Let $\ep =p/q$ with $p,q$ positive integers. Consider the function
$$A^{\ep}_{\btheta}({\bf X},{\bf Y}; v, w) = \ep\sum_{r=0}^{\ell
-1} \Tr(X_rX_r^{\dagger}) - \sum_{r=0}^{\ell -1} \Tr
(Y_rY_r^{\dagger}).$$ The function $qA_{\btheta}^{\ep}$ is the analogue of the function $f_{\btheta}$ of \ref{morse} for the group $T_{p,q}$ instead of the standard $\C^*$-action, see also \cite[5.2]{nakhilb}. Arguing exactly as in the proof Lemma \ref{geomorder} and using the fact the attracting sets for the $T_{p,q}$-action and the standard $\C^*$-action agree, we see that if $\bmu \prec_{\bf h} \blambda$ then $A_{\btheta}^{\ep}(\bmu)  > A_{\btheta}^{\ep}(\blambda)$. Since $\lim_{\ep \rightarrow 0}
A_{\btheta}^{\ep} = - A_{\btheta},$ we deduce that \begin{equation} \label{aineq} A_{\btheta}(\bmu) \leq A_{\btheta}(\blambda).\end{equation}

The inequality \eqref{aineq} holds for all $\btheta$ in the alcove containing $\bf h$. Since $A_{\btheta}(\blambda)$ and $A_{\btheta}(\bmu)$ are linear on this alcove by Proposition \ref{afnrel}(i), it follows that either the inequality is strict or $A_{\btheta}(\blambda) = A_{\btheta}(\bmu)$ for all $\btheta$ in the alcove. It is easy to see that the second possibility cannot occur. By \eqref{whatAis} we have $$A_{\bf 1}(\blambda)  \doteq \sum_{(p,q) \in \tau_{\bf s}({}^t \blambda)}\!\!\!\!\!\! q  \,\, = n(\tau_{\bf s}({}^t \blambda)) \qquad \text{for all }\blambda \in \mathcal{P}(\ell, n).$$ Then a quick calculation, or \cite[Theorem B and Proposition 1.6]{shi}, shows that if $\tau_{\bf s}({}^t \blambda) \lhd \tau_{\bf s}({}^t \bmu)$ then $n(\tau_{\bf s}({}^t \bmu))< n(\tau_{\bf s}({}^t \blambda))$. In other words, thanks to Proposition \ref{bij}(ii), if $\bmu \prec_{\bf h} \blambda$ then $A_{\bf 1}(\bmu) = n(\tau_{\bf s}({}^t \bmu))< n(\tau_{\bf s}({}^t \blambda)) = A_{\bf 1}(\blambda)$, as required.

We have shown that if $\bmu \prec_{\bf h} \blambda$ then $A_{\btheta}(\bmu) < A_{\btheta}(\blambda)$. The theorem follows from Proposition \ref{afnrel}(ii).
\end{proof}
\subsection{}
\label{a2order}
By the same argument as the proof of Proposition \ref{gotofacs}(ii) we see
that Theorem \ref{aorderthm} extends to the walls with the
statement: $$\text{if ${\bf h}\in {\bf H}$ then $\bmu \prec_{\bf
h} \blambda$ implies $a_{\bf h}(\mu) \leq a_{\bf h}(\blambda)$.}$$
It is not true, however, that one of either the $c$-ordering or
the $a$-ordering refines the other.

\subsection{}\label{Brmi} We now have an elementary result which generalises \cite[Section 4.21 and Proposition 4.1]{brmi}. (Note that $n$ (respectively $n(n-1)$) is the number of reflections in $G_n(2)$ associated to the parameter $H_1$ (respectively $h$).)
\begin{cor}
Let $n=2$, ${\bf h} \in \bf H$ and $\blambda \in \mathcal{P}(2 ,n)$. Then $$c_{\bf h}({\blambda}) = a_{{\bf h}}({}^t \blambda) + (nH_1-n(n-1)h - a_{\bf h}(\blambda)).$$
\end{cor}
\begin{proof} To start with we will let $\ell\geq 2$ and set ${\bf h}' = (h, H_{\ell-1}, \ldots , H_{1})$. Note that when $\ell=2$ we have ${\bf h}' = {\bf h}$. A little later we will have to insist $\ell=2$.

Assume that $h=1$. Then we get $$\bpsi = (1-H_0 - s_1 + s_\ell,
-H_1 - s_2 + s_{1} , \ldots , -H_1 - s_{\ell} + s_{\ell -1})$$
where ${\bf s}$ defines the chamber for $\overline{h}$. We begin
with a couple of simple observations. Let $\nu \in \mathcal{P}(K)$
and let $\bpsi \in \Q^{\ell}$. We set $\overline{\bpsi } = (\psi_0, \psi_{\ell -1}, \psi_{\ell -2}, \ldots , \psi_{1} ) $ and
$\overline{\bf s} = (-s_{\ell}, -s_{\ell -1}, \ldots , -s_1)$.
Then the chamber of $(H_0, H_{\ell -1}, \ldots , H_1)$ corresponds
to $\overline{\bf s}$ and then we find that $\overline{\bpsi}$
corresponds to ${\bf h}'$. Moreover ${}^t\tau_{\bf
s}({}^t\blambda) = \tau_{\overline{\bf s}}(\blambda)$ where
$\overline{\bf s} = (-s_{\ell}, -s_{\ell -1}, \ldots ,
-s_1)$.

We have \begin{eqnarray*} \sum_{(p,q)\in \nu} \psi_{\res(p,q)} p &=& \sum_{i=1}^{K} (0\psi_{-i} + 1\psi_{-i+1} + \cdots + (\nu_i -1) \psi_{-i + \nu_i - 1}) \\& = & \sum_{i=1}^{K} (i-1)(\psi_{i} + \psi_{i-1} + \cdots + \psi_{i + ({}^t \nu)_i -1}) \\ & = & -\sum_{i=1}^{K} (i-1)(\overline{\psi}_{-i} + \overline{\psi}_{-i+1} + \cdots + \overline{\psi}_{-i - ({}^t \nu)_i +1}). \end{eqnarray*}

It thus follows from the proof of Proposition \ref{bij} that \begin{eqnarray*}c_{\bf h}({\blambda}) &\doteq& \sum_{(p,q)\in \tau_{\bf s}({}^t\blambda)} \psi_{\res(p,q)}\res(p,q)\\ &\doteq& - a_{\bf h}(\blambda)+\sum_{(p,q)\in \tau_{\bf s}({}^t\blambda)} \psi_{\res(p,q)}p\\ & = & -a_{\bf h}(\blambda) + \sum_{(p,q)\in \tau_{\overline{\bf s}}(\blambda)} \overline{\psi}_{\res(p,q)}q \\ & \doteq& -a_{\bf h}(\blambda) + a_{{\bf h}'}({}^t\blambda).\end{eqnarray*}

Now when $\ell=2$ we get the equality in the statement of the corollary by observing that when $\blambda = ((n), \emptyset, \ldots , \emptyset)$ we have $c_{\bf h}({\blambda}) = 0 = a_{\bf h}(\blambda)$ while, by a calcuation left to the reader, $$a_{{\bf h}'}({}^t\blambda) =  -nH_1 + n(n-1)h.$$
\end{proof}

\subsection{Remarks} For ${\bf h}\in \bf H$ it would be very interesting
 to know whether the ordering given here on ($J$-classes of) $2$-multipartitions
  agrees with the ordering on unequal parameter two-sided cells of the Weyl group of type $B_n$. The chamber picture has already appeared in the original work of \cite{lus1}, and conjectural claims on constancy within chambers were made in \cite[Conjecture 2.17]{geck}. Moreover in the asymptotic case everything agrees thanks to \ref{asy} and \cite[(3.8)]{bon}.
It is shown in \cite{GM} that the $J$-classes agree are in natural bijection with the conjectured combinatorial description of the two-sided
cells in \cite[Conjectures A and B]{BGIL} and related work in
\cite{piet}.

\section{Connections to other topics}

\subsection{Category $\mathcal{O}_{\bf h}$} \label{Oord}
We first ask whether there is a more refined ordering than the $c$-ordering on $\mathcal{O}_{\bf h}$ described in \ref{O}.
\begin{question}
Is $\mathcal{O}_{\bf h}$ is a highest weight category with ordering given by $L_{\bf h}(\blambda) < L_{\bf h}(\bmu)$ if $\blambda \prec_{\bf h} \bmu$ and $c_{\bf h}({\blambda}) - c_{\mu}({\bf h}) \in \Z$?
\end{question}

\noindent
This is true in the asymptotic case of Section \ref{asy}, \cite[Theorems 6.6 and 6.8]{rou}. In particular, it is true for $\ell =1$ and any value of $h$. An important consequence of a positive answer to this question would be the strengthening of Rouquier's Theorem \ref{rouq} where $c$-chambers are replaced by G.I.T. chambers in the statement of the theorem. Thanks to Theorem \ref{aorderthm} and \ref{a2order} it would also imply that $\{ KZ_{\bf h}(L_{\bf h}(\lambda)): KZ_{\bf h}(L_{\bf h}(\lambda))\neq 0\}$ would be a canonical basis set for $\mathcal{H}_{\bf q}(W)$, where $KZ_{\bf h}$ denotes the KZ-functor on $\mathcal{O}_{\bf h}$ of \cite{GGOR}. Moreover if $\ell = 2$ then $\{ KZ_{\bf h}(\Delta_{\bf h}(\blambda)): \blambda \in \mathcal{P}(2,n)\}$ would be the set of cell modules with respect to the cellular algebra structure found on $\mathcal{H}_{\bf q}(W)$ in \cite[Theorem 1.1]{geckcell}.

\subsection{}\label{cycle} In this section we assume that ${\bf h}\in \Hr$ and the $\Z$-algebra $B_{\bf h}$ appearing in Theorem \ref{Zalg} is a Morita $\Z$-algebra, i.e. the ``shift functors" defined in the proof of theorem are Morita equivalences.
Then we can mimic the construction of characteristic cycles in \cite[2.7]{GS2} and thus associate to any finitely generated $H_{1, {\bf h}}$-module $M$ a cycle in $\mathcal{M}_{\btheta}(n)$, written $\ch (M)$. As in \cite[Proposition 4.8]{GS2}, the characteristic cycles of objects from $\mathcal{O}_{\bf h}$ will lie in $\mathcal{Z}_{\btheta}$ as defined in \ref{geomorder}. Let $\Delta_{\bf h}(\blambda)$ be the standard module in $\mathcal{O}_{\bf h}$ with simple head $L_{\bf h}(\blambda)$, \cite[Section 3.2]{GGOR}.
\begin{question}
Is $\ch (\Delta_{\bf h}(\blambda)) = \sum_{\bmu\prec_{\bf h} \blambda} a_{\blambda, \bmu} [\mathcal{Z}_{\bmu}]$ with $a_{\bmu, \blambda} \in \Z_{\geq 0}$ and $a_{\blambda, \blambda} = 1$?
\end{question}

\noindent
This is true when $\ell =1$ by \cite[Theorem 6.7]{GS2} and it has been established for $n=1$ and any $\ell$ in \cite{kuw}. In \ref{genn!} we present a possible symmetric function theoretic interpretation of the $a_{\bmu, \blambda}$.
\begin{lem} Assuming that $B_{\bf h}$ is a Morita $\Z$-algebra then a positive answer to Question \ref{cycle} implies a positive answer to Question \ref{Oord}.
\end{lem}
\begin{proof}
If Question \ref{cycle} has a positive answer then the argument of \cite[Corollary 6.8]{GS2} shows that $\rch (L_{\bf h}(\bmu)) = [\mathcal{Z}_{\bmu}] + \sum_{\bnu} b_{\bmu, \bnu}[\mathcal{Z}_{\bnu}]$ where $\bnu <_{\bf h} \bmu$ ($c$-ordering) and $\rch$ denotes the restricted characteristic cycle defined in \cite[2.8]{GS2}. It then follows from \cite[Lemma 2.8]{GS2} that if $L_{\bf h}(\bmu)$ is a composition factor of $\Delta_{\bf h}(\blambda)$ then $[\mathcal{Z}_{\bmu}]$ appears with non-zero multiplicity in $\rch \Delta_{\bf h}(\blambda)$. Thus $\bmu \preceq_{\bf h} \blambda$. This, together with \cite[2.6.2 and Proposition 3.3]{GGOR} show that $(\mathcal{O}_{\bf h},\prec_{\bf h})$ is a highest weight category.
\end{proof}

\subsection{Generalised $n!$ conjecture} \label{genn!} We can mimic the construction of \cite[5.5]{GS} and \cite[4.2]{GS2} to associate  a coherent sheaf on $\mathcal{M}_{\btheta}(n)$ to any $H_{1,{\bf h}}$-module with a good filtration $(M,\Lambda)$: we label it by $\hat{\Phi}_{\Lambda}(M)$.
\begin{question} Assume that ${\bf h}\in \Hr$ belongs to an alcove labelled by $({\bf s},\id , +) \in \Z_0^{\ell}\times \mathfrak{S}_\ell\times \{ \pm\}$. Let $(M, \Lambda)$ be $H_{1, {\bf h}}$ with its filtration by order of differential operators. Is $\hat{\Phi}_{\Lambda}(H_{1,{\bf h}})$ the $G$-equivariant vector bundle on $\mathcal{M}_{\btheta}(n)$ whose existence is predicted in \cite[Conjecture 7.2.13]{hai}?
\end{question}

\noindent
When $\ell =1$ this question has a positive answer thanks to \cite[Theorem 4.5]{GS2}.

\subsection{} The question above also suggests candidates for the numbers $a_{\bmu,\blambda}$ appearing in Question \ref{cycle} since the standard modules $\Delta_{\bf h}({\blambda})$ are obtained from $H_{1,\bf h}$ by factoring out $Y_1, \ldots , Y_n$. We state it here only for ${\bf h}\in \Hr$ belonging to alcoves labelled by $({\bf s},\id, +)$. The generalisation is not complicated, but notationally awkward to state.
\begin{question} Are the integers $a_{\blambda, \bmu}$ which appear in Question \ref{cycle} the wreath Kostka numbers, that is the coefficient of the irreducible representation ${}^t\blambda$ of $G$ in the expansion of the specialised wreath Hall-Littlewood polynomial $H_{{}^t\bmu}(1)$ of \cite[Very end of Section 7.2]{hai}?
\end{question}
In the case $\ell =1$ this question has a positive answer thanks to \cite[Theorem 6.7]{GS2}. Moreover it is consistent with the conjectural ordering property of wreath Macdonald polynomials presented in \cite[Conjecture 7.2.19, (i)]{hai}.

\subsection{Derived equivalences} Let ${\bf h}\in \Hr$, so that the corresponding variety $\mathcal{M}_{\btheta}(n)$ is a symplectic resolution of $V/G$. By \cite{BK} there are equivalences of bounded derived categories $$\begin{CD} \Phi: D^{b}(\C[V]\ast G) @>\sim >> D^b (\coh \mathcal{M}_{\btheta}(n)).\end{CD}$$ As the $\mathcal{M}_{\btheta}(n)$ are isomorphic in a chamber we may as well pick the $\mathfrak{S}_{\ell}$-translates of our standard representatives $\btheta = {\bf 1} + (s_1-s_{\ell}, s_1-s_1, \ldots , s_{\ell}-s_{\ell -1})$ and consider only the corresponding equivalences $\Phi_{({\bf s},w)}$ where ${\bf s}\in \Z^{\ell}_0$ and $w\in \mathfrak{S}_{\ell}$. Now given a wall between two adjacent alcoves whose representatives are labelled by $({\bf s}, w, +)$ and $({\bf s}', w', +)$ we could ask for a geometrically defined wall-crossing functor $\mathcal{T}$ and then find an endomorphism  $\Phi_{({\bf s}',w')}^{-1} \circ \mathcal{T} \circ  \Phi_{({\bf s},w)}$ of $D^b(\C[V]\ast G)$.
\begin{question} Is there an action of the affine braid group of type $\tilde{A}_{\ell -1}$ on $D^b(\C[V]\ast G)$ which arises from a wall-crossing action of $\tilde{\mathfrak{S}}_{\ell}$ on alcoves?
\end{question}

\subsection{} We hope that such equivalences would be $\C^*$-equivariant and so preserve the attracting sets we have studied. We could then study what happens to the cycles $[\mathcal{Z}_{\blambda}]$ under these derived equivalences. Understanding combinatorially that picture would then give a model for derived equivalences for category $\mathcal{O}$, and hence also for the structure of the corresponding quasi-hereditary algebras (generalised cyclotomic $q$-Schur algebras to be precise!).
\begin{question}[Rouquier] Let ${\bf h}, {\bf h}' \in \Hr$ differ by an element of $\Z^{\ell}$. Are $\mathcal{O}_{\bf h}$ and $\mathcal{O}_{\bf h'}$ derived equivalent (in a way that reflects the geometric equivalences)?
\end{question}

\subsection{Generalised cyclotomic $q$-Schur algebras}
Yvonne has conjectured in \cite[Section 2]{yvo} that the matrix of composition multiplicities $[\Delta_{\bf h}(\blambda) : L_{\bf h}(\bmu)]$ in $\mathcal{O}_{\bf h}$ equals the transition matrix (specialised at $q=1$) between the standard basis and canonical basis of a higher level Fock space whose multicharge ${\bf s}$ depends on $\bf h$. The ordering on multipartitions that Yvonne studies depends not only on ${\bf s}$, but also on a a root of unity $\xi$. If his conjecture is true, this ordering will define an ordering on $\mathcal{O}_{\bf h}$ which will respect the highest weight structure. When $\xi = 1$ (an unusual choice in Yvonne's regime!) his ordering agrees with our geometric ordering by Theorem \ref{portmanteau}. We expect that for any choice of $\xi$ the corresponding ordering is refined by the geometric, or equivalently $\xi=1$, ordering. 

The results here show that the quiver varieties $\mathcal{M}_{\btheta}(n)$ are degenerations of rational Cherednik algebras {\it and} that they contain some of the relevant combinatorial information of the Fock space. The fine structure of the canonical basis of the Fock space is not seen by the geometry; that should be found in the rigidity provided by the quantisation to the world of the (noncommutative geometry of the) rational Cherednik algebra. Quantisations of $\mathcal{M}_{\btheta}(n)$ should hold the key to the combinatorics of Cherednik algebras.

\appendix
\section{Calculation for proof of Theorem \ref{bij}}
We will keep the same notation as in the proof.

\subsection{} We have to show that equation \eqref{hilbequ} is true, i.e. that $$ \sum_{(p,q)\in \tau_{\bf s}({}^t\blambda)} \psi_{\res (p,q)} \res(p,q) \doteq c_{\blambda} ({\bf h})  \qquad \text{for all }\blambda \in \mathcal{P}(\ell, n).$$ Here $\bpsi = \bf 1 + \bep$ and $\bf h$ corresponds to $\btheta + \bep$ where $\bep \in \Q^{\ell}_0$ and by Lemma \ref{basic} $$\btheta = (s_1 - s_{\ell} + \frac{1}{\ell} , \ldots, s_{\ell} - s_{\ell-1} + \frac{1}{\ell}).$$ We split this into separate tasks. We prove first that
\begin{equation} \label{firsttask}
 \sum_{(p,q)\in \tau_{\bf s}(\blambda)} \frac{1}{\ell}\res(p,q) \doteq   \sum_{r=2}^{\ell } |{}^t\lambda^{(r)}|(\ell s_{r} - \ell s_ 1 + r-1) + \ell \left( \frac{n(n-1)}{2} + \sum_{r=1}^{\ell} n({}^t\lambda^{(r)}) - n(\lambda^{(r)}) \right)
\end{equation}
and then that
\begin{equation}
\label{secondtask}
\sum_{(p,q)\in \tau_{\bf s}(\blambda)} \ep_{\res (p,q)} \res(p,q) \doteq \ell \sum_{r=2}^{\ell} |{}^t\lambda^{(r)}|(\ep_1+\cdots + \ep_{r-1}).
\end{equation}
Using \ref{cdef}, this will complete the calculation and hence the proof of Theorem \ref{bij}.

\subsection{Confirming \eqref{firsttask}} Suppose that we have a charge $s$ and a partition $\lambda$, so that the corresponding $\beta$-numbers are $x_j = \lambda_j + s + 1 - j$ for $j\geq 1$. Let $T> |\lambda|$ be a positive integer. Note that for all $t\geq T$ we have that $x_t = s+ 1-t$. We have $$ n(\lambda) = \sum_{j=1}^T (j-1) \lambda_j  =  \sum_{j=1}^T (j-1) ( x_j + j - 1 - s)$$ and $$ n({}^t\lambda) = \sum_{j=1}^T \frac{1}{2} \lambda_j (\lambda_j - 1)=  \sum_{j=1}^T \frac{1}{2} (x_j + j -1 -s) (x_j + j -2 -s).$$ Now combining these we find \begin{eqnarray*} n(\lambda) - n({}^t\lambda) & = & \frac{1}{2} \sum_{j=1}^T x_j(2s +1 -x_j) - \frac{1}{2} \sum_{j=1}^T (s+ 1-j)(s+j) .\end{eqnarray*} Of course, this expression is independent of the choice of large $T$. It is also independent of the ordering of the $x_j$: this will allow us to calculate both sides of \eqref{firsttask} knowing only the $\beta$-numbers, and not the order they come in.

\subsection{}
In order to remove the ambiguity in \eqref{firsttask} concerning equality up to scalar we will normalise by using the special partition $\bmu = (\emptyset, \ldots , \emptyset, (1^n))\in \mathcal{P}(\ell ,n)$. We will prove that when we subtract terms corresponding to $\bmu$ from the left and right hand sides we get equality.

The $\beta$-numbers corresponding to $\bmu$ are thus $$y_j^{(i)} = \begin{cases} \ell (s_{i}  - j) + i \quad &1\leq i \leq \ell - 1 \\ \ell (s_{\ell} + 1 - j) + \ell  & i = \ell , 1\leq j \leq n \\ \ell (s_{\ell}  - j) + \ell  & i = \ell , j > n. \end{cases} $$

\subsection{}  Following the construction in \ref{multpart} let $A_j^{(i)} = \ell(\lambda_j^{(i)} + s_i - j) + i$ for $1\leq i \leq \ell$, and set $B_j^{(i)} = \ell (s_{i} - j) + i$. Let's go!
\begin{eqnarray*}\sum_{(p,q)\in \tau_{\bf s}(\blambda)} \!\!\!\!\!\res(p,q) -  \!\!\!\!\! \sum_{(p,q)\in \tau_{\bf s}(\bmu)} \!\!\!\!\! \res(p,q)& = &  (n( {}^t\tau_{\bf s}(\blambda)) - n(\tau_{\bf s}(\blambda))) - (n({}^t\tau_{\bf s}(\bmu)) - n(\tau_{\bf s}(\bmu)))  \\  & = &  \frac{1}{2} \sum_{i=1}^{\ell}
 \sum_{j=1}^T \left[ A_j^{(i)}(A_j^{(i)}-1 )
- y_j^{(i)}(y_j^{(i)}-1)\right] \\ & = & \frac{1}{2} \sum_{i=1}^{\ell}
 \sum_{j=1}^T \left[ A_j^{(i)}(A_j^{(i)}-1) -  (\ell( s_{i} - j) + i)(\ell (s_{i} - j) + i-1 )\right] + \\ && + \frac{1}{2} \sum_{j=1}^n \ell \left( 1 - 3\ell - 2\ell (s_{\ell}  -j) \right) . \end{eqnarray*}
Thus the term we would like calculate equals
\begin{eqnarray}\label{yougotit} \frac{1}{2\ell} \sum_{i=1}^{\ell}
 \sum_{j=1}^T \left[ A_j^{(i)}(A_j^{(i)}-1) - B_j^{(i)}(B_j^{(i)}-1)\right] + \frac{1}{2}n(1- 2\ell s_{\ell} - 2\ell + \ell n). \end{eqnarray}

\subsection{} Now we calculate the difference between the $c$-function of ${}^t\blambda$ and of ${}^t\bmu$. We get \begin{eqnarray} \label{wewantit} \qquad \sum_{i=2}^{\ell} |{}^t \lambda^{(i)} |(\ell s_{i} - \ell s_1 + i-1) - \sum_{i= 1}^{\ell } \ell (n(\lambda^{(i)}) - n({}^t\lambda^{(i)})) - n (\ell s_{\ell} - \ell s_1 + \ell -1) + \frac{1}{2} \ell n(n-1). \end{eqnarray}
Obviously $$\lambda_j^{(i)} = \frac{1}{\ell} ( A_j^{(i)} - i) +j - s_{i}  = \frac{1}{\ell} ( A_j^{(i)} - B_j^{(i)})$$ and thus $$ \ell n = \sum_{i=1}^{\ell} \ell |\lambda^{(i)}| =  \sum_{i=1}^{\ell } \sum_{j=1}^T  ( A_j^{(i)} - B_j^{(i)}).$$  We now have \begin{eqnarray*} \eqref{wewantit} & = &
  \sum_{i=1}^{\ell}\left[ | {}^t\lambda^{(i)} |(\ell s_{i} - \ell s_1 + i-1) - \ell (n(\lambda^{(i)}) - n({}^t\lambda^{(i)})) \right] - n (\ell s_{\ell} - \ell s_1 + \ell -1) + \frac{1}{2} \ell n(n-1) \\
 & = & \sum_{i=1}^{\ell} \sum_{j=1}^T  \left[ \lambda_j^{(i)} ( \ell s_{i} - \ell s_1 + i-1) - \ell (j-1)\lambda_j^{(i)} + \frac{\ell}{2} \lambda_j^{(i)}(\lambda_j^{(i)}-1)\right] - n (\ell s_{\ell} - \ell s_1 + \ell -1) + \frac{1}{2} \ell n(n-1) \\
 & = & \sum_{i=1}^{\ell} \sum_{j=1}^T  \left[ \lambda_j^{(i)} ( \ell s_{i}  + i-1) - \ell (j-1)\lambda_j^{(i)} + \frac{\ell}{2} \lambda_j^{(i)}(\lambda_j^{(i)}-1)\right] - n (\ell s_{\ell}  +\ell -1) + \frac{1}{2} \ell n(n-1) \\
 & = & \sum_{i=1}^{\ell} \sum_{j=1}^T  \frac{1}{2} \lambda_j^{(i)}\left[ 2\ell s_{i}  + 2i - 2- 2\ell j + \ell + \ell \lambda_j^{(i)}\right] + \frac{n}{2} (\ell n - 2\ell s_{\ell}  + 2 - 3\ell )
  \\
  & = & \sum_{i=1}^{\ell} \sum_{j=1}^T  \frac{1}{2} \lambda_j^{(i)}\left[ A_{j}^{(i)} + \ell s_{i}  + i +\ell -2 - \ell j \right]+ \frac{n}{2} (\ell n - 2\ell s_{\ell}  + 2 - 3\ell )
 \\
 & = &  \sum_{i=1}^{\ell} \sum_{j=1}^T  \frac{1}{2} \lambda_j^{(i)}\left[ A_{j}^{(i)} + B_j^{(i)}\right] + \frac{n}{2}(\ell -2) + \frac{n}{2} (\ell n - 2\ell s_{\ell}  + 2 - 3\ell )
\\
& = & \sum_{i=1}^{\ell} \sum_{j=1}^T \frac{1}{2\ell}\left[ ( A_j^{(i)} - B_j^{(i)})(A_j^{(i)}+ B_j^{(i)} )\right] + \frac{n}{2} (\ell n - 2\ell s_{\ell} -2 \ell ) \\
& = & \sum_{i=1}^{\ell } \sum_{j=1}^T  \frac{1}{2\ell}\left[ A_j^{(i)}(A_j^{(i)} -1)  - B_j^{(i)}(B_j^{(i)} - 1) + (A_j^{(i)} - B_j^{(i)} )\right] + \frac{n}{2} (\ell n - 2\ell s_{\ell} -2 \ell )\\
&=& \sum_{i=1}^{\ell} \sum_{j=1}^T  \frac{1}{2\ell}\left[ A_j^{(i)}(A_j^{(i)} -1)  - B_j^{(i)}(B_j^{(i)} - 1) \right] +\frac{n}{2} + \frac{n}{2} (\ell n - 2\ell s_{\ell} -2 \ell )\\ &= & \eqref{yougotit}.
\end{eqnarray*}

\subsection{Confirming \eqref{secondtask}} Since $\bep \in \Q_0^{\ell}$ we have $$ \sum_{(p,q)\in \tau_{\bf s}(\blambda)} \ep_{\res (p,q)} \res(p,q) = \sum_{i=1}^{\ell -1} \ep_i \left(  \sum_{(p,q)\in \tau_{\bf s}(\blambda)\atop \res (p,q) \equiv i}\!\!\!\!\! \res(p,q)  - \!\!\!\!\! \sum_{(p,q)\in \tau_{\bf s}(\blambda) \atop \res (p,q)\equiv 0}  \!\!\!\!\! \res(p,q)\right) . $$ Therefore we only need to show for $1\leq i \leq \ell -1$ that \begin{equation} \label{what} \sum_{(p,q)\in \tau_{\bf s}(\blambda)\atop \res (p,q) \equiv i} \!\!\!\!\!\res(p,q)  -  \!\!\!\!\!\sum_{(p,q)\in \tau_{\bf s}(\blambda) \atop \res (p,q)\equiv 0} \!\!\!\!\!\res(p,q) \doteq {\ell} \sum_{r=i+1}^{\ell } |\lambda^{(r)}|.\end{equation}

\subsection{} By \cite[Proof of Proposition 9.3]{AKT} we have for any $i\in \Z$ the following formula for $\cont_i(\tau_{\bf s}(\blambda))$  $$\cont_i(\tau_{\bf s}(\blambda))  = \sum_{j \in \Z} \left[ \#\{ y\leq \ell j + i : y\notin \beta_0(\tau_{\bf s}(\blambda))\} - S(\ell j+ i)\right ].$$ Here for any integer $k$ we define $S(k) = k$ if $k\geq 0$ and $S(k) = 0$ otherwise. Similarly, the sum of the residues congruent to $i$ modulo $\ell$ is given by $$\sum_{j \in \Z} (\ell j+i)\left[ \#\{ y\leq \ell j + i : y\notin \beta_0(\tau_{\bf s}(\blambda))\} - S(\ell j+ i)\right ].$$ This means that we can calculate \eqref{what} from knowledge of the $\beta$-numbers only.

\subsection{} We have \begin{eqnarray*}  \sum_{(p,q)\in \tau_{\bf s}(\blambda)\atop \res (p,q) \equiv i} \!\!\!\!\!\res(p,q)  -  \!\!\!\!\!\sum_{(p,q)\in \tau_{\bf s}(\blambda) \atop \res (p,q)\equiv 0} \!\!\!\!\!\res(p,q) & = & \sum_{j \in \Z} (\ell j+i)\left[ \#\{ y\leq \ell j + i : y\notin \beta_0(\tau_{\bf s}(\blambda))\} - S(\ell j+ i)\right ] \\ & & - \sum_{j \in \Z} \ell j\left[ \#\{ y\leq \ell j  : y\notin \beta_0(\tau_{\bf s}(\blambda))\} - S(\ell j)\right ] \\ & = & \sum_{j\in \Z} \ell j \left[ \#\{\ell j< y \leq \ell j + i : y\notin \beta_0(\tau_{\bf s}(\blambda)) \} - S(\ell j +i) +S(\ell j) \right] \\ && + \sum_{j\in \Z} i \left[ \# \{ y \leq \ell j + i : y\notin \beta_0(\tau_{\bf s}(\blambda) \} - S(\ell j + i) \right] . \end{eqnarray*} The second term here is just $i\cont_i(\tau_{\bf s}(\blambda))$, which by is independent of $\blambda$. So we now need $$ \sum_{j\in \Z} \ell j \left[ \#\{\ell j< y \leq \ell j + i : y\notin \beta_0(\tau_{\bf s}(\blambda)) \} - S(\ell j +i) +S(\ell j) \right] \doteq  \ell (|\lambda^{(i+1)}| + \cdots + |\lambda^{(\ell)}|). $$
Define $T(j) = 1$ if $j\geq 0$ and $T(j) = 0$ otherwise and let $U<<0$. Now we apply repeatedly the following formula with $1\leq k \leq i$ \begin{eqnarray*} \sum_{j\in \Z} j \left[ \# \{ y = j  : y=\ell j + k\notin \beta_{{0}}(\tau_{\bf s}(\blambda)) \} - T(j) \right] & = & \sum_{j<0} j \left[ \# \{ y = \ell j + k  : y\notin \beta_{0}(\tau_{\bf s}(\blambda)) \}\right] + \\ && + \sum_{j\geq 0} j\left[ \# \{ y = \ell j + k  : y\notin \beta_{0}(\tau_{\bf s}(\blambda)) \}-1\right]  \\ & = & \sum_{U\leq j<0} j - \sum_{U\leq j} j\left[ \# \{ y = \ell j + k : y \in \beta_{0}(\tau_{\bf s}(\blambda)) \} \right] \\ & \doteq & -\!\!\!\!\!\sum_{r\atop U\leq \lambda^{(k)} + s_k - r } \!\!\!\!\!(\lambda^{(k)} + s_k -r) \\ & \doteq & -|\lambda^{(k)}|. \end{eqnarray*}
This gives equality with $-|\lambda^{(1)}| - \cdots - |\lambda^{(i)}|$ and hence with $|\lambda^{(i+1)}|+\cdots + |\lambda^{(\ell)}| - \ell n$. This gives \eqref{what} and hence \eqref{secondtask}.

\section{Calculation for proof of Theorem \ref{afnrel}}
We will keep the same notation as in the proof.
\subsection{} We have to show that \eqref{wanta} is true, i.e. that \begin{eqnarray*}&& - \sum_{1\leq i, j \leq \ell \atop {1\leq u\leq n \atop 1\leq k \leq \lambda^{(i)}_u}} \!\!\!\!\!\!\min\{n+k-u+ (H_1 + \cdots H_{i-1}), (H_1 + \cdots H_{j-1})\} \\ &&\qquad \qquad \doteq \sum_{1\leq i \neq j \leq \ell \atop {1\leq u \leq n \atop n+1\leq v \leq n+s_j - S }}\min \{ B_u^{(i)}(\blambda), B_v^{(j)}(\blambda) \} + \sum_{1\leq i \leq \ell  \atop {1\leq u \leq n \atop n+1\leq v \leq n+s_i - S}} \min \{ B_u^{(i)}(\blambda), B_v^{(i)}(\blambda) \}. \end{eqnarray*} Here $B_u^{(i)}(\blambda) = n+\lambda_u^{(i)}-u + H_1 + \cdots + H_{i-1}$ and ... For ease of notation let $M^{(i)} = H_1 + \cdots + H_{i-1}$

We will call the left hand side of this equality $L(\blambda)$ and the right hand side $R(\blambda)$. We carry out the comparison of $L(\blambda)$ and $R(\blambda)$ in two steps. First let $\bmu \in \mathcal{P}(\ell,n)$ be chosen such that $\mu^{(j)} = \lambda^{(j)}$ for all $j \neq i$, that $\mu^{(i)}_u = \lambda^{(i)}_u$ for all $u\neq r, r+t$, and that $\mu^{(i)}_r = \lambda^{(i)}_r -1$ and  $\mu^{(i)}_{r+t} = \lambda^{(i)}_{r+t}+1$ where $t>0$. We will show that $L(\blambda) - L(\bmu) =  R(\blambda) - R(\bmu)$. Second let $\btau, \bnu \in \mathcal{P}(\ell, n)$ be chosen such that $\btau = ((\tau_1), \ldots , (\tau_{\ell}))$ and $\bnu = ((\nu_1), \ldots , (\nu_\ell))$ where $\sum \tau_i = \sum\nu_i = n$ and $\tau_j = \nu_j$ for $j\neq i,i+1$ whilst $\tau_i = \nu_i -1$ and $\tau_{i+1} = \nu_{i+1}+1$. We will show that  $L(\bnu) - L(\btau) =  R(\bnu) - R(\btau)$. Since every $\ell$-multipartition of $n$ can be obtained from $((n), \emptyset, \ldots, \emptyset)$ by a sequence of the two moves above it will follow that $$L(\blambda) \doteq R(\blambda) \qquad \text{for all }\blambda \in \mathcal{P}(\ell ,n).$$

\subsection{} Let $\bmu$ be as above. We have \begin{eqnarray*} L(\blambda) - L(\bmu) &=&  \sum_{j=1}^{\ell} \min\{n + \lambda_{r+t}^{(i)}+1 -r-t + M^{(i)}, M^{(j)}\}  - \min\{n+ \lambda_{r}^{(i)} - r + M^{(i)}, M^{(j)}\} \\ & = & \sum_{j=1}^{\ell} \min\{B_{r+t}^{(i)}(\blambda) + 1, M^{(j)}\} - \min\{B_{r}^{(i)}(\blambda), M^{(j)}\}. \end{eqnarray*} On the other hand \begin{eqnarray*} R(\blambda) - R(\bmu) = \sum_{1\leq j \leq \ell \atop n+1 \leq v \leq n+s_j - S} && \min\{B_r^{(i)}(\blambda), B_v^{(j)}(\blambda) \} -  \min\{B_r^{(i)}(\blambda)-1, B_v^{(j)}(\blambda) \} \\ &&
 +  \min\{B_{r+t}^{(i)}(\blambda), B_v^{(j)}(\blambda) \}
 -  \min\{B_{r+t}^{(i)}(\blambda)+1, B_v^{(j)}(\blambda) \}. \end{eqnarray*}

Fix $j$. For any $\bnu \in \mathcal{P}(\ell ,n)$ we have $B^{(j)}_{v}(\bnu) < M^{(j)} \leq B^{(j)}_u(\bnu)$ for $1\leq j \leq \ell$, $1\leq u \leq n$ and $v\geq n+1$. Thus there are three cases to consider.

\subsection{Case 1: $M^{(j)}\geq B_r^{(i)}(\blambda)$} We have $$\min\{B_{r+t}^{(i)}(\blambda) + 1, M^{(j)}\} - \min\{B_{r}^{(i)}(\blambda), M^{(j)}\} = B_{r+t}^{(i)}(\blambda) + 1 - B_r^{(i)}(\blambda).$$ Let $x$ and $y$ be defined by the following inequalities: $$B_{x+1}^{(j)} < B_{r}^{(i)}\leq B_x^{(j)} \qquad \text{and} \qquad B_{y+1}^{(j)} < B_{r+t}^{(i)}+1 \leq B_y^{(j)}.$$ Then we have
\begin{eqnarray*}  &&\!\!\!\!\!\!\sum_{n+1\leq v\leq n+s_j - S} \!\!\!\!\!\!\!\min\{B_r^{(i)}(\blambda), B_v^{(j)}(\blambda) \} -  \min\{B_r^{(i)}(\blambda)-1, B_v^{(j)}(\blambda) \} \\ && \qquad \qquad
 +  \min\{B_{r+t}^{(i)}(\blambda), B_v^{(j)}(\blambda) \}
 -  \min\{B_{r+t}^{(i)}(\blambda)+1, B_v^{(j)}(\blambda) \} \\ & = & (x-n) + (B_{x+1}^{(j)}(\blambda) - (B_r^{(i)}(\blambda) - 1)) - (y-n) + (B_{r+t}^{(i)}(\blambda)-B_{y+1}^{(j)}(\blambda) ) \\ & = & (x-y) + (B_{x+1}^{(j)}(\blambda) - B_{y+1}^{(j)}(\blambda)) + B_{r+t}^{(i)}(\blambda) + 1 - B_r^{(i)}(\blambda) \\ & = & B_{r+1}^{(i)}(\blambda) + 1 - B_r^{(i)}(\blambda) =\min\{B_{r+t}^{(i)}(\blambda) + 1, M^{(j)}\} - \min\{B_{r}^{(i)}(\blambda), M^{(j)}\} . \end{eqnarray*}
\subsection{Case 2: $B_r^{(i)} (\blambda) > M^{(j)} \geq B_{r+t}^{(i)}(\blambda) + 1$}We have $$\min\{B_{r+t}^{(i)}(\blambda) + 1, M^{(j)}\} - \min\{B_{r}^{(i)}(\blambda), M^{(j)}\} = B_{r+t}^{(i)}(\blambda) + 1 - M^{(j)}.$$ Note that $B_{r}^{(i)}(\blambda) -1 > B_{n+1}^{(j)}(\blambda)$ since $B_r^{(i)} (\blambda) > M^{(j)}$. Let $y$ be defined by the inequality $ B_{y+1}^{(j)} < B_{r+t}^{(i)}+1 \leq B_y^{(j)}.$ Then we have
\begin{eqnarray*}  &&\!\!\!\!\!\!\sum_{n+1\leq v\leq n+s_j - S} \!\!\!\!\!\!\!\min\{B_r^{(i)}(\blambda), B_v^{(j)}(\blambda) \} -  \min\{B_r^{(i)}(\blambda)-1, B_v^{(j)}(\blambda) \} \\ && \qquad \qquad
 +  \min\{B_{r+t}^{(i)}(\blambda), B_v^{(j)}(\blambda) \}
 -  \min\{B_{r+t}^{(i)}(\blambda)+1, B_v^{(j)}(\blambda) \} \\ & = & 0 - (y-n) + (B_{r+t}^{(i)}(\blambda)-B_{y+1}^{(j)}(\blambda) ) \\ & = & n-y - (n-(y+1) + M^{(j)}) +B_{r+t}^{(i)}(\blambda)  \\ & = & B_{r+t}^{(i)}(\blambda) + 1 - M^{(j)} =\min\{B_{r+t}^{(i)}(\blambda) + 1, M^{(j)}\} - \min\{B_{r}^{(i)}(\blambda), M^{(j)}\}  \end{eqnarray*}

\subsection{Case 3: $B_{r+t}^{(i)}(\blambda) + 1 > M^{(j)}$}We have $$\min\{B_{r+t}^{(i)}(\blambda) + 1, M^{(j)}\} - \min\{B_{r}^{(i)}(\blambda), M^{(j)}\} = 0.$$ Since $B_{r+t}^{(i)}(\blambda) > B_{n+1}^{(j)} (\blambda)$ by our assumption. Thus
\begin{eqnarray*}  &&\!\!\!\!\!\!\sum_{n+1\leq v\leq n+s_j - S} \!\!\!\!\!\!\!\min\{B_r^{(i)}(\blambda), B_v^{(j)}(\blambda) \} -  \min\{B_r^{(i)}(\blambda)-1, B_v^{(j)}(\blambda) \} \\ && \qquad \qquad
 +  \min\{B_{r+t}^{(i)}(\blambda), B_v^{(j)}(\blambda) \}
 -  \min\{B_{r+t}^{(i)}(\blambda)+1, B_v^{(j)}(\blambda) \} \\ & = & 0. \end{eqnarray*}

\subsection{} This ends the analysis of the three possible cases and proves that $L(\blambda) - L(\bmu) = R(\blambda) - R(\mu)$. The proof of the equality $L(\bnu) - L(\btau) = R(\bnu) - R(\btau)$ is very similar and involves no new ideas, so we leave it to the reader.


 
\section*{Index of Notation}\label{index}  
\begin{multicols}{2}
{\small  \baselineskip 14pt 
$\lhd$, dominance order\hfill\eqref{dom-defn}

$<_{\bf h}$, the $c$-order \hfill\eqref{cord-defn}

$\prec_{\bf h}$, the geometric order\hfill\eqref{geomorder}

$\lhd_J$, refinement of $\lhd$ depending on $J$\hfill\eqref{Jord-defn}

$\ast$, shifted action of $\tilde{\mathfrak{S}}_{\ell}$\hfill\eqref{ast-defn}

$a_{\bf h}(\blambda)$, the $a$-function\hfill\eqref{a-defn}

$\beta_s(\lambda)$, $\beta$ numbers of $\lambda$\hfill\eqref{beta-defn}

$c_{\bf h}(\blambda)$, the $c$-function \hfill\eqref{cdef}

$\res(p,q)$, content of node at $(p,q)$\hfill\eqref{cont-defn}

$G = G_n(\ell)$, the group $\mathfrak{S}_n\ltimes (\mu_{\ell})^n$ \hfill\eqref{Gnl-defn} 

${\bf h}$, parameters $(h, H_1, \ldots , H_{\ell-1})\in {\bf H}$ \hfill\eqref{h-defn}

${\bf H}$, parameter space for $H_{1, {\bf h}}$ \hfill\eqref{H-defn}

$\hil{\nu}$, component of invariant Hilbert scheme\hfill\eqref{hiv-defn}

$\Hr$, G.I.T. chambers\hfill\eqref{Hr-defn}

$H_{t, {\bf h}}$, the rational Cherednik algebra\hfill\eqref{RCA-defn}

$I_{\lambda}$, monomial ideal of $\C[A,B]$\hfill\eqref{I-defn}

$J$-class, partitions with same $J$-heart\hfill\eqref{Jclass-defn}

$J$-heart, removal of $j$-nodes\hfill\eqref{heart-defn}

$\mathcal{M}_{\btheta}({\bf d}), \mathcal{M}_{\btheta}(n)$, G.I.T. quotient of $R({\bf d'})$\hfill\eqref{M-defn},\eqref{M-2defn}

$\cont_i(\lambda)$, count of content of $\lambda$\hfill\eqref{N-defn}

$\mathcal{O}_{\bf h}$, category $\mathcal{O}$ for $H_{1,{\bf h}}$\hfill\eqref{O}

$\pi_{\btheta}$, (partial) resolution of $V/G$ \hfill\eqref{pi-defn}

$\mathcal{P}(\ell, n)$,  $\ell$-multipartitions of $n$ \hfill\eqref{multpart-defn}

$\mathcal{P}(n)$, partitions of $n$\hfill\eqref{part-defn}

$\mathcal{P}_{\nu}(n)$, partitions of $n$ with core $\nu$\hfill\eqref{partnu-defn}

$Q, \overline{Q}$, cyclic quiver and its double\hfill\eqref{Q-defn}

$Q_{\infty}$, $\overline{Q}_{\infty}$, extended cyclic quiver and its double\hfill\eqref{Qi-defn}

$R({\bf d'})$, ${\bf d'}$-dimensional $\overline{Q}_{\infty}$-representation space \hfill\eqref{rep-defn}

$\tilde{\mathfrak{S}}_{\ell}$, affine symmetric group\hfill\eqref{affine-defn}

$\tau_{\bf s}$, bijection between $\mathcal{P}(\ell,n)$ and $\mathcal{P}(n)$\hfill\eqref{tau-defn}

$\Theta_1$, affine hyperplane of stability parameters\hfill\eqref{slice-defn}

${}^t\lambda, \, {}^t\blambda$, transpose of a (multi)partition\hfill\eqref{trpart-defn}

type $J$, the parameter stabiliser\hfill\eqref{type-defn}
 
$\mathcal{X}_{\btheta}({\bf d}), \mathcal{X}_{\btheta}(n)$, algebraic quotient of $R({\bf d'})$ \hfill\eqref{X-defn}, \eqref{X-2defn}

$x_{\btheta}(\blambda)$, fixed point of $\mathcal{M}_{\btheta}(n)$\hfill\eqref{xt-defn}

$\mathcal{Z}_{\blambda}, \mathcal{Z}_{\btheta},$ (components of) attracting subvariety\hfill\eqref{Zt-defn}

$Z_{0, {\bf h}}$, the centre of $H_{0,{\bf h}}$\hfill\eqref{Z-defn}

}
\end{multicols}


\end{document}